\documentclass{alggeom}
\usepackage[utf8]{inputenc}

\usepackage{amsmath}
\usepackage{amssymb}
\usepackage{amsfonts}
\usepackage{mathrsfs}
\usepackage{xcolor}
\usepackage{dsfont}
\usepackage{stmaryrd}
\usepackage{bm}
\usepackage{hyperref}
\hypersetup{
     colorlinks=true,
     linkcolor=purple,
     filecolor=purple,
     citecolor=purple,      
     urlcolor=blue,
     }
\usepackage{tikz-cd}
\usepackage{tikz}
\usepackage{tikzsymbols}
\usepackage{physics}
\usetikzlibrary{matrix}

\usepackage[nosubsections]{mytheorems2}

\newcommand{\q}{\; / \;}

\DeclareMathOperator{\codim}{\operatorname{codim}}
\DeclareMathOperator{\rk}{\operatorname{rank}}
\DeclareMathOperator{\characteristic}{\operatorname{char}}
\DeclareMathOperator{\lcm}{\operatorname{lcm}}
\newcommand{\gp}{\textup{gp}}

\newcommand{\NN}{\mathbb{N}} 
\newcommand{\kk}{\mathds{k}} 
\newcommand{\ZZ}{\mathbb{Z}} 
\newcommand{\QQ}{\mathbb{Q}} 
\newcommand{\GG}{\mathbb{G}} %

\newcommand{\fs}{\mathfrak{s}}

\newcommand{\SE}{\mathscr{E}}

\newcommand{\SO}{\mathscr{O}}

\newcommand{\SM}{\mathscr{M}}

\newcommand{\SD}{\mathscr{D}}

\newcommand{\SR}{\mathscr{R}}

\newcommand{\SJ}{\mathscr{J}}
\newcommand{\SQ}{\mathscr{Q}}

\DeclareMathOperator{\Spec}{\operatorname{Spec}}
\DeclareMathOperator{\GSpec}{\underline{\Spec}}
\DeclareMathOperator{\ZR}{\operatorname{ZR}}
\DeclareMathOperator{\val}{\operatorname{val}}
\DeclareMathOperator{\Supp}{\operatorname{Supp}}

\DeclareMathOperator{\ord}{\operatorname{ord}}
\DeclareMathOperator{\logord}{\operatorname{log-ord}}
\DeclareMathOperator{\Hom}{\operatorname{Hom}}
\DeclareMathOperator{\Proj}{\operatorname{Proj}} 
\DeclareMathOperator{\SProj}{\underline{\mathscr{P}\textup{roj}}} 
\DeclareMathOperator{\Bl}{\operatorname{Bl}} 

\newcommand{\Gm}{\GG_m}
\newcommand{\Gmu}{\pmb{\mu}}

\makeatletter
\newcommand*{\relrelbarsep}{.386ex}
\newcommand*{\relrelbar}{%
  \mathrel{%
    \mathpalette\@relrelbar\relrelbarsep
  }%
}
\newcommand*{\@relrelbar}[2]{%
  \raise#2\hbox to 0pt{$\m@th#1\relbar$\hss}%
  \lower#2\hbox{$\m@th#1\relbar$}%
}
\providecommand*{\rightrightarrowsfill@}{%
  \arrowfill@\relrelbar\relrelbar\rightrightarrows
}
\providecommand*{\leftleftarrowsfill@}{%
  \arrowfill@\leftleftarrows\relrelbar\relrelbar
}
\providecommand*{\xrightrightarrows}[2][]{%
  \ext@arrow 0359\rightrightarrowsfill@{#1}{#2}%
}
\providecommand*{\xleftleftarrows}[2][]{%
  \ext@arrow 3095\leftleftarrowsfill@{#1}{#2}%
}
\makeatother

\DeclareMathOperator{\inv}{\operatorname{inv}}
\DeclareMathOperator{\MC}{\operatorname{MC}}

\begin{document}

\title{Logarithmic Resolution via Weighted Toroidal Blow-ups}
\author{Ming Hao Quek}
\email{ming\_hao\_quek@brown.edu}
\address{Department of Mathematics, Brown University, Box 1917, 151 Thayer Street, Providence, RI 02906}

\thanks{This research work was supported in part by funds from BSF grant 2014365 and NSF grant DMS-1759514.}

\begin{abstract}
Let $X$ be a fs logarithmic scheme that is generically logarithmically smooth, and that admits a strict closed embedding into a logarithmically smooth scheme $Y$ over a field $\kk$ of characteristic zero. We construct a simple and fast procedure to functorial logarithmic resolution of $X$, where the end result is in particular a stack-theoretic modification $X' \rightarrow X$ such that $X'$ is logarithmically smooth over $\kk$. In particular, if $X$ is a closed subscheme of a smooth $\kk$-scheme $Y$, the procedure not only shares the same desirable features as the ``dream resolution algorithm'' of Abramovich-Temkin-W{\l}odarczyk \cite{ATW19}, but also accounts for a key feature of Hironaka's \cite[Main Theorem I]{Hir64} which was not addressed in \cite{ATW19}. As a consequence, we recover a different and simpler approach to Hironaka's resolution of singularities in characteristic zero.
\end{abstract}

\maketitle

\section{Introduction}\label{chpt1}

\subsection{Statement of main theorem}\label{1.1}

Consider a fs (\ref{B101}) logarithmic scheme $Y$ which is logarithmically smooth over a field $\kk$ of characteristic zero, or equivalently a \emph{toroidal} $\kk$-scheme $Y$ (\ref{B106}), as well as a reduced closed subscheme $X \subset Y$ of pure codimension $c$. More generally, we consider a reduced closed substack $X$ of pure codimension $c$ in a \emph{toroidal Deligne-Mumford stack $Y$ over $\kk$} (\ref{B201}). We will always regard $X$ (without mention) as a logarithmic Deligne-Mumford stack over $\kk$ by pulling the logarithmic structure of $Y$ back to $X$. Such pairs $X \subset Y$ form the objects of a category, where a morphism between pairs $(\widetilde{X} \subset \widetilde{Y}) \rightarrow (X \subset Y)$ is a cartesian square \begin{equation*}
    \begin{tikzcd}
    \widetilde{X} = X \times_Y \widetilde{Y} \arrow[to=1-2, hookrightarrow] \arrow[to=2-1] & \widetilde{Y} \arrow[to=2-2, "f"] \\
    X \arrow[to=2-2, hookrightarrow] & Y
    \end{tikzcd}
\end{equation*} 
where $f \colon \widetilde{Y} \rightarrow Y$ is logarithmically smooth and surjective. We refer to such morphisms as logarithmically smooth, surjective morphisms of pairs. Note, however, that in certain situations below, we do not demand surjectivity in our morphisms of pairs. 

The goal of this paper is to define a \emph{logarithmic embedded resolution functor} on the aforementioned category, which assigns to each pair $X \subset Y$ as above, a proper birational morphism $\Pi \colon Y' \rightarrow Y$ such that both $Y'$ and the proper transform $X' \subset X \times_Y Y'$ are toroidal. Moreover, $\Pi$ satisfies two properties which resemble those in Hironaka's \cite[Main Theorem I]{Hir64}: \begin{enumerate}
    \item $\Pi$ is an isomorphism over the logarithmically smooth locus $X^{\textup{log-sm}}$ of $X$.
    \item We are able to control $\Pi^{-1}(X \setminus X^{\textup{log-sm}})$: namely, $\Pi^{-1}(X \setminus X^{\textup{log-sm}})$ will be contained in the toroidal divisor (\ref{B105}(ii)) of $X'$. 
\end{enumerate}

We will explicitly construct the proper birational morphism $\Pi$ as a composition of stack-theoretic blow-ups along \emph{toroidal centers} (which are the \emph{weighted toroidal blow-ups} in the title of this paper). The notion of a \emph{toroidal center} will be defined in \S\ref{3.2}. A convenient tool for bookkeeping the information carried by toroidal centers is the notion of \emph{idealistic exponents}, which we study in detail in \S\ref{chpt2}. In \S\ref{4.4}, we will also explicate the charts of the weighted toroidal blow-ups appearing in $\Pi$.

In addition, for a point $p \in \abs{Y}$, \S\ref{6.1} defines a (logarithmic) \emph{invariant} of $X \subset Y$ at $p$ (motivated by the invariants in \cite{ATW20a} and \cite{ATW19}), denoted by $\inv_p(X \subset Y)$, which is an nondecreasing truncated sequence of nonnegative rational numbers, but we allow for the last entry to be $\infty$. There is a total order on the set consisting of all such invariants of $X \subset Y$ at a point $p$ by the lexicographic order $<$, which turns out to be a well-ordering (\ref{6103}(i)). There is a caveat here: our lexicographic order considers truncated sequences to be strictly larger. For example, \begin{align*}
    (0) < (1,1,2,2) < (1,1,3) < (1,1) < (1,2,5) < (1,3,4) < (1,\infty) < (1) < (),
\end{align*}
where $()$ is the empty sequence. The invariant satisfies the following properties: \begin{enumerate}
    \item[(a)] Assuming $c \geq 1$, it detects logarithmic smoothness at any $p \in \abs{X}$: $\inv_p(X \subset Y)$ is bounded below (via the lexicographic order) by the sequence $(1,\dotsc,1)$ of length $c$, and equality holds if and only if $X$ is logarithmically smooth at $p$.
    \item[(a$'$)] $\inv_p(X \subset Y) = (0)$ if and only if $p \notin \abs{X}$. 
    \item[(a$''$)] If $c=0$ (i.e. $X = Y$), $\inv_p(X \subset Y) = ()$ for all $p \in \abs{Y} = \abs{X}$.
    \item[(b)] It is upper semi-continuous on $Y$.
    \item[(c)] It is functorial for logarithmically smooth morphisms of pairs $X \subset Y$, whether or not surjective. See (\ref{6103}(iii)) for a precise statement.
    \item[(d)] The first term of $\inv_p(X \subset Y)$ is the logarithmic order (\ref{B114}) at $p$ of the ideal $\mathcal{I}_{X \subset Y}$ of $X$ embedded in $Y$. In particular, it is in $\NN \cup \lbrace \infty \rbrace$.
\end{enumerate}
This invariant is constructed via logarithmic analogues of the classical notions of maximal contact elements and coefficient ideals, which we study in \S\ref{chpt5}. We set the \emph{maximal invariant} of $X \subset Y$ to be $\max\inv(X \subset Y) = \max_{p \in \abs{X}}{\inv_p(X \subset Y)}$: this is functorial for logarithmically smooth and surjective morphisms of pairs $X \subset Y$, and is equal to the sequence $(1,\dotsc,1)$ of length $c$ if and only if $X$ is toroidal. We can now state the main result:

\begin{theorem}\label{1101}
There is a functor $F_{\textup{log-ER}}$ associating to \begin{quote}
    a reduced, closed substack $X$ of pure codimension in a toroidal Deligne-Mumford stack $Y$ over a field $\kk$ of characteristic zero, such that $X$ is generically toroidal but not toroidal
\end{quote}
a toroidal center $\overline{\SJ}$ on $Y$, with weighted toroidal blow-up $\pi \colon Y' \rightarrow Y$ along $\overline{\SJ}$, and proper transform $F_{\textup{log-ER}}(X \subset Y) = X' \subset Y'$, such that: \begin{enumerate}
    \item $Y'$ is again a toroidal Deligne-Mumford stack over $\kk$;
    \item $\max\inv(X' \subset Y') < \max\inv(X \subset Y)$;
    \item $\pi$ is an isomorphism away from the closed locus consisting of points $p \in \abs{X}$ with $\inv_p(X \subset Y) = \max\inv(X \subset Y)$;
    \item the exceptional divisor underlying $\pi$ is contained in the toroidal divisor of $Y'$.
\end{enumerate}
Functoriality here is with respect to logarithmically smooth, surjective morphisms of pairs $X \subset Y$, as described before the theorem. 

In particular, one stops at an integer $N \geq 1$ where the iterated application $(X_N \subset Y_N) = F_{\textup{log-ER}}^{\circ N}(X \subset Y)$ is accompanied by a sequence of weighted toroidal blow-ups $\Pi \colon Y_N \xrightarrow{\pi_{N-1}} \dotsb \xrightarrow{\pi_1} Y_1 \xrightarrow{\pi_0} Y_0 = Y$ such that: \begin{enumerate}
    \item[(a)] $X_N$ and $Y_N$ are both toroidal Deligne-Mumford stacks over $\kk$;
    \item[(b)] $\Pi$ is an isomorphism over the logarithmically smooth locus $X^{\textup{log-sm}}$ of $X$;
    \item[(c)] $\Pi^{-1}(X \setminus X^{\textup{log-sm}})$ is contained in the toroidal divisor of $X_N$.
\end{enumerate} 
This stabilized functor $F_{\textup{log-ER}}^{\circ \infty}$, together with the sequence of weighted toroidal blow-ups $\Pi$ after removing empty blow-ups, is functorial for all logarithmically smooth morphisms of pairs $X \subset Y$, whether or not surjective.
\end{theorem}

The toroidal center $\overline{\SJ}$ associated to $X \subset Y$ in the first paragraph of (\ref{1101}), will be defined and studied in \S\ref{chpt6}. In the words of \cite{Kol07}, (\ref{1101}) will be proven, seemingly by accident, via the logarithmic analogue of \emph{principalization} (\ref{7101}). This is the content of \S\ref{chpt7}.

Finally, we remark that parts (ii) and (iii) of (\ref{1101}) are precisely the two key features of the ``dream resolution algorithm'' in \cite{ATW19}: namely, \begin{enumerate}
    \item each step of that algorithm improves singularities immediately in a visible way,
    \item and does so by blowing up the most singular locus. 
\end{enumerate}
It is well-known that besides the case of curves, these two features are in general not plausible for Hironaka's resolution algorithm. There are plenty of examples corroborating this observation --- for example, see \cite[1.7]{ATW19}.

\subsection{Recovering Hironaka's resolution of singularities}\label{1.2}

In this section, we recover Hironaka's \cite[Main Theorem I]{Hir64} from (\ref{1101}) in three steps. The first step is to deduce \emph{logarithmic resolution} from (\ref{1101}):

\begin{theorem}[Logarithmic Resolution]\label{1201}
There is a functor $F_{\textup{log-res}}$ associating to \begin{quote}
    a pure-dimensional, reduced, fs\footnote{Appendix \ref{B.2}} logarithmic Deligne-Mumford stack $X$ of finite type over a field $\kk$ of characteristic zero
\end{quote}
a proper and birational morphism $\Pi \colon F_{\textup{log-res}}(X) \rightarrow X$, such that: \begin{enumerate}
    \item $F_{\textup{log-res}}(X)$ is a pure-dimensional, toroidal Deligne-Mumford stack over $\kk$;
    \item $\Pi$ is an isomorphism over the logarithmically smooth locus $X^{\textup{log-sm}}$ of $X$;
    \item $\Pi^{-1}(X \setminus X^{\textup{log-sm}})$ is contained in the toroidal divisor of $F_{\textup{log-res}}(X)$.
\end{enumerate}
Functoriality here is with respect to logarithmically smooth morphisms: if $\widetilde{X} \rightarrow X$ is a logarithmically smooth morphism, then $F_{\textup{log-res}}(\widetilde{X}) = F_{\textup{log-res}}(X) \times_X \widetilde{X}$.
\end{theorem}

We emphasize that the fiber product at the end of (\ref{1201}) should be taken in the same category of the theorem. Note this differs from standard notation (for example, in \cite[III.2.1]{Ogu18}), where we would instead write $F_{\textup{log-res}}(\widetilde{X}) = (F_{\textup{log-res}}(X) \times_X \widetilde{X})^{\textup{sat}}$.

\begin{proof}
The proof of this theorem follows the strategy in \cite[Theorem 8.1.1]{ATW19}, with minor modifications. Let $X$ be as in (\ref{1201}). Since \'etale locally $X$ can always be embedded in pure codimension in a toroidal $\kk$-scheme, the theorem follows once we show the following: \begin{quote}
    Given two strict closed embeddings of $X$ into pure-dimensional, toroidal Deligne-Mumford stacks $Y_i$ over $\kk$ (where $i=1,2$), the logarithmic resolutions of $X$ obtained from $F_{\textup{log-ER}}(X \subset Y_i)$ (for $i=1,2$) coincide.
\end{quote}
First assume that $\dim(Y_1) = \dim(Y_2)$: in this case, the two embeddings are \'etale locally isomorphic. By functoriality, the logarithmic embedded resolutions $F_{\textup{log-ER}}^{\circ \infty}(X \subset Y_i)$ (for $i=1,2$) are isomorphic, whence the resulting logarithmic resolutions of $X$ coincide. In general, this reduces to the earlier case, by a repeated application of (\ref{1202}).
\end{proof}

\begin{lemma}[Re-embedding Principle \textrm{\cite[Proposition 2.9.3]{ATW20a}}]\label{1202}
Let $X$ be a reduced, closed substack of pure codimension in a toroidal Deligne-Mumford stack $Y$ over a field $\kk$ of characteristic zero. Let $Y_1$ be the fiber product $Y \times_{\kk} \mathbb{A}^1_\kk$ in the category of logarithmic schemes, where $\mathbb{A}^1_\kk$ and $\kk$ are given the trivial logarithmic structure. Then: \begin{enumerate}
    \item For every $p \in \abs{X}$, $\inv_p(X \subset Y_1)$ is the concatenation $(1,\inv_p(X \subset Y))$.
    \item If $F_{\textup{log-ER}}(X \subset Y) = (X' \subset Y')$ and $F_{\textup{log-ER}}(X \subset Y_1) = (X_1' \subset Y_1')$, then $Y'$ is canonically identified with the proper transform of $Y = Y \times \lbrace 0 \rbrace \subset Y_1$ in $Y_1'$, under which $X' = X_1'$.
\end{enumerate}
\end{lemma}

We prove ({\ref{1202}}) in \S\ref{7.5}. The second step is to deduce the following theorem from (\ref{1201}) via \emph{resolution of toroidal singularities}:

\begin{theorem}[Resolution]\label{1203}There is a functor $F_{\textup{res}}$ associating to \begin{quote}
    a pure-dimensional, reduced Deligne-Mumford stack $X$ of finite type over a field $\kk$ of characteristic zero
\end{quote}
a proper and birational morphism $\Pi \colon F_{\textup{res}}(X) \rightarrow X$, such that: \begin{enumerate}
    \item $F_{\textup{res}}(X)$ is a pure-dimensional, smooth Deligne-Mumford stack over $\kk$.
    \item $\Pi$ is an isomorphism over the smooth locus $X^{\textup{sm}}$ of $X$.
    \item $\Pi^{-1}(X \setminus X^{\textup{sm}})$ is a simple normal crossing divisor on $F_{\textup{res}}(X)$.
\end{enumerate}
Functoriality here is with respect to smooth morphisms: if $\widetilde{X} \rightarrow X$ is a smooth morphism, then $F_{\textup{res}}(\widetilde{X}) = F_{\textup{res}}(X) \times_X \widetilde{X}$.
\end{theorem}

\begin{proof}
Let $X$ be as in the theorem, give $X$ the trivial logarithmic structure, and apply (\ref{1201}) to obtain $F_{\textup{log-res}}(X) \rightarrow X$. Note $X^{\textup{log-sm}} = X^{\textup{sm}}$ in this case. Next, apply \cite[Theorem 6.5.1]{Wlo19}: there is a projective birational, morphism $\phi \colon X'' \rightarrow X' = F_{\textup{log-res}}(X)$, where $X''$ is a pure-dimensional, smooth Deligne-Mumford stack over $\kk$, $\phi$ is an isomorphism over the smooth locus $(X')^{\textup{sm}}$ of $X'$, the preimage of the toroidal divisor on $X'$ under $\phi$ is a simple normal crossing divisor, and $\phi$ is functorial with respect to strict, smooth morphisms of toroidal Deligne-Mumford stacks over $\kk$. Take $\Pi \colon F_{\textup{res}}(X) \rightarrow X$ to be the composition $X'' \xrightarrow{\phi} X' = F_{\textup{log-res}}(X) \rightarrow X$.
\end{proof}

We remark that that if $X$ happens to be a scheme in (\ref{1203}), $F_{\textup{res}}(X)$ is, more often than not, a stack. Therefore, the final step involves \emph{Bergh's destackification theorem}:

\begin{theorem}[Coarse Resolution]\label{1204}
There is a functor $F_{\textup{c-res}}$ associating to \begin{quote}
    a pure-dimensional, reduced Deligne-Mumford stack $X$ of finite type over a field $\kk$ of characteristic zero
\end{quote}
a \emph{\textbf{projective}} and birational morphism $\Pi \colon F_{\textup{c-res}}(X) \rightarrow X$, such that: \begin{enumerate}
    \item $F_{\textup{c-res}}(X)$ is a pure-dimensional, smooth Deligne-Mumford stack over $\kk$.
    \item $\Pi$ is an isomorphism over the smooth locus $X^{\textup{sm}}$ of $X$.
    \item $\Pi^{-1}(X \setminus X^{\textup{sm}})$ is a simple normal crossing divisor on $F_{\textup{c-res}}(X)$.
\end{enumerate}
Functoriality here is with respect to smooth morphisms: if $\widetilde{X} \rightarrow X$ is a smooth morphism, then $F_{\textup{c-res}}(\widetilde{X}) = F_{\textup{c-res}}(X) \times_X \widetilde{X}$. 

In particular, if we restrict to the full subcategory whose objects are pure-dimensional, reduced schemes of finite type over $\kk$, we recover Hironaka's \emph{\cite[Main Theorem I]{Hir64}}.
\end{theorem}

\begin{proof}
This proof follows verbatim as in the proof of \cite[Theorem 8.12]{ATW19}. Let $X$ be as in the theorem, and apply \cite[Theorem 7.1]{BR19} to the standard pair $(F_{\textup{res}}(X),D)$ (where $D$ is the simple normal crossing divisor in (\ref{1203}(iii))) and $F_{\textup{res}}(X) \rightarrow X \rightarrow \Spec(\kk)$ (where the $F_{\textup{res}}(X) \rightarrow X$ is given in (\ref{1203})). This provides a projective morphism $F_{\textup{res}}(X)' \rightarrow F_{\textup{res}}(X) \rightarrow X$, functorial for all smooth morphisms, such that the relative coarse moduli space $F_{\textup{res}}(X)' \rightarrow \underline{F_{\textup{res}}(X)'} \rightarrow X$ is projective over $X$, and such that $F_{\textup{res}}(X)'$ and $\underline{F_{\textup{res}}(X)'}$ are smooth over $\kk$. We then take $\Pi \colon F_{\textup{c-res}}(X) \rightarrow X$ to be $\underline{F_{\textup{res}}(X)'} \rightarrow X$.
\end{proof}

\subsection{Adapting methods in \texorpdfstring{\cite{ATW19}}{[ATW19]} to the logarithmic setting}\label{1.3}

We recall the set-up in \cite{ATW19}: let $\kk$ be a field of characteristic zero (as before), but we instead consider a smooth $\kk$-scheme $Y$, and a reduced closed subscheme $X \subset Y$ of pure codimension $c$ --- or more generally, a reduced closed substack $X$ of pure codimension $c$ in a smooth Deligne-Mumford stack $Y$ over $\kk$. Then \cite{ATW19} proposes a faster and simpler approach to embedded resolution of singularities of $X$ in $Y$, where each step immediately and visibly improves the singularities --- by considering a \emph{broader} notion of blow-up centers \cite[2.4]{ATW19}. However, as the example in \cite[8.3]{ATW19} demonstrates, \begin{quote}
    $(\ast)$ at each step of the resolution, the chosen blow-up center does not necessarily have simple normal crossings with the exceptional loci obtained at that step,
\end{quote}
and hence, \begin{quote}
    $(\dagger)$ the exceptional loci at subsequent steps of the resolution may not be simple normal crossing divisors.
\end{quote} Consequently, this does not address a key feature of Hironaka's \cite[Main Theorem I]{Hir64}: \begin{quote}
    $(\diamondsuit)$ namely, the preimage of the singular locus of $X$ under the resolution in \cite{ATW19} is not always a simple normal crossing divisor.
\end{quote}
Our result (\ref{1101}) on logarithmic embedded resolution can be seen as a resolution to the aforementioned issue as follows: \begin{enumerate}
    \item[(I)] Give $Y$ (and hence $X$) the trivial logarithmic structure. At each step of the resolution, we will first encode the exceptional divisor obtained at that step into the logarithmic structure (cf. (\ref{1101}(iv))).
    \item[(II)] With respect to these logarithmic structures, we then adapt the methods used in \cite[\S 5]{ATW19} to obtain a \emph{toroidal} blow-up center at each step. 
\end{enumerate}
We remark that (II) does not resolve $(\ast)$ or $(\dagger)$: in fact, the Deligne-Mumford stacks $Y_i$ obtained in this modified resolution $Y_N \rightarrow \dotsb \rightarrow Y_1 \rightarrow Y_0 = Y$ may not even be smooth (unlike in \cite{ATW19}). For an example of this, see \S\ref{8.1}. Nevertheless, (II) assures that the $Y_i$ will be toroidal. Moreover, (I) and (II) will give us some control over the exceptional loci obtained in the process: namely, the exceptional loci at each step will be always contained in the toroidal divisor (\ref{B105}(ii)) of the toroidal Deligne-Mumford stack $Y_i$ at that step. Consequently, the preimage of the singular locus of $X$ under this modified resolution $X_N = Y_N \times_Y X \rightarrow X$ is contained in the toroidal divisor of $X_N$ (cf. (\ref{1101}(c))). One then resolves the issue in $(\diamondsuit)$ via \emph{resolution of toroidal singularities} as in (\ref{1203}).

This justifies our need to work in the logarithmic setting as outlined in \S\ref{1.1}. We also remark that the above strategy (I)$+$(II) was already pursued earlier in \cite{ATW20a}, although with respect to Hironaka's classical resolution algorithm.

\subsection{Acknowledgements}\label{1.4}

I would like to express my gratitude to my advisor Dan Abramovich for suggesting this project, and his continued patience and guidance. I would also like to thank Jonghyun Lee, David Rydh, Michael Temkin, Jaros{\l}aw W{\l}odarczyk, and the referee for clarifications, corrections, comments, and questions, as well as Tangli Ge and Stephen Obinna for discussions during early stages of this project.

\tableofcontents

\section{Idealistic exponents}\label{chpt2}

In this chapter, let $\kk$ be a field, and unless otherwise stated, $Y$ is usually a $\kk$-variety\footnote{Following \cite{Har77}, a $\kk$-variety is an integral, separated scheme of finite type over $\kk$.} with field of fractions $K$. Let $\ZR(Y)$ denote the Zariski-Riemann space of $Y$, as defined in \cite[2.1]{ATW19} or Appendix \ref{appA} of this paper. $\ZR(Y)$ is a locally ringed space, whose elements are valuation rings $R_\nu$ of $K$ containing $\kk$ which possess a center on $Y$ \cite[Exercise II.4.5]{Har77}. We usually denote $R_\nu$ by its corresponding valuation $\nu \colon K^\ast \twoheadrightarrow G_\nu$, where $G_\nu$ is the value group of $\nu$. The monoid of non-negative elements of $G_\nu$ is denoted $(G_\nu)_+$. Let us fix some related notation for this chapter:
\renewcommand{\arraystretch}{1.2}
\begin{quote}
\begin{tabular}{ c c l }
 $\SO_{\ZR(Y)}$ & --- & sheaf of rings carried by $\ZR(Y)$, whose stalk at $\nu$ is $R_\nu$ \\ 
 $\Gamma_Y$ & --- & sheaf of ordered groups $K^\ast/\SO_{\ZR(Y)}^\ast$ on $\ZR(Y)$, whose stalk at $\nu$ is $G_\nu$ \\  
 $\Gamma_{Y,+}$ & --- & subsheaf of $\Gamma_Y$ consisting of non-negative sections of $\Gamma$ \\
 $y_\nu$ & --- & the (unique) center of $\nu$ on $Y$ \\
 $\pi_Y$ & --- & the canonical morphism $\ZR(Y) \rightarrow Y$ which maps $\nu \mapsto y_\nu$
\end{tabular}
\end{quote}
See Appendix \ref{appA} for a self-contained exposition of the aforementioned notions.

\subsection{Valuative ideals}\label{2.1}

Following \cite[2.2]{ATW19}, a \emph{valuative ideal} over $Y$ is defined to be a section \[
    \gamma \in H^0(\ZR(Y),\Gamma_{Y,+}).
\]
A (coherent) ideal $0 \neq \mathcal{I} \subset \SO_Y$ determines a valuative ideal $\gamma_\mathcal{I}$ over $Y$ as follows. For every $\nu \in \ZR(Y)$, remember that $y_\nu$ denotes the center of $\nu$ on $Y$, and let $f_\nu^\# \colon \SO_{Y,y_\nu} \rightarrow R_\nu$ denote the corresponding local $\kk$-homomorphism. We then set \begin{align*}
    \gamma_{\mathcal{I},\nu} := \min\lbrace \nu(g) \colon g \textup{ is a non-zero section of } \mathcal{I} \subset \SO_Y \textup{ over an open set containing $y_\nu$}\rbrace,
\end{align*}
where $\nu(g)$ is an abbreviation for $\nu(f_\nu^\#(g_{y_\nu}))$ (this abbreviation will persist in this paper). Note that since $\mathcal{I}$ is coherent, this minimum exists in $(G_\nu)_+$. Indeed, if $\mathcal{I}_{y_\nu}$ is generated by $g_1,\dotsc,g_r \in \SO_{Y,y_\nu}$, then $\gamma_{\mathcal{I},\nu} = \min\lbrace \nu(g_i) \colon 1 \leq i \leq r \rbrace$. 

Moreover, if we let $1 \leq j \leq r$ be such that $\nu(g_j) = \gamma_{\mathcal{I},\nu}$, then $f_\nu^\#(\mathcal{I}_{y_\nu})R_\nu$ is the principal ideal $(f_\nu^\#(g_j))R_\nu$ of $R_\nu$. For \begin{align*}
    (\gamma_{\mathcal{I},\nu})_{\nu \in \ZR(Y)} \in \prod_{\nu \in \ZR(Y)}{(G_\nu)_+}
\end{align*}
to define a valuative ideal $\gamma_\mathcal{I}$ over $Y$, we need to check that it is a compatible collection of germs of $\Gamma_{Y,+}$. Indeed, fix an arbitrary $\nu \in \ZR(Y)$, and assume that $g_1,\dotsc,g_r \in \SO_{Y,y_\nu}$ generate $\mathcal{I}_{y_\nu}$, with $\nu(g_j) = \min\lbrace \nu(g_i) \colon 1 \leq i \leq r \rbrace$. There exists an affine open neighbourhood $V_\nu$ of $y_\nu$ in $Y$ such that $g_1,\dotsc,g_r$ extend to sections of $\mathcal{I}$ over $V_\nu$ which generate the stalk of $\mathcal{I}$ at every point in $V_\nu$. Then $U_\nu = \pi_Y^{-1}(V_\nu) \cap U\big(\frac{g_i}{g_j} \colon i \neq j\big)$ is an open neighbourhood of $\nu$ in $\ZR(Y)$ such that for all $\nu' \in U_\nu$, $\gamma_{\nu'} = \nu'(g_j)$.

In fact, the same argument shows that any valuative ideal over $Y$ arising from an ideal on $Y$ is locally represented by generators of that ideal:

\begin{lemma}\label{2101}
Let the notation be as above, and let $\mathcal{I}$ be a non-zero ideal on $Y$. There exist: \begin{enumerate}
    \item a finite open affine cover $\mathcal{V} = \lbrace V_\ell \colon 1 \leq \ell \leq m \rbrace$ of $Y$;
    \item for each $1 \leq \ell \leq m$, a finite open cover $\mathcal{U}_\ell = \lbrace U_{\ell,j} \colon 1 \leq j \leq r_\ell \rbrace$ of $\pi_Y^{-1}(V_\ell)$;
    \item for each $1 \leq \ell \leq m$, sections $\lbrace g_{\ell,j} \colon 1 \leq j \leq r_\ell \rbrace$ of $\mathcal{I}$ over $V_\ell$, which generate $\mathcal{I}$ at every point of $V_\ell$,
\end{enumerate} 
such that for each $1 \leq \ell \leq m$, each $1 \leq j \leq r_\ell$, and every $\nu \in U_{\ell,j}$, we have $\gamma_{\mathcal{I},\nu} = \nu(g_{\ell,j})$.
\end{lemma}

\begin{proof}
For every $y \in Y$, there exists $g_1,\dotsc,g_r \in \mathcal{I}_y$ and an open affine neighbourhood $y \in V_y \subset Y$ such that $g_1,\dotsc,g_r$ extend to sections of $\mathcal{I}$ over $V_y$ generating $\mathcal{I}$ at every point of $V_y$. Since $Y$ is quasi-compact, there exists a finite open subcover of $\lbrace V_y \colon y \in Y \rbrace$, say $\mathcal{V} = \lbrace V_\ell \colon 1 \leq \ell \leq m \rbrace$. For each $\ell$, let $g_{\ell,1},\dotsc,g_{\ell,r_\ell} \in \mathcal{I}(V_\ell)$ be the sections chosen earlier.

For each $1 \leq j \leq r_\ell$, let $U_{\ell,j} = \pi_Y^{-1}(V_\ell) \cap U\big(\frac{g_{\ell,i}}{g_{\ell,j}} \colon i \neq j\big)$. For all $\nu \in \pi_Y^{-1}(V_\ell)$, we have $y_\nu \in V_\ell$, whence $\mathcal{I}_{y_\nu}$ is generated by $\lbrace g_{\ell,j} \colon 1 \leq j \leq r_\ell \rbrace$, so $\gamma_{\mathcal{I},\nu} = \min\lbrace \nu(g_{\ell,j}) \colon 1 \leq j \leq r_\ell \rbrace$. From this, it is immediate that $\pi_Y^{-1}(V_\ell) = \bigcup_{j=1}^{r_\ell}{U_{\ell,j}}$. The conclusion is also immediate.
\end{proof}

\begin{definition}[Idealistic classes]\label{2102}
Let $Y$ be a $\kk$-variety. A valuative ideal $\gamma$ over $Y$ associated to a non-zero ideal $\mathcal{I}$ on $Y$ is called an \emph{idealistic class} over $Y$.
\end{definition}

Conversely, every valuative ideal $\gamma$ over $Y$ determines an ideal $\mathcal{I}_\gamma$ on $Y$: we let $\mathcal{I}_\gamma$ be the subsheaf of $\SO_Y$ whose sections $g$ over an open set $U$ satisfy $\nu(g) \geq \gamma_\nu$ for every $\nu \in \pi_Y^{-1}(U) \subset \ZR(Y)$ (namely those $\nu$ such that $y_\nu \in U$).

Before moving on, we recall the following definition \cite[9.6.A]{Laz04}: if $I$ is an ideal of a ring $A$, the \emph{integral closure} $\overline{I}$ of $I$ in $A$ consists of elements $x \in A$ which satisfy a ``weighted integral equation'' \begin{align*}
    x^n + a_1x^{n-1}+\dotsb+a_{n-1}x+a_n=0, \qquad \textup{where } a_i \in I^i.
\end{align*}
We say $I$ is \emph{integrally closed} in $A$ if $I = \overline{I}$. Observe that $I \subset \overline{I} \subset \sqrt{I}$ (where $\sqrt{I}$ is the radical of $I$). In \S\ref{2.2}, we will prove that: \begin{enumerate}
    \item[(a)] $\overline{I}$ is an ideal of $A$;
    \item[(b)] if $\mathcal{I}$ is an ideal on a $\kk$-variety $Y$, the presheaf on $Y$ given by $U \mapsto \overline{\mathcal{I}(U)}$ is a sheaf, denoted by $\overline{\mathcal{I}}$;
    \item[(c)] and the following lemma, which was essentially noted in \cite{Hir77}:
\end{enumerate}

\begin{lemma}\label{2103}
Let the notation be as above. \begin{enumerate}
    \item If $\gamma$ is a valuative ideal over $Y$, then $\mathcal{I}_\gamma$ is integrally closed in $\SO_Y$.
    \item Let $0 \neq \mathcal{I} \subset \SO_Y$ be an ideal, with associated idealistic class $\gamma = \gamma_\mathcal{I}$ over $Y$. Then $\mathcal{I}_\gamma = \overline{\mathcal{I}}$.
\end{enumerate}
\end{lemma}

\begin{corollary}\label{2104}
Let $Y$ be a $\kk$-variety. The above describes a one-to-one correspondence between non-zero, integrally closed ideals of $\SO_Y$ and idealistic classes over $Y$. \qed
\end{corollary}

\subsection{Rees algebras and valuative \texorpdfstring{$\QQ$}{Q}-ideals}\label{2.2}

\begin{definition}[Rees algebras]\label{2201}
Given a scheme $Y$, a \emph{Rees algebra} on $Y$ is a finitely generated, quasi-coherent, graded $\SO_Y$-subalgebra $\SR = \bigoplus_{m \in \NN}{\mathcal{I}_m \cdot T^m} \subset \SO_Y[T]$, satisfying $\mathcal{I}_0 = \SO_Y$ and $\mathcal{I}_m \supset \mathcal{I}_{m+1}$ for every $m \in \NN$. We say $\SR$ is \emph{non-zero} if $\mathcal{I}_m \neq 0$ for some $m \geq 1$.
\end{definition}

Recall that we can associate a Rees algebra to every ideal $\mathcal{I}$ on $Y$, namely $\bigoplus_{m \in \NN}{\mathcal{I}^m \cdot T^m}$. This sets up a one-to-one correspondence: \begin{align*}
    \lbrace \textup{ideals of $\SO_Y$} \rbrace \leftrightarrow \lbrace \textup{Rees algebras generated in degree $1$}\rbrace.
\end{align*}

For the remainder of this section, let $Y$ be a $\kk$-variety $Y$. Accompanying the notion of a Rees algebra on $Y$ is the notion of a \emph{valuative $\QQ$-ideal} over $Y$ \cite[2.2]{ATW19}. To define this notion, consider the sheaf of ordered groups $\Gamma_{Y,\QQ} = \QQ \otimes \Gamma_Y$. We denote the sheaf of monoids consisting of non-negative sections of $\Gamma_{Y,\QQ}$ by $\Gamma_{Y,\QQ+}$. A valuative $\QQ$-ideal over $Y$ is a section $\gamma$ in $H^0(\ZR(Y),\Gamma_{Y,\QQ+})$. Note that since $\gamma$ is locally constant and $\ZR(Y)$ is quasi-compact, there exists a sufficiently large natural number $N \geq 1$ such that $N \cdot \gamma$ is a valuative ideal over $Y$.

A non-zero Rees algebra $\SR$ on $Y$ determines a valuative $\QQ$-ideal $\gamma_\SR$ over $Y$ by: \begin{align*}
    \gamma_\SR := (\gamma_{\SR,\nu})_{\nu \in \ZR(Y)} \in \prod_{\nu \in \ZR(Y)}{(\QQ \otimes G_\nu)_+}
\end{align*}
where \begin{align*}
    \gamma_{\SR,\nu} := \min\left\lbrace \frac{1}{n} \cdot \nu(g) \colon 0 \neq gT^n \textup{ is a section of } \SR \textup{ over an open set containing } y_\nu \textup{ (for $n \geq 1$}) \right\rbrace.
\end{align*}
Again, we have to show (i) this minimum exists in $(\QQ \otimes G_\nu)_+$, and (ii) $(\gamma_{\SR,\nu})_{\nu \in \ZR(Y)}$ defines a compatible collection of germs, and hence, defines a valuative $\QQ$-ideal over $Y$. Indeed, fix $\nu \in \ZR(Y)$, and suppose $g_1T^{n_1},\dotsc,g_rT^{n_r}$ generate $\SR_{y_\nu}$ as a $\SO_{Y,y_\nu}$-algebra. Then we claim: \begin{align*}
    \gamma_{\SR,\nu} = \min\left\lbrace \frac{1}{n_i} \cdot \nu(g_i) \colon 1 \leq i \leq r \right\rbrace
\end{align*}
from which (i) is immediate. Indeed, suppose $gT^n \in \mathcal{I}_{y_\nu}$. Then we can write \begin{align*}
    gT^n &= \sum_{k_1n_1+\dotsb+k_rn_r=n}{a_{\vec{k}} \cdot \prod_{i=1}^r{(g_iT^{n_i})^{k_i}}} \qquad \textup{in } \SO_{Y,y_\nu}[T]
\end{align*}
which means \begin{align*}
    g = \sum_{k_1n_1+\dotsb+k_rn_r=n}{a_{\vec{k}} \cdot \prod_{i=1}^r{g_i^{k_i}}} \qquad \textup{in } \SO_{Y,y_\nu}.
\end{align*}
Consequently, \begin{align*}
    \frac{1}{n} \cdot \nu(g) &\geq \min\left\lbrace \frac{1}{n}\sum_{i=1}^r{k_in_i \cdot \left(\frac{1}{n_i} \cdot \nu(g_i)\right)} \colon k_1n_1+\dotsb+k_rn_r=n \right\rbrace \\
    &\geq \min\left\lbrace \left(\frac{1}{n}\sum_{i=1}^r{k_in_i}\right) \cdot \min\left\lbrace \frac{1}{n_i} \cdot \nu(g_i) \colon 1 \leq i \leq r \right\rbrace \colon k_1n_1+\dotsb+k_rn_r = n \right\rbrace \\
    &= \min\left\lbrace \frac{1}{n_i} \cdot \nu(g_i) \colon 1 \leq i \leq r\right\rbrace
\end{align*}
as desired. 

For (ii), there exist an affine open neighbourhood $V_\nu$ of $y_\nu$ in $Y$ such that $g_1T^{n_1},\dotsc,g_rT^{n_r}$ extend to sections of $\SR$ over $V_\nu$ which generate the stalk of $\SR$ at every point in $V_\nu$. Let $1 \leq j \leq r$ such that $\gamma_{\SR,\nu} = \frac{1}{n_j} \cdot \nu(g_j)$. Then $U_\nu = \pi_Y^{-1}(V_\nu) \cap U\left(\frac{g_i^{n_j}}{g_j^{n_i}} \colon i \neq j \right)$ is an open neighbourhood of $\nu$ in $\ZR(Y)$ such that for all $\nu' \in U_\nu$, $\gamma_{\nu'} = \frac{1}{n_j} \cdot \nu'(g_j)$.

Note that if $\SR$ is the Rees algebra of an ideal $0 \neq \mathcal{I} \subset \SO_Y$, then $\gamma_\SR = \gamma_\mathcal{I}$.

Lastly, imitating the proof of (\ref{2101}) yields the analogous lemma: 

\begin{lemma}\label{2202}
Let the notation be as above, and let $\SR$ be a non-zero Rees algebra on $Y$. There exist: \begin{enumerate}
    \item a finite open affine cover $\mathcal{V} = \lbrace V_\ell \colon 1 \leq \ell \leq m \rbrace$ of $Y$;
    \item for each $1 \leq \ell \leq m$, a finite open cover $\mathcal{U}_\ell = \lbrace U_{\ell,j} \colon 1 \leq j \leq r_\ell \rbrace$ of $\pi_Y^{-1}(V_\ell)$;
    \item for each $1 \leq \ell \leq m$, sections $\lbrace g_{\ell,j}T^{n_{\ell,j}} \colon 1 \leq j \leq r_\ell \rbrace$ of $\SR$ over $V_\ell$, which generate $\SR$ at every point of $V_\ell$ (as a $\SO_{Y,y}$-algebra),
\end{enumerate} 
such that for each $1 \leq \ell \leq m$, each $1 \leq j \leq r_\ell$, and every $\nu \in U_{\ell,j}$, we have $\gamma_{\SR,\nu} = \frac{1}{n_{\ell,j}} \cdot \nu(g_{\ell,j})$. \qed
\end{lemma}

\begin{definition}[Idealistic exponents, cf. \textup{\cite[Definition 3]{Hir77}}]\label{2203}Let $Y$ be a $\mathds{k}$-variety. A valuative $\QQ$-ideal $\gamma$ over $Y$ associated to some non-zero Rees algebra on $Y$ is called an \emph{idealistic exponent} over $Y$.
\end{definition}

Conversely, let $\gamma$ be a valuative $\QQ$-ideal over $Y$. As in \S\ref{2.1}, $\gamma$ also determines an ideal $\mathcal{I}_\gamma$ on $Y$ whose sections $g$ over an open set $U$ satisfy $\nu(g) \geq \gamma_\nu$ for every $\nu \in \pi_Y^{-1}(U) \subset \ZR(Y)$ (namely those $\nu$ such that $y_\nu \in U$). But $\gamma$ also determines a $\SO_Y$-subalgebra of $\SO_Y[T]$: \begin{align*}
    \SR_\gamma = \bigoplus_{m \in \NN}{\mathcal{I}_{m \cdot \gamma} \cdot T^m} \subset \SO_Y[T]
\end{align*}
where $\mathcal{I}_{m \cdot \gamma}$ is the ideal of $\SO_Y$ associated to the multiple $m \cdot \gamma$ (which was just described). In general, $\SR_\gamma$ is not a Rees algebra on $Y$, but (\ref{2205}) below says $\SR_\gamma$ is a Rees algebra on $Y$ whenever $\gamma$ is an idealistic exponent over $Y$. Note that $\SR_\gamma$ contains the Rees algebra of $\mathcal{I}_\gamma$, but of course these are rarely equal --- see (\ref{2206}).

\begin{lemma}\label{2204}
Let the notation be as above, and let $\gamma$ be a valuative $\QQ$-ideal over $Y$. The corresponding $\SO_Y$-subalgebra $\SR_\gamma$ of $\SO_Y[T]$ is integrally closed in $\SO_Y[T]$.
\end{lemma}

\begin{proof}
(cf. \cite[Chapter I, Lemma 1]{KKMSD73}) Since the integral closure of $\SR_\gamma$ in $\SO_Y[T]$ is a subring of $\SO_Y[T]$, it suffices to show that whenever a non-zero homogeneous section $gT^r$ of $\SO_Y[T]$ over an open set $U \subset Y$ satisfies an equation of the form \begin{align*}
    (gT^r)^n + a_1(gT^r)^{n-1} + \dotsb + a_{n-1}(gT^r)+a_n = 0, \qquad a_i \in \SR_\gamma(U),
\end{align*}
then $gT^r$ is a section of $\SR_\gamma$ over $U$. By writing each $a_i$ as a sum of homogeneous sections in $\SR_\gamma(U)$ and comparing degrees, we may assume that each $a_i$ is $\alpha_i T^{ir}$ for some $\alpha_i \in \mathcal{I}_{ir \cdot \gamma}(U)$. If $r = 0$, there is nothing to show. If $r > 0$, we have \begin{align*}
    g^n + \alpha_1 g^{n-1} + \dotsb + \alpha_{n-1} g + \alpha_n = 0 \qquad \textup{in $\SO_Y(U)$}.
\end{align*}
Let $\nu \in \pi_Y^{-1}(U) \subset \ZR(Y)$. We claim that there must exist some $1 \leq j \leq n$ such that $j \cdot \nu(g) \geq \nu(\alpha_j)$. Indeed, if not, then $i \cdot \nu(g) < \nu(\alpha_i)$ for all $1 \leq i \leq n$, so $\nu(g^n) < \nu(\alpha_i g^{n-i})$ for all $1 \leq i \leq n$. This implies $g^n + \alpha_1g^{n-1} + \dotsb + \alpha_{n-1}g + \alpha_n \neq 0$, a contradiction. Now our claim implies $\nu(g) \geq \frac{1}{j}\nu(\alpha_j) \geq r \cdot \gamma_\nu$, so $g \in \mathcal{I}_{r \cdot \gamma}(U)$. Since $\nu \in \pi_Y^{-1}(U)$ is arbitrary, $gT^r \in \SR_\gamma(U)$.
\end{proof}

A special case of the next theorem is noted in \cite[3.4]{ATW19}:

\begin{proposition}\label{2205}
Let the notation be as above, and let $\gamma = \gamma_\SR$ be the idealistic exponent over $Y$ associated to a non-zero Rees algebra $\SR$ on $Y$. Then $\SR_\gamma$ is the integral closure of $\SR$ in $\SO_Y[T]$. In particular, $\SR_\gamma$ is a finite $\SR$-module, and hence a Rees algebra on $Y$.
\end{proposition}

\begin{proof}
(cf. \cite[Theorem 10.4]{Mat89}) By (\ref{2204}), $\SR_\gamma$ contains the integral closure of $\SR$ in $\SO_Y[T]$. We can check the converse on stalks: let $y \in Y$, and it suffices to show that whenever a homogeneous element $gT^n$ of $\SO_{Y,y}[T]$ is not integral over $\SR_y$, then $gT^n$ is not in $(\SR_\gamma)_y$. Fix a set of generators $g_1T^{n_1},\dotsc,g_rT^{n_r}$ of $\SR_y$ as a $\SO_{Y,y}$-algebra: then our goal is to find $\nu \in \ZR(Y)$ whose center $y_\nu$ on $Y$ is $y$, and such that \begin{align*}
    \frac{1}{n}\nu(g) < \min\left\lbrace\frac{1}{n_i}\nu(g_i) \colon 1 \leq i \leq r\right\rbrace.
\end{align*}
Let $A = \SO_{Y,y}\big[\frac{g_i^n}{g^{n_i}} \colon 1 \leq i \leq r\big]$, a subring of $K$ containing $\mathds{k}$. Let $I$ be the ideal of $A$ generated by $\big\lbrace\frac{g_i^n}{g^{n_i}} \colon 1 \leq i \leq r\big\rbrace$ and the maximal ideal $\mathfrak{m}_{Y,y}$ of $\SO_{Y,y}$. We claim that $1 \notin I$. If not, \begin{align*}
    1 = \alpha + \sum_{\substack{J=(j_1,\dotsc,j_r) \\ j_1+\dotsb+j_r \geq 1}}{\beta_J\prod_{i=1}^r\left(\frac{g_i^n}{g^{n_i}}\right)^{j_i}}
\end{align*}
where $\alpha \in \mathfrak{m}_{Y,y}$ and only finitely many $\beta_J \in \SO_{Y,y}$ are non-zero. Since $1-\alpha$ is a unit in $\SO_{Y,y}$, we may assume $\alpha = 0$. For each $1 \leq i \leq r$, let $t_i = \max\lbrace j_i \colon$there exists $J = (j_1,\dotsc,j_r)$ such that $\beta_J \neq 0 \rbrace$. Multiplying the above equation throughout by $\prod_{i=1}^r(g^{n_i})^{t_i} = g^{\sum_{i=1}^r{n_it_i}}$, we get \begin{align*}
    g^{\sum_{i=1}^r{n_it_i}} = \sum_{\substack{J=(j_1,\dotsc,j_r) \\ j_1+\dotsb+j_r \geq 1}}{\beta_J\prod_{i=1}^r{\left(g_i^{nj_i} \cdot g^{n_i(t_i-j_i)}\right)}} = \sum_{\substack{J=(j_1,\dotsc,j_r) \\ j_1+\dotsb+j_r \geq 1}}{\left(\beta_J\prod_{i=1}^r{g_i^{nj_i}}\right) \cdot g^{\sum_{i=1}^r{n_i(t_i-j_i)}}}
\end{align*}
which implies \begin{align*}
    (gT^n)^{\sum_{i=1}^r{n_it_i}} - \sum_{\substack{J=(j_1,\dotsc,j_r) \\ j_1+\dotsb+j_r \geq 1}}{\left(\beta_J\prod_{i=1}^r{(g_iT^{n_i})^{nj_i}}\right) \cdot \left(gT^n\right)^{\sum_{i=1}^r{n_i(t_i-j_i)}}} = 0
\end{align*}
which is an integral equation for $gT^n$ over $\SR_y = \SO_{Y,y}[g_iT^{n_i} \colon 1 \leq i \leq r]$, a contradiction. Therefore, $I$ is a proper ideal of $A$, so there exists a maximal ideal $\mathfrak{p}$ of $A$ containing $I$. By \cite[Theorem 10.2]{Mat89}, there exists $\nu \in \ZR(K,\mathds{k})$\footnote{Recall that $K$ denotes the field of fractions of $Y$, and an element $\nu \in \ZR(K,\mathds{k})$ is a valuation ring $R_\nu$ of $K$ containing $\mathds{k}$, as defined in Appendix \ref{A.1}.} such that $R_\nu \supset A$ and $\mathfrak{m}_\nu \cap A = \mathfrak{p}$. Consequently, $\lbrace \frac{g_i^n}{g^{n_i}} \colon 1 \leq i \leq r \rbrace \subset \mathfrak{p} \subset \mathfrak{m}_\nu$, whence for each $1 \leq i \leq r$, $\frac{g^{n_i}}{g_i^n} \notin R_\nu$. This means that for each $1 \leq i \leq r$, \begin{align*}
    \nu\left(\frac{g^{n_i}}{g_i^n}\right) < 0, \qquad \textup{which implies} \qquad \frac{1}{n}\nu(g) < \frac{1}{n_i}\nu(g_i),
\end{align*}
as desired. Moreover, $\mathfrak{p} \cap \SO_{Y,y} = \mathfrak{m}_{Y,y}$, so $\mathfrak{m}_\nu \cap \SO_{Y,y} = \mathfrak{m}_{Y,y}$. Thus, the center of $\nu$ on $Y$ is necessarily $y$ (in particular, $\nu \in \ZR(Y)$).
\end{proof}

\begin{corollary}\label{2206}
Let $\mathcal{I}$ be a non-zero ideal on a $\mathds{k}$-variety $Y$, with associated idealistic class $\gamma = \gamma_\mathcal{I}$ over $Y$. Then the Rees algebra $\SR_\gamma$ associated to $\gamma$ is the integral closure of the Rees algebra of $\mathcal{I}$ in $\SO_Y[T]$.
\end{corollary}

\begin{proof}
If $\SR$ is the Rees algebra of $\mathcal{I}$, we noted earlier that $\gamma_\SR = \gamma_\mathcal{I}$. Apply (\ref{2205}).
\end{proof}

\begin{corollary}\label{2207}
Let $Y$ be a $\mathds{k}$-variety. The above describes a one-to-one correspondence between non-zero, integrally closed Rees algebras on $Y$ and idealistic exponents over $Y$. \qed
\end{corollary}

\begin{notation}\label{2208}
In light of (\ref{2207}), the following notation in \cite{ATW20a} makes sense, and we adopt it moving ahead: if $\SR$ is the integral closure of a non-zero Rees algebra generated by sections $g_1^{a_1}T^{b_1},\dotsc,g_r^{a_r}T^{b_r}$, we record $\SR$ as $\SR = (g_1^{q_1},\dotsc,g_r^{q_r})$, where $q_i = \frac{a_i}{b_i}$ for $1 \leq i \leq r$. Note that since $\SR$ is integrally closed, this expression is well-defined, independent of the presentation of $q_i$ as a quotient of two positive integers. Moreover, if we write $\SR = (g_1^{q_1},\dotsc,g_r^{q_r},\mathcal{I}^{q})$ for an ideal $\mathcal{I} \subset \SO_Y$ and a positive rational number $q=\frac{a}{b}$, we mean that $\SR$ is the integral closure of a Rees algebra generated by sections $g_1^{a_1}T^{b_1},\dotsc,g_r^{a_r}T^{b_r}$ and $\lbrace g^aT^b \colon g$ is a section of $\mathcal{I} \rbrace$. For a positive rational number $s$, we write $\SR^s$ to mean $(g_1^{q_1s},\dotsc,g_r^{q_rs},\mathcal{I}^{qs})$. By convention, we shall write $\SR^0$ to mean the trivial Rees algebra $(1) = \SO_Y[T]$.
\end{notation}

Finally, let us tie some loose ends from the end of \S\ref{2.1}. Note that if $I$ is an ideal of a ring $A$, the Rees algebra of $I$ is integrally closed in $A[T]$ if and only if $I^r$ is integrally closed in $A$ for all $r \geq 1$. In particular, $\overline{I}$ is the degree $1$ part of the integral closure of the Rees algebra of $I$ in $A[T]$, so it must be an ideal of $A$. This is assertion (a) before (\ref{2103}), and assertion (b) is proven similarly. We also deduce (\ref{2103}) from results in this section:

\begin{proof}[Proof of \emph{(\ref{2103})}]
Let $\gamma$ be a valuative ideal over $Y$. By (\ref{2204}), $\SR_\gamma$ is integrally closed in $\SO_Y[T]$. Hence, $\mathcal{I}_\gamma^r$ is integrally closed in $\SO_Y$ for all $r \geq 1$. In particular, we get (i).

For (ii), let $\SR$ be the Rees algebra of $\mathcal{I}$, and apply (\ref{2206}): $\SR_\gamma$ is the integral closure of $\SR$ in $\SO_Y[T]$. In particular, the degree $1$ part of $\SR_\gamma$ is $\overline{\mathcal{I}}$, so $\mathcal{I}_\gamma = \overline{\mathcal{I}}$.
\end{proof}

\subsection{Functoriality with respect to dominant morphisms}\label{2.3}

Let $f \colon Y' \rightarrow Y$ be a dominant morphism of $\mathds{k}$-varieties. In Appendix \ref{A.3}, we noted that $f$ naturally induces a morphism $\ZR(f) \colon \ZR(Y') \rightarrow \ZR(Y)$ of locally ringed spaces, which induces a morphism of ordered groups $\Gamma_Y \rightarrow \ZR(f)_\ast\Gamma_{Y'}$, as well as a morphism of sheaves of monoids $\Gamma_{Y,+} \rightarrow \ZR(f)_\ast\Gamma_{Y',+}$. Tensoring with $\QQ$, we also get a morphism of ordered groups $\Gamma_{Y,\QQ} \rightarrow \ZR(f)_\ast\Gamma_{Y',\QQ}$, and a morphism of sheaves of monoids $\Gamma_{Y,\QQ+} \rightarrow \ZR(f)_\ast\Gamma_{Y',\QQ+}$. In particular, for every valuative ideal (resp. valuative $\QQ$-ideal) $\gamma$ over $Y$, we can consider the pullback of $\gamma$ to $Y'$, denoted $\gamma\SO_{Y'}$ (following \cite{ATW19}). If $\gamma = \gamma_\mathcal{I}$ for some ideal $0 \neq \mathcal{I}$ on $\SO_Y$, then $\gamma\SO_{Y'}$ is simply $\gamma_{\mathcal{I}\SO_{Y'}}$. Likewise, if $\gamma = \gamma_\SR$ for some non-zero Rees algebra $\SR$, then $\gamma\SO_{Y'}$ is simply $\gamma_{\SR\SO_{Y'}}$. More generally, whenever $Y' \rightarrow Y$ is a morphism of $\mathds{k}$-varieties with $\mathcal{I}\SO_{Y'} \neq 0$ (resp. $\SR\SO_{Y'}$ non-zero), the pullback $\gamma\SO_Y$ of $\gamma=\gamma_\mathcal{I}$ (resp. $\gamma = \gamma_\SR$) is well-defined.

\section{Toroidal centers}\label{chpt3}

\subsection{Reminders}\label{3.1}

For the remainder of this paper, $\kk$ denotes a field of characteristic zero. Although toroidal Deligne-Mumford stacks over $\kk$ (\ref{B201}) are the main objects of study in our paper (as mentioned in \S\ref{1.1}), a significant portion of the paper instead deals with strict toroidal $\kk$-schemes (\ref{B106}). There are two reasons for this: \begin{enumerate}
    \item[(a)] \'Etale locally a toroidal Deligne-Mumford stack over $\kk$ is a strict toroidal $\kk$-scheme (see paragraph after (\ref{B201})).
    \item[(b)] The constructions and discussions in this paper are \'etale-local. This was hinted in (\ref{1101}).
\end{enumerate}
Henceforth, we shall assume $Y$ is a strict toroidal $\kk$-scheme (with the exception of \S\ref{6.4}), and denote its logarithmic structure by $\alpha_Y \colon \SM_Y \rightarrow \SO_Y$. All ideals $\mathcal{I}$ of $\SO_Y$ considered in the sequel are always assumed to be coherent. Let us recall the following notions from Appendix \ref{appB}:
\begin{quote}
\begin{tabular}{ c c l }
 $\overline{\SM}_Y$ & --- & characteristic of $\SM_Y$, defined as $\SM_Y/\SO_Y^\ast$ \\
 $\fs_y$ & --- & logarithmic stratum through a point $y \in Y$ \\
 $\SD^1_Y$ & --- & logarithmic tangent sheaf of $Y$ over $\kk$ \\
 $\SD^{\leq n}_Y$ & --- & sheaf of logarithmic differential operators on $Y$ of order $\leq n$ \\
 $\SD^\infty_Y$ & --- & total sheaf of logarithmic differential operators on $Y$
\end{tabular}
\end{quote}
For an ideal $\mathcal{I}$ on $Y$, we also have the following notions:
\begin{quote}
\begin{tabular}{ c c l }
 $\SD^{\leq n}_Y(\mathcal{I})$ or $\SD^{\leq n}(\mathcal{I})$ & --- & ideal on $Y$ generated by the image of $\mathcal{I}$ under $\SD^{\leq n}_Y$ \\
 $\SD^\infty_Y(\mathcal{I})$ or $\SD^\infty(\mathcal{I})$ & --- & ideal on $Y$ generated by the image of $\mathcal{I}$ under $\SD^\infty_Y$ \\
 $\SM(\mathcal{I})$ & --- & monomial saturation of $\mathcal{I}$ \\
 $\logord_y(\mathcal{I})$ & --- & logarithmic order of $\mathcal{I}$ at a point $y \in Y$
\end{tabular}
\end{quote}
These notions (and more) are discussed in Appendix \ref{appB}. In particular, we would also like to bring the reader's attention to the notion of logarithmic coordinates and parameters in (\ref{B108}), as well as (\ref{B109}) and (\ref{B110}). They will play a crucial role in the remainder of this paper.

\subsection{Toroidal centers}\label{3.2}

In this section, we discuss the notion of toroidal centers on a strict toroidal $\kk$-scheme $Y$. These are the ``blow-up centers'' for the resolution algorithm in this paper.

\begin{definition}[Toroidal centers, \textrm{cf. \cite[2.4]{ATW19}}]\label{3201}
Fix a natural number $k \geq 1$, and a nondecreasing sequence \begin{align*}
    (a_1,\dotsc,a_k) \in \QQ_{>0}^{k-1} \times \left(\QQ_{>0} \cup \lbrace \infty \rbrace\right).
\end{align*}
A \emph{toroidal center} $\SJ$ on $Y$, with \emph{invariant} \begin{align*}
    \inv(\SJ) = (a_1,\dotsc,a_k),
\end{align*}
is defined to be an integrally closed Rees algebra on $Y$ (equivalently, an idealistic exponent over $Y$) such that at each point $y$ in $Y$, there exists an (irreducible\footnote{Recall $Y$ is a disjoint union of its irreducible components (\ref{B105}(iii)), so the assertion that $U_y$ is irreducible is equivalent to the assertion that $U_y$ is contained in the component of $Y$ containing $y$.}) open affine neighbourhood $U_y \subset Y$ of $y$ on which either: \begin{enumerate}
    \item[(a)] $\SJ|_{U_y} = \SO_Y[T]|_{U_y}$,
    \item[(b)] or there exist: \begin{enumerate}
    \item[(i)] a choice of logarithmic parameters $\big((x_1^{(y)},\dotsc,x_{n_y}^{(y)}), M_y = \overline{\SM}_{Y,y} \xrightarrow{\beta^{(y)}} H^0(U_y,\SM_Y|_{U_y})\big)$ at $y$, which defines a strict, smooth morphism $U_y \rightarrow \Spec(M_y \rightarrow \kk[M_y \oplus \NN^{n_y}])$ (as in (\ref{B110}(ii)));
    \item[(ii)] if $a_k = \infty$, an nonempty ideal $Q_y$ of $M_y = \overline{\SM}_{Y,y}$, whose image under $\beta^{(y)}$ generates a monomial ideal (\ref{B112}) $\SQ_y$ on $U_y$,
    \end{enumerate}  
    such that \begin{align*}
    \SJ|_{U_y} = \begin{cases}
        \big((x_1^{(y)})^{a_1},\dotsc,(x_k^{(y)})^{a_k}\big) \qquad & \textup{if } a_k \in \QQ_{>0} \\
        \big((x_1^{(y)})^{a_1},\dotsc,(x_{k-1}^{(y)})^{a_{k-1}},\SQ_y^r\big) \qquad & \textup{if } a_k = \infty
    \end{cases}
    \end{align*}
    for some positive rational number $r \in \QQ_{>0}$ independent of $y$. (Note that in particular, $k \leq n_y$ if $a_k \in \QQ_{>0}$, and $k-1 \leq n_y$ if $a_k = \infty$.)
\end{enumerate}
 Given a toroidal center $\SJ$ on $Y$, a choice of data as above for each $y \in Y$ is called a ``\emph{presentation}'' of $\SJ$. We mimic the notation in \cite{ATW19} and record the aforementioned presentation of $\SJ$ as \begin{align*}
    \SJ = \begin{cases}
        (x_1^{a_1},\dotsc,x_k^{a_k}) \qquad & \textup{if } a_k \in \QQ_{>0} \\
        (x_1^{a_1},\dotsc,x_{k-1}^{a_{k-1}},(Q \subset M)^r) \qquad & \textup{if } a_k = \infty.
    \end{cases}
\end{align*}
By the \emph{support}\footnote{Note this is different from, and should not be confused with, the notion of the support of a Rees algebra defined in \cite[Definition 5.1]{Ryd13}.} of a toroidal center $\SJ$, we mean the complement of the Zariski open subset of points $y \in Y$ such that $\SJ_y = \SO_{Y,y}[T]$.

A toroidal center $\SJ^{(y)}$ at a point $y \in Y$, with invariant $\inv(\SJ^{(y)}) = (a_1,\dotsc,a_k)$, is an integrally closed Rees algebra on an open affine neighbourhood $U_y \subset Y$ of $y$ on which we have (i) and (ii) above satisfying the aforementioned properties.
\end{definition}

Observe that we chose to drop the index $y$ in the notation of a toroidal center $\SJ$ on $Y$. This choice of notation would make more sense later: it is justified by the expectation that the resolution algorithm in this paper would be done locally around each $y \in Y$, and patched up afterwards. Some of our results later are written this way, i.e. without making reference to the index $y$ (see for example, \S\ref{4.4}).

It is not immediate that the invariant of a toroidal center is well-defined, i.e. independent of the choice of presentation of $\SJ$. This is the content of the next lemma.

\begin{lemma}\label{3202}
The invariant $\inv(\SJ^{(y)})$ of a toroidal center $\SJ^{(y)}$ at a point $y \in Y$ is independent of choice of presentation for $\SJ^{(y)}$, and hence, is well-defined.
\end{lemma}

We postpone the proof of (\ref{3202}) till \S\ref{4.6}.

\begin{remark}\label{3203}\ 
\begin{enumerate}
    \item Another equivalent definition of a toroidal center on $Y$ (resp. at a point $y \in Y$) is an integrally closed Rees algebra on $Y$ (resp. on an open affine neighbourhood $U_y \subset Y$ of $y$) with a presentation $(x_1^{a_1},\dotsc,x_k^{a_k},(Q \subset M)^r)$ as in (\ref{3201}), but this time allowing $Q$ to be the empty ideal of $M$. In this case, one defines the invariant as $(a_1,\dotsc,a_k,\infty)$ if $Q \neq \emptyset$, and $(a_1,\dotsc,a_k)$ if $Q = \emptyset$.
    \item While the invariant of a toroidal center is well-defined, the positive rational number $r$ appearing in the exponent of $Q$ is evidently not. For example, by replacing $Q$ by $m \cdot Q$ (or $Q^m$ if the monoid $M$ is written multiplicatively), one can replace $r$ by $\frac{r}{m}$. In particular, one can always adjust $Q$ so that $r = \frac{1}{N}$ for some natural number $N \geq 1$.
\end{enumerate}
\end{remark}

\begin{definition}[Reduced toroidal centers]\label{3204}\
\begin{enumerate}
    \item A toroidal center $\SJ$ on $Y$ is \emph{reduced} if the finite entries in $\inv(\SJ)$ are $\frac{1}{n_i}$ for some positive integer $n_i$, and the $\gcd$ of the $n_i$ is $1$.
    \item Given a toroidal center $\SJ$ on $Y$, let $s$ be the unique positive rational number such that $\SJ^s$ is reduced. We denote $\SJ^s$ by $\overline{\SJ}$ and call it the \emph{unique reduced toroidal center associated to $\SJ$}.
\end{enumerate}
One can also define the aforementioned notions for a toroidal center $\SJ^{(y)}$ at a point $y \in Y$.
\end{definition}

Akin to how one can adjust $Q$ in (\ref{3203}(ii)), one can also adjust the $x_i$ appearing in the presentation of a toroidal center, without changing the toroidal center:

\begin{remark}\label{3205}
Let $y \in Y$, and let $\SJ^{(y)} = (x_1^{a_1},\dotsc,x_k^{a_k},(Q \subset M)^r)$ be a toroidal center at $y$, with $k \geq 1$. For each $1 \leq i \leq k$, replace $x_i$ by \begin{align*}
    x_i' = (\lambda_{i,1}x_1+\dotsb+\lambda_{i,i-1}x_{i-1}) + x_i
\end{align*}
where $\lambda_{i,j}$ are sections of $\SO_{Y,y}$. Then we claim that $\SJ^{(y)} = ((x_1')^{a_1},(x_2')^{a_2},\dotsc,(x_k')^{a_k},(Q \subset M)^r)$. While it is possible to prove this from the standpoint of integrally closed Rees algebras, we find it easier to tackle this from the equivalent standpoint of idealistic exponents, where this assertion is reduced to checking the following equality: \begin{align*}
    &\min\big(\lbrace a_i \cdot \nu(x_i) \colon 1 \leq i \leq k \rbrace \cup \lbrace r \cdot \nu(q) \colon q \in Q \rbrace\big) \\
    = \ &\min\big(\lbrace a_i \cdot \nu(x_i') \colon 1 \leq i \leq k \rbrace \cup \lbrace r \cdot \nu(q) \colon q \in Q \rbrace\big).
\end{align*}
More generally, note that one can replace each $x_i$ by \begin{align*}
    x_i' = (\lambda_{i,1}x_1+\dotsb+\lambda_{i,i-1}x_{i-1}) + x_i + (\lambda_{i,i+1}x_{i+1}+\dotsb+\lambda_{i,\ell}x_\ell)
\end{align*}
where $\ell = \max\lbrace 1 \leq j \leq k \colon a_j = a_i \rbrace$, and once again $\lambda_{i,j}$ are sections of $\SO_{Y,y}$. 
\end{remark}

\begin{definition}[Admissibility]\label{3206}
Let $\mathcal{I} \subset \SO_Y$ be an ideal on $Y$, and let $y \in Y$. \begin{enumerate}
    \item A toroidal center $\SJ$ on $Y$ is \emph{$\mathcal{I}$-admissible} if $\SJ$ contains the Rees algebra of $\mathcal{I}$.
    \item A toroidal center $\SJ^{(y)}$ at $y$ is \emph{$\mathcal{I}$-admissible} if, after passing to a smaller affine neighbourhood of $y$ on which $\SJ^{(y)}$ is defined, $\SJ^{(y)}$ contains the Rees algebra of $\mathcal{I}$.
\end{enumerate}
\end{definition}

Note that the support of an $\mathcal{I}$-admissible toroidal center $\SJ$ is always contained in the vanishing locus $V(\mathcal{I})$ of $\mathcal{I}$: indeed, if $y \notin V(\mathcal{I})$, then $\mathcal{I}_y = \SO_{Y,y}$, so $\SJ_y = \SO_{Y,y}[T]$.

Before stating the next lemma, we re-visit (\ref{3201}): there we have that each $U_y$ is a $\kk$-variety, so by \S\ref{2.2}, $\SJ|_{U_y}$ defines an idealistic exponent over $U_y$, which we denote by $\gamma_\SJ^{(y)}$, and refer to it as the idealistic exponent at $y$ associated to $\SJ$ and the affine open neighbourhood $U_y$ of $y$. We can express the notion of admissibility in terms of these idealistic exponents:

\begin{lemma}[Valuative criterion for admissibility]\label{3207}Let the notation be as above. Let $\SJ$ be a toroidal center on $Y$. For a nowhere zero ideal $\mathcal{I}$ on $Y$, the following are equivalent: \begin{enumerate}
    \item $\SJ$ is $\mathcal{I}$-admissible.
    \item For every $y \in Y$ and every open affine neighbourhood $U_y$ of $y$ as in (\ref{3201}), we have $\gamma_{\SJ}^{(y)} \leq \gamma_{\mathcal{I}|_{U_y}}$.
    \item For every $y \in Y$, there exists an open affine neighbourhood $U_y$ of $y$ as in (\ref{3201}) such that $\gamma_{\SJ}^{(y)} \leq \gamma_{\mathcal{I}|_{U_y}}$.
\end{enumerate}
\end{lemma}

\begin{proof}
$\SJ$ contains the Rees algebra of $\mathcal{I}$ if and only if $\SJ$ contains the integral closure of the Rees algebra of $\mathcal{I}$. Fix a choice of open affine neighbourhoods $(U_y)_{y \in Y}$ as in (\ref{3201}). For every $y \in Y$, we deduce from (\ref{2206}) that $\SJ|_{U_y}$ contains the Rees algebra of $\mathcal{I}|_{U_y}$ if and only if $\SJ|_{U_y}$ contains the Rees algebra on $U_y$ associated to $\gamma_{\mathcal{I}|_{U_y}}$. Passing to idealistic exponents over each $U_y$, we see that $\SJ$ is $\mathcal{I}$-admissible if and only if $\gamma_\SJ^{(y)} = \gamma_{\SJ|_{U_y}} \leq \gamma_{\mathcal{I}|_{U_y}}$ for every $y \in Y$.
\end{proof}

Fix a choice of affine open neighbourhoods $(U_y)_{y \in Y}$ as in (\ref{3201}). Then $(\gamma_\SJ^{(y)})_{y \in Y}$ is called the idealistic exponent over $Y$ associated to $\SJ$ and $(U_y)_{y \in Y}$. We will only denote $(\gamma_\SJ^{(y)})_{y \in Y}$ by $\gamma_\SJ$ whenever the discussion does not depend on the choice of $(U_y)_{y \in Y}$. For example:

\begin{notation}\label{3208}We write $\gamma_\SJ \leq \gamma_\mathcal{I}$ to mean either statement (ii) or (iii) in (\ref{3207}). Thus, $\SJ$ is $\mathcal{I}$-admissible if and only if $\gamma_\SJ \leq \gamma_\mathcal{I}$.
\end{notation}

Given a toroidal center $\SJ^{(y)}$ at $y$, let $\widehat{\SJ}^{(y)}$ denote the $\widehat{\SO}_{Y,y}$-subalgebra of $\widehat{\SO}_{Y,y}[T]$ generated by the image of $\SJ^{(y)}$ under \begin{align*}
    \SO_Y[T] \rightarrow \SO_{Y,y}[T] \rightarrow \widehat{\SO}_{Y,y}[T].
\end{align*}
Equivalently, $\widehat{\SJ}^{(y)}$ is the completion $\varprojlim_{k}{\SJ^{(y)}_y/\mathfrak{m}_{Y,y}^k\SJ^{(y)}_y}$, where $\SJ^{(y)}_y$ is the stalk of $\SJ^{(y)}$ at $y$. The next lemma says we can check admissibility by passing to completions:

\begin{lemma}\label{3209}
Let the notation be as above. Let $\SJ^{(y)}$ be a toroidal center at $y$. For an ideal $\mathcal{I}$ on $Y$, $\SJ^{(y)}$ is $\mathcal{I}$-admissible if and only if $\widehat{\SJ}^{(y)}$ is $\widehat{\mathcal{I}}$-admissible.
\end{lemma}

\begin{proof}
Indeed, $\SJ^{(y)}$ is $\mathcal{I}$-admissible if and only if the stalk of $\SJ^{(y)}$ at $y$ contains the Rees algebra of $\mathcal{I}_y$. Since $\SO_{Y,y}[T] \rightarrow \widehat{\SO}_{Y,y}[T]$ is faithfully flat, the latter is equivalent to $\widehat{\SJ}^{(y)}$ being $\widehat{\mathcal{I}}$-admissible \cite[Theorem 7.5]{Mat89}.
\end{proof}

We conclude this section with some easy properties pertaining to admissibility:

\begin{lemma}\label{3210}
Fix a toroidal center $\SJ$ on $Y$, let $\mathcal{I}$ and $\mathcal{I}_j$ be ideals on $Y$, and let $r_j$ be positive rational numbers. Then: \begin{enumerate}
    \item $\SJ$ is $\sum_j{\mathcal{I}_j}$-admissible if and only if $\SJ$ is $\mathcal{I}_j$-admissible for every $j$.
    \item If $\SJ^{r_j}$ is $\mathcal{I}_j$-admissible for every $j$, then $\SJ^{\sum_j{r_j}}$ is $\prod_j{\mathcal{I}_j}$-admissible.
    \item For an integer $\ell \geq 1$, $\SJ$ is $\mathcal{I}$-admissible if and only if $\SJ^\ell$ is $\mathcal{I}^\ell$-admissible.
\end{enumerate}
\end{lemma}

\begin{proof}
Part (i) can be seen directly from (\ref{3206}), and it is easier to deduce part (ii) using the criterion in (\ref{3207}) (after replacing $Y$ with the support of $\mathcal{I}$): if $r_j \cdot \gamma_{\SJ} = \gamma_{\SJ^{r_j}} \leq \gamma_{\mathcal{I}_j}$ for each $j$, then $\sum_j{r_j \cdot \gamma_\SJ} \leq \sum_j{\gamma_{\mathcal{I}_j}} = \gamma_{\prod_j{\mathcal{I}_j}}$. Part (iii) is also clear using (\ref{3207}).
\end{proof}

\begin{lemma}\label{3211}
Let $y \in Y$, and let $\SJ^{(y)} = (x_1^{a_1},\dotsc,x_k^{a_k},(Q \subset M)^r)$ be a toroidal center at $y$, with $k \geq 1$. Let $H$ be the hypersurface $x_1 = 0$ defined on a neighbourhood of $y$ on which $\SJ^{(y)}$ is defined, and let $\mathcal{I}$ be an ideal on $H$. If the restriction of $\SJ^{(y)}$ to $H$, namely $\SJ^{(y)}_H = (x_2^{a_2},\dotsc,x_k^{a_k},(Q \subset M)^r)$, is $\mathcal{I}$-admissible, then $\SJ^{(y)}$ is $(\mathcal{I}\SO_Y)$-admissible.
\end{lemma}

We remark that the ordinary parameters\footnote{(\ref{B108}) (also refer back to (\ref{3201}))} $x_2,\dotsc,x_k$ appearing in the restricted toroidal center $\SJ_H^{(y)}$ are, more precisely, the reduction of $x_2,\dotsc,x_k$ modulo $x_1=0$. Note that if $x_1,\dotsc,x_n$ is a system of ordinary parameters at $y$, then the reduction of $x_2,\dotsc,x_n$ modulo $x_1 = 0$ is a system of ordinary parameters on $H$ at $y$.

\begin{proof}
This can also be verified using (\ref{3206}).
\end{proof}

\section{Weighted toroidal blow-ups}\label{chpt4}

\subsection{Stack-theoretic Proj}\label{4.1}

Let $Y$ be a scheme, or more generally, an algebraic stack, and let $\SR = \bigoplus_{m \in \NN}{\SR_m}$ be a quasi-coherent sheaf of graded $\SO_Y$-algebras on $Y$. In this paper, we will be using the construction $\SProj_Y(\SR)$ in \cite[10.2.7]{Ols16}, called the \emph{stack-theoretic} (or \emph{stacky}) \emph{Proj} of $\SR$ on $Y$. This construction was also recalled briefly in \cite[3.1]{ATW19}, and will be pursued in greater depth in \cite{QR22}. For brevity, we will not repeat the full construction here, but instead recall some of its properties which are relevant to this paper: \begin{enumerate}
    \item When $Y$ is a scheme, $\SProj_Y(\SR)$ is the quotient stack $[(\GSpec_Y(\SR) \setminus S_0)\q\Gm]$, where the grading on $\SR$ defines a $\Gm$-action $(T,s) \mapsto T^m \cdot s$ for $s \in \SR_m$, and the vertex $S_0$ is the closed subscheme defined by the irrelevant ideal $\bigoplus_{m \geq 1}{\SR_m}$ of $\SR$.
    \item When $\SR_1$ is coherent and generates $\SR$ over $\SR_0$, this coincides with the construction in \cite[Chapter II, Section 7, page 160]{Har77}.
    \item $\SProj_Y(\SR)$ carries an invertible sheaf $\SO_{\SProj_Y(\SR)}(1)$ corresponding to the graded $\SO_Y$-algebra $\SR(1)$. If $\SR$ is a Rees algebra on $Y$, then the inclusion $\SO_{\SProj_Y(\SR)}(1) \hookrightarrow \SO_{\SProj_Y(\SR)}$ defines the ideal sheaf of a divisor\footnote{(possibly reducible)} $\SE$ on $\SProj_Y(\SR)$, called the \emph{exceptional divisor} in \cite{QR22}.
    \item When $\SR$ is finitely generated as a $\SO_Y$-algebra with coherent graded components, the resulting morphism $\SProj_Y(\SR) \rightarrow Y$ is proper.
    \item If $f \colon Y' \rightarrow Y$ is a morphism of schemes (or algebraic stacks), $\SProj_{Y'}(f^\ast\SR) = \SProj_Y(\SR) \times_Y Y'$. If $f$ is flat and $\SR$ is a Rees algebra on $Y$, $\SProj_{Y'}(\SR\SO_{Y'}) = \SProj_Y(\SR) \times_Y Y'$.
\end{enumerate}

\subsection{Blow-up along a Rees algebra}\label{4.2}

If $\SR = \bigoplus_{m \in \NN}{\mathcal{I}_m \cdot T^m} \subset \SO_Y[T]$ is a Rees algebra on $Y$, the blow-up $\Bl_Y(\SR)$ of $Y$ along $\SR$ is $\SProj_Y(\SR)$. If $\SR$ is the Rees algebra of an ideal $\mathcal{I} \subset \SO_Y$, $\Bl_Y(\SR)$ is the usual blow-up of $Y$ along the ideal $\mathcal{I}$ (see \cite[page 163]{Har77}).

\subsection{Blow-up along an idealistic exponent}\label{4.3}

Let $\gamma$ be an idealistic exponent over a reduced, separated scheme $Y$ of finite type over $\mathds{k}$, with associated Rees algebra $\SR_\gamma$ on $Y$ (\S\ref{2.2}). The blow-up $\Bl_Y(\gamma)$ of $Y$ along $\gamma$ is defined as $\Bl_Y(\SR_\gamma)$.

\subsection{Weighted toroidal blow-ups: local charts and logarithmic structures}\label{4.4}

From \S\ref{4.4} to \S\ref{4.6}, let $\mathds{k}$ be a field of characteristic zero, and let $Y$ be a strict toroidal $\mathds{k}$-scheme.

Consider a toroidal center $\SJ^{(y)}$ at a point $y \in Y$ of the form $\SJ = (x_1^{1/n_1},\dotsc,x_k^{1/n_k},(Q \subset M)^{1/d})$, where $n_i \geq 1$ and $d \geq 1$ are integers. For this section \S\ref{4.4} only, we replace $Y$ with the open affine neighbourhood of $y$ on which $\SJ^{(y)}$ is defined, and write $\SJ = \SJ^{(y)}$ (so $\SJ$ is now defined on $Y$). Unless otherwise stated, we do not assume $\SJ$ is reduced, and we allow $Q = \emptyset$ (\ref{3203}(i)). As in (\ref{3201}), \begin{enumerate}
    \item[(a)] $x_1,\dotsc,x_k$ is part of a system of ordinary parameters $x_1,\dotsc,x_n$ on $Y$ at $y$ (with $n = \codim_{\fs_y}\overline{\lbrace y \rbrace} \geq k$),
    \item[(b)] $M \rightarrow H^0(Y,\SM_Y)$ is a chart which is neat at $y$,
\end{enumerate}
and together they induce a morphism $Y \rightarrow \Spec(M \rightarrow \mathds{k}[x_1,\dotsc,x_n,M])$ which is strict and smooth of relative dimension $\dim\overline{\lbrace y \rbrace}$ (as in (\ref{B110}(ii))). It is notationally more convenient to identify the ideal $Q \subset M$ with its image of $Q$ in $\SO_Y$, and hence, write $Q$ multiplicatively.

In this section, we study the weighted toroidal blow-up $Y' = \Bl_Y(\SJ) \rightarrow Y$. Since $\SJ$ is the integral closure of the simpler Rees algebra generated by $\lbrace x_1T^{n_1},\dotsc,x_kT^{n_k} \rbrace \cup \lbrace mT^d \colon m \in Q \rbrace$, $Y'$ is covered by the $(x_iT^{n_i})$-charts (for $1 \leq i \leq k$) and the $(mT^d)$-charts (as $m$ varies over a fixed finite set of generators for $Q$). Our first task is to explicate these charts:

\begin{lemma}\label{4401}
The $(x_1T^{n_1})$-chart of $Y'$ is the pullback of the square \begin{equation*}
    \begin{tikzcd}
    & & & \textup{$[U_{x_1} \q \Gmu_{n_1}] = [\Spec(M_{x_1} \rightarrow \mathds{k}[x_2',\dotsc,x_n',M_{x_1}]) \q \Gmu_{n_1}]$} \arrow[to=2-4] \\
    Y \arrow[to=2-4, "\textup{smooth, strict}"] & & & \Spec(M \rightarrow \mathds{k}[x_1,\dotsc,x_n,M])
    \end{tikzcd}
\end{equation*}
where: \begin{enumerate}
    \item $x_1 = u^{n_1}$,
    \item $x_i' = x_i/u^{n_i}$ for $2 \leq i \leq k$,
    \item $x_i' = x_i$ for $i > k$,
    \item $M_{x_1}$ is the saturation of the submonoid of $M \oplus \ZZ \cdot u$ generated by $u$, $M$, and $\lbrace q' = q/u^d \colon q \in Q \rbrace$,
    \item the group $\Gmu_{n_1} = \langle \zeta_{n_1} \rangle$ acts through $\zeta_{n_1} \cdot u = \zeta_{n_1}^{-1}u$, $\zeta_{n_1} \cdot x_i' = \zeta_{n_1}^{n_i}x_i'$ for $2 \leq i \leq k$, trivially on $x_i'$ for $i > k$, and trivially on $M$ (so $\zeta_{n_1} \cdot q' = \zeta_{n_1}^d \cdot q'$ for $q \in Q$).
\end{enumerate}
\end{lemma}

\begin{proof}
Since $Y \rightarrow \Spec(\mathds{k}[x_1,\dotsc,x_n,M])$ is flat, and stacky \emph{Proj} commutes with pullbacks, it suffices to assume $Y = \Spec(\mathds{k}[x_1,\dotsc,x_n,M])$. Set $y_1 = x_1T^{n_1}$. The $y_1$-chart of $Y'$ is the stack $[\Spec(\SJ[y_1^{-1}])\q\Gm]$. By \cite[Lemma 1.3.1]{QR22}, the homomorphism $\SJ[y_1^{-1}] \to \SJ[y_1^{-1}]/(y_1-1)$ of $(\ZZ/n_1\ZZ)$-graded $\SO_Y$-algebras induces an isomorphism of algebraic stacks \begin{equation}\label{eq4.1}
    \left[\Spec\left(\frac{\SJ[y_1^{-1}]}{(y_1-1)}\right)\q\Gmu_{n_1}\right] \xrightarrow{\simeq} [\Spec(\SJ[y_1^{-1}])\q\Gm].
\end{equation}
We sketch the proof presented in \cite[Lemma 1.3.1]{QR22}. On $W := \Spec(\frac{\SJ[y_1^{-1}]}{(y_1-1)}) \times \Gm$, there is the diagonal $\Gmu_{n_1}$-action given by $(y,t) \cdot s = (ys,s^{-1}t)$, and there is also the $\Gm$-action on the second factor given by $(y,t) \cdot s = (y,ts)$. These two actions are free, and commute with each other (and hence, together they induce a free $(\Gmu_{n_1} \times \Gm)$-action on $W$). Then the left hand side of (\ref{eq4.1}) is isomorphic to $[W \q (\Gmu_{n_1} \times \Gm)] = [(W \q \Gmu_{n_1}) \q \Gm]$. One then checks that there is a natural $\Gm$-equivariant isomorphism from $(W \q \Gmu_{n_1})$ to $\Spec(\SJ[y_1^{-1}])$, which yields the desired isomorphism $[(W \q \Gmu_{n_1}) \q \Gm] \xrightarrow{\simeq} [\Spec(\SJ[y_1^{-1}]) \q \Gm]$.

Thus, it remains to show the left hand side of (\ref{eq4.1}) has the desired description. Since $(T^{-1})^{n_1} = y_1^{-1}x_1 \in \SJ[y_1^{-1}]$ and $\SJ[y_1^{-1}]$ is integrally closed in $\SO_Y[T,T^{-1}]$ (\ref{2204}), we see that $T^{-1} \in \SJ[y_1^{-1}]$. Let $u = T^{-1}$. Restricting to $W_1$, we get $u^{n_1} = x_1$, $x_iT^{n_i} = x_i/u^{n_i}$ for $2 \leq i \leq k$, and $qT^d = qu^{-d}$ for every $q \in Q$. Therefore, $\mathds{k}[x_2',\dotsc,x_n',M_{x_1}] \subset \SJ[y_1^{-1}]/(y_1-1)$ is a finite birational extension, and since both are integrally closed in $\SO_Y[T,T^{-1}]$, that inclusion is actually an equality.
\end{proof}

A similar proof explicates the $(mT^d)$-charts of $Y'$. We first fix a notation. Given a (multiplicative) monoid $M$, with an element $m \in M$ and an integer $d > 1$, we write $M[{m^{1/d}}]$ for the pushout of the diagram \begin{equation*}
    \begin{tikzcd}
    \NN \arrow[to=1-3, rightarrow, "1 \mapsto m"] \arrow[to=2-1, rightarrow, swap, "d \cdot"] & & M \arrow[to=2-3, dotted] \\
    \NN \arrow[to=2-3,dotted] & & M[m^{1/d}]
    \end{tikzcd}
\end{equation*}
in the category of monoids, or equivalently, the monoid $M \oplus \NN$ modulo the congruence generated by $(m,0_{\NN}) \thicksim (0_M,d)$. In the lemma below, we shall denote the image of $1$ under the horizontal dotted arrow $\NN \rightarrow M[m^{1/d}]$ as $u$ (so $u^d = m$). Note $M[m^{1/d}]$ may not be torsion-free in general, even if $M$ is torsion-free.

\begin{lemma}\label{4402}
The $(mT^d)$-chart of $Y'$ is the pullback of the square \begin{equation*}
    \begin{tikzcd}
    & & & \textup{$[U_m \q\Gmu_d] = [\Spec(M_m \rightarrow \mathds{k}[x_1',x_2',\dotsc,x_n',M_m]) \q \Gmu_d]$} \arrow[to=2-4] \\
    Y \arrow[to=2-4, "\textup{smooth, strict}"] & & & \Spec(M \rightarrow \mathds{k}[x_1,\dotsc,x_n,M])
    \end{tikzcd}
\end{equation*}
where: \begin{enumerate}
    \item $M_m = \overline{M}'$, and $M'$ is the saturation of the submonoid of $M[{m^{1/d}}]^{\gp}$ generated by $M[{m^{1/d}}]$ and $\lbrace q' = q/m = q/u^d \colon q \in Q \rbrace$,
    \item $x_i' = x_i/u^{n_i}$ for $1 \leq i \leq k$,
    \item $x_i' = x_i$ for $i > k$,
    \item the group $\Gmu_d = \langle \zeta_d \rangle$ acts through $\zeta_d \cdot u = \zeta_d^{-1}u$, $\zeta_d \cdot x_i' = \zeta_d^{n_i}x_i'$ for $1 \leq i \leq k$, trivially on $x_i'$ for $i > k$, and trivially on $M$ (so $\zeta_d$ also acts trivially on $q'$, for $q \in Q$). \qed
\end{enumerate}
\end{lemma}

Together, (\ref{4401}) and (\ref{4402}) present a natural choice of an \'etale cover $U$ of $Y' = \Bl_Y(\SJ)$. Each $(x_iT^{n_i})$-chart of $\Bl_Y(\SJ)$ admits an \'etale cover from the pullback of $U_{x_i}$ to $Y$, and each $(mT^d)$-chart of $\Bl_Y(\SJ)$ admits an \'etale cover from the pullback of $U_m$ to $Y$. For the remainder of this paper, \begin{quote}
    $U$ denotes the disjoint union of the pullbacks of $U_{x_i}$ and $U_m$ to $Y$ (where $1 \leq i \leq k$, and $m$ varies over a fixed finite set of generators for $Q$).
\end{quote}
Note that the composition $U \rightarrow Y' = \Bl_Y(\SJ) \rightarrow Y$ is an alteration. In addition, the principal ideal $E = (u)$ on $U$ descends to give the exceptional ideal $\SE$ on $Y'=\Bl_Y(\SJ)$.

In (\ref{4401}) and (\ref{4402}), we have also specified logarithmic structures on the $U_{x_i}$ and the $U_m$, such that the exceptional ideal $E = (u)$ is encoded in the logarithmic structures (this should be compared to (\ref{1101}(iv))). These pull back, via the strict morphism $Y \rightarrow \Spec(M \rightarrow k[x_1,\dotsc,x_n,M])$, to define a logarithmic structure on $U$, which manifests $U$ as a strict toroidal $\mathds{k}$-scheme. 

Finally, the logarithmic structure on $U$ descends to a logarithmic structure on $Y' = \Bl_Y(\SJ)$ (see Appendix \ref{B.2}). Observe (from the charts) that the \'etale cover $U$ is a strict toroidal $\mathds{k}$-scheme, whence $Y'$ is a \emph{toroidal Deligne-Mumford stack over $\mathds{k}$} (see (\ref{B201})). If $k=0$, observe too that the morphism $Y' \rightarrow Y$ is logarithmically smooth (because $U \twoheadrightarrow Y' \rightarrow Y$ is logarithmically smooth). This is not true if $k \geq 1$.

\begin{lemma}\label{4403}
Let the notation be as above, and let $\gamma = \gamma_\SJ$ be the idealistic exponent over $Y$ associated to $\SJ$. Then $\gamma\SO_U$ is the idealistic exponent over $U$ associated to the exceptional ideal $E$, i.e. $\gamma\SO_U = \gamma_E$.
\end{lemma}

\begin{proof}
This is a simple computation: for example, over $U_{x_1}$, we have for every $\nu \in \ZR(U_{x_1})$, \begin{enumerate}
    \item $\nu(u) = \frac{1}{n_1} \cdot \nu(x_1)$.
    \item $\nu(u) = \frac{1}{n_i} \cdot \nu(x_i) - \frac{1}{n_i} \cdot \nu(x_i') \leq \frac{1}{n_i}\nu(x_i)$ for $2 \leq i \leq k$.
    \item $\nu(u) = \frac{1}{d} \cdot \nu(q) - \frac{1}{d}\nu(q') \leq \frac{1}{d} \cdot \nu(q)$ for $q \in Q$.
\end{enumerate}
Therefore, $\min\big(\lbrace \frac{1}{n_i} \cdot \nu(x_i) \colon 1 \leq i \leq k \rbrace \cup \lbrace \frac{1}{d} \cdot \nu(q) \colon q \in Q \rbrace\big) = \nu(u)$. This computation persists in the other $U_{x_i}$ and $U_m$.
\end{proof}

\begin{proposition}\label{4404}
Let $\mathcal{I} \subset \SO_Y$ be a nowhere zero ideal on $Y$, and $\SJ = (x_1^{1/n_1},\dotsc,x_k^{1/n_k},(Q \subset M)^{1/d})$ be a toroidal center on $Y$ (as in this section), where $n_i,d \geq 1$ are integers. Let $\SE$ be the exceptional ideal of the weighted toroidal blow-up $Y' = \Bl_Y(\SJ) \rightarrow Y$. \begin{enumerate}
    \item If $\ell \geq 1$ is an integer such that $\SJ^\ell$ is $\mathcal{I}$-admissible, then $\mathcal{I}\SO_{Y'}$ factors as $\SE^\ell \cdot \mathcal{I}'$ for some ideal $\mathcal{I}'$ on $\SO_{Y'}$. 
    \item The converse holds as well: if $\mathcal{I}\SO_{Y'}$ factors as $\SE^\ell \cdot \mathcal{I}'$ for some ideal $\mathcal{I}'$ on $\SO_{Y'}$ and some integer $\ell \geq 1$, then $\SJ^\ell$ is $\mathcal{I}$-admissible.
\end{enumerate}
\end{proposition}

\begin{proof}
Let $U$ be the \'etale cover of $Y'$ defined earlier, with principal ideal $E=(u)$ on $U$. For (i), use (\ref{4403}): $\gamma_\SJ\SO_U = \gamma_E$, so $\gamma_{\SJ^\ell}\SO_U = \gamma_{E^\ell}$. But $\SJ^\ell$ is $\mathcal{I}$-admissible, so $\gamma_{\SJ^\ell} \leq \gamma_\mathcal{I}$, whence $\gamma_{\mathcal{I}\SO_U} = \gamma_\mathcal{I}\SO_U \geq \gamma_{\SJ^\ell}\SO_U = \gamma_{E^\ell}$. But $U$ is normal (\ref{B105}(iii)), so (\ref{A201}), coupled with the inequality $\gamma_{\mathcal{I}\SO_U} \geq \gamma_{E^\ell}$, implies that the fractional ideal $E^{-\ell}(\mathcal{I}\SO_U)$ is an ideal $I'$ on $\SO_U$. Moreover, since $E$ is principal, $\mathcal{I}\SO_U = E^\ell \cdot I'$. By descent, we get $\mathcal{I}\SO_{Y'} = \SE^\ell \cdot \mathcal{I}'$ for some ideal $\mathcal{I}'$ on $\SO_{Y'}$. 

For (ii), the hypothesis says $\gamma_\mathcal{I}\SO_U = \gamma_{\mathcal{I}\SO_U} \geq \gamma_{E^\ell} = \gamma_{\SJ^\ell}\SO_U$. Pulling idealistic exponents back to $\SO_U$ is order-preserving, whence $\gamma_\mathcal{I} \geq \gamma_{\SJ^\ell}$. Thus, $\SJ^\ell$ is $\mathcal{I}$-admissible.
\end{proof}

\begin{definition}\label{4405}
Take $\ell = \max\lbrace n \in \NN \colon \SJ^n$ is $\mathcal{I}$-admissible$\rbrace$ in (\ref{4404}(i)). The corresponding ideal $\mathcal{I}'$ is called the weak (or birational) transform of $\mathcal{I}$ under the weighted toroidal blow-up $Y' = \Bl_Y(\SJ) \rightarrow Y$.
\end{definition}

By considering the charts in (\ref{4401}) and (\ref{4402}), we get part (i) of the following lemma:

\begin{lemma}\label{4406}
Let $\SJ = (x_1^{1/n_1},\dotsc,x_k^{1/n_k},(Q \subset M)^{1/d})$ be a toroidal center on $Y$ (as in this section), where $n_i,d \geq 1$ are integers. Fix a natural number $c \geq 1$, and set $\widetilde{\SJ} = \SJ^{1/c} = (x_1^{1/cn_1},\dotsc,x_k^{1/cn_k},(Q \subset M)^{1/cd})$. \begin{enumerate}
    \item If $Y' \rightarrow Y$ and $\widetilde{Y}' \rightarrow Y$ are weighted toroidal blow-ups corresponding to $\SJ$ and $\widetilde{\SJ}$, with respective exceptional ideals $\SE$ and $\widetilde{\SE}$, then $\widetilde{Y}' = Y'(\sqrt[c]{\SE})$ is the root stack of $Y'$ along $\SE$.
    \item Assume $k \geq 1$. Write $H$ for the hypersurface $x_1 = 0$ on $Y$, and let $\overline{H}' \rightarrow H$ be the weighted toroidal blow-up along the reduced toroidal center $\overline{\SJ}_H$ associated to the restricted toroidal center $\SJ_H = (x_2^{1/n_2},\dotsc,x_k^{1/n_k},\SQ^{1/d})$, with exceptional ideal $\overline{\SE}_H$. Then the proper transform $\widetilde{H}' \rightarrow H$ of $H$ via the weighted toroidal blow-up along $\widetilde{\SJ}$ is the root stack $\overline{H}'\Big(\sqrt[(cc')]{\overline{\SE}_H}\Big)$ of $\overline{H}'$ along $\overline{\SE}_H \subset \overline{H}'$, where $c' = \gcd(n_2,\dotsc,n_k)$. In other words, $\widetilde{H}' \rightarrow H$ is the weighted toroidal blow-up along $\overline{\SJ}_H^{1/(cc')}$.
\end{enumerate}
\end{lemma}

\begin{proof}[Proof of \emph{(ii)}]
Let $H' \rightarrow H$ be the weighted toroidal blow-up along $\SJ_H$, with exceptional ideal $\SE_H$. By (i), we have $H' = \overline{H}'\Big(\sqrt[c']{\overline{\SE}_H}\Big)$. Next, note that $H' \rightarrow H$ coincides with the proper transform of $H$ via the weighted toroidal blow-up along $\SJ$ --- this can be seen from the charts in (\ref{4401}) and (\ref{4402}). Now apply (i) again: it says $\widetilde{H}' = H'(\sqrt[c]{\SE_H})$. Combining our observations, we get $\widetilde{H}' = \overline{H}'\Big(\sqrt[(cc')]{\overline{\SE}_H}\Big)$, as desired.
\end{proof}

\subsection{Admissibility of toroidal centers: further results}\label{4.5}

(\ref{4404}) provides a convenient method to verify more intricate results on admissibility of toroidal centers. Before doing so, we state a key lemma:

\begin{lemma}\label{4501}
Let $\SJ = (x_1^{1/n_1},\dotsc,x_k^{1/n_k},(Q \subset M)^{1/d})$ be a toroidal center on $Y$, where $k \geq 1$, and $n_i,d \geq 1$ are integers. Let $\SE$ be the exceptional ideal of the weighted toroidal blow-up $Y' = \Bl_Y(\SJ) \rightarrow Y$. For an ideal $\mathcal{I}$ on $Y$, we have $\SD^{\leq 1}_Y(\mathcal{I})\SO_{Y'} \subset \SE^{-n_1} \cdot \SD^{\leq 1}_{Y'}(\mathcal{I}\SO_{Y'})$.
\end{lemma}

\begin{proof}
We can check the lemma over a point $y \in Y$, and hence, it suffices to assume $\SJ$ is a toroidal center at a fixed $y \in Y$ (as in the beginning of \S\ref{4.4}). We shall also work on the \'etale cover $U$ of $Y' = \Bl_Y(\SJ)$ defined before (\ref{4403}), where $\SE$ pulls back to the principal ideal $E=(u)$ on $U$. We can also pass to completion at $y$, i.e. work in $\widehat{\SO}_{Y,y} \simeq \kappa\llbracket x_1,\dotsc,x_n,M \rrbracket$, where $x_1,\dotsc,x_n$ are ordinary parameters at $y$, $M = \overline{\SM}_{Y,y}$, and $\kappa = \kappa(y)$ is the residue field at $y$. Extend $x_1,\dotsc,x_n$ to ordinary coordinates $x_1,\dotsc,x_N$ (\ref{B108}) at $y$, and fix a basis $m_1,\dots,m_r \in M$ for $M^{\gp}$. Let $u_i = \exp(m_i)$ (\ref{B108}) for $1 \leq i \leq r$. By (\ref{B109}), $\SD^1_{Y,y}$ admits a basis given by $\frac{\partial}{\partial x_1},\dotsc,\frac{\partial}{\partial x_N},u_1\frac{\partial}{\partial u_1},\dotsc,u_r\frac{\partial}{\partial u_r}$. For a point $y'$ in the $(x_1T^{n_1})$-chart $U_{x_1}$ over $y$ (\ref{4401}), the same lemma says $\SD^1_{U_{x_1},y'}$ admits a basis given by $\frac{\partial}{\partial x_2'},\dotsc,\frac{\partial}{\partial x_k'},\frac{\partial}{\partial x_{k+1}},\dotsc,\frac{\partial}{\partial x_N},u\frac{\partial}{\partial u},u_1\frac{\partial}{\partial u_1},\dotsc,u_r\frac{\partial}{\partial u_r}$ (where $x_i = u^{n_i}x_i'$ for $2 \leq i \leq k$). For $f = f(x_1,\dotsc,x_n,u_1,\dotsc,u_r) \in \widehat{\mathcal{I}} \subset \widehat{\SO}_{Y,y}$, we compute, on $U_{x_1}$, the following equations: \begin{enumerate}
    \item For $2 \leq i \leq k$, \begin{align*}
    &\frac{\partial}{\partial x_i'}\big(f(u^{n_1},u^{n_2}x_2',\dotsc,u^{n_k}x_k',x_{k+1}',\dotsc,x_n',u_1,\dotsc,u_r)\big) \\
    = \ &\frac{\partial f}{\partial x_i}(u^{n_1},u^{n_2}x_2',\dotsc,u^{n_k}x_k',x_{k+1}'\dotsc,x_n',u_1,\dotsc,u_r) \cdot u^{n_i}.
    \end{align*}
    \item $\big(u\frac{\partial}{\partial u}\big)\big(f(u^{n_1},u^{n_2}x_2',\dotsc,u^{n_k}x_k',x_{k+1}',\dotsc,x_n')\big) = \frac{\partial f}{\partial x_1}(u^{n_1},u^{n_2}x_2',\dotsc,u^{n_k}x_k',x_{k+1}',\dotsc,x_n') \cdot (n_1 u^{n_1}) \; + \sum_{i=2}^k{\frac{\partial f}{\partial x_i}(u^{n_1},u^{n_2}x_2',\dotsc,u^{n_k}x_k',x_{k+1}',\dotsc,x_n') \cdot (n_i u^{n_i}x_i')}$.
\end{enumerate}
Rewriting the equation in (ii), we get: \begin{align*}
    &\frac{\partial f}{\partial x_1}(u^{n_1},u^{n_2}x_2',\dotsc,u^{n_k}x_k',x_{k+1}',\dotsc,x_n') \\
    = \ &\frac{1}{n_1} \cdot u^{-n_1}\Bigg(\big(u\frac{\partial}{\partial u}\big)\big(f(u^{n_1},u^{n_2}x_2',\dotsc,u^{n_k}x_k',x_{k+1}',\dotsc,x_n')\big) \\
    &\qquad \qquad \qquad \qquad \qquad \qquad -\sum_{i=2}^k{n_i \cdot x_i' \cdot \frac{\partial}{\partial x_i'}\big(f(u^{n_1},u^{n_2}x_2',\dotsc,u^{n_k}x_k',x_{k+1}',\dotsc,x_n')\big)}\Bigg).
\end{align*}
Since $n_1 \geq n_2 \geq \dotsb \geq n_k$, these equations suffice to show the lemma on the $(x_1T^{n_1})$-chart. This computation persists for points $y'$ over $y$ in the remaining $(x_iT^{n_i})$-charts, as well as the $(mT^d)$-charts (\ref{4402}).
\end{proof}

The key proposition in this section is:

\begin{proposition}\label{4502}Let $\SJ = (x_1^{a_1},\dotsc,x_k^{a_k},(Q \subset M)^r)$ be a toroidal center on $Y$, where $k \geq 1$. Let $\mathcal{I}$ be a nowhere zero ideal on $Y$. If $\SJ$ is $\mathcal{I}$-admissible, then we have the following statements: \begin{enumerate}
    \item If $a_1 \geq 1$, then $\SJ^{\frac{a_1-1}{a_1}}$ is $\SD^{\leq 1}(\mathcal{I})$-admissible.
    \item $\SJ^{\frac{a_1+1}{a_1}}$ is $(x_1\mathcal{I})$-admissible.
\end{enumerate}
\end{proposition}

\begin{proof}
Before delving into the proof, let us fix some notation. Let $N$ be a natural number such that $\widetilde{\SJ} = \SJ^{1/N}$ is of the form $(x_1^{1/n_1},\dotsc,x_k^{1/n_k},(Q \subset M)^{r/N})$ for positive integers $n_i$. Note that in particular, $N=a_1n_1$. By replacing $Q$ by some multiple $m \cdot Q$ (or $Q^m$ if the monoid is written multiplicatively), we may assume $\frac{r}{N} = \frac{1}{d}$ for some integer $d \geq 1$. Let $Y' \rightarrow Y$ be the weighted toroidal blow-up along $\widetilde{\SJ}$, with exceptional ideal $\SE$. By (\ref{4404}(i)), since $\widetilde{\SJ}^N = \SJ$ is $\mathcal{I}$-admissible, $\mathcal{I}\SO_{Y'}$ factors as $\SE^N \cdot \mathcal{I}' = \SE^{a_1n_1} \cdot \mathcal{I}'$ for some ideal $\mathcal{I}'$ on $\SO_{Y'}$. If we can show that $\SD^{\leq 1}(\mathcal{I})\SO_{Y'} = \SE^{(a_1-1)n_1} \cdot \mathcal{I}_1$ for some ideal $\mathcal{I}_1$ on $\SO_{Y'}$, part (i) follows from (\ref{4404}(ii)). This is the $j=1$ case of the following lemma:

\begin{lemma}\label{4503}
Assume the hypotheses of \emph{(\ref{4502})}, and adopt the set-up above. Then for every integer $1 \leq j \leq a_1$, $\SD^{\leq j}_Y(\mathcal{I})\SO_{Y'}$ factors as $\SE^{(a_1-j)n_1} \cdot  \mathcal{I}_j$ for some ideal $\mathcal{I}_j$ on $Y'$, with $\mathcal{I}_j \subset \SD^{\leq j}_{Y'}(\mathcal{I}')$.
\end{lemma}

\begin{proof}[Proof of \emph{(\ref{4503})} for the case $j=1$] 
We use the product rule to obtain: \begin{align*}
    \SD^{\leq 1}_{Y'}(\mathcal{I}\SO_{Y'}) = \SD^{\leq 1}_{Y'}(\SE^N \cdot \mathcal{I}') &\subset \SD^{\leq 1}_{Y'}(\SE^N) \cdot \mathcal{I}'+ \SE^N \cdot \SD^{\leq 1}_{Y'}(\mathcal{I}') \\
    &= \SE^N \cdot \mathcal{I}' + \SE^N \cdot \SD^{\leq 1}_{Y'}(\mathcal{I}') = \SE^N \cdot \SD^{\leq 1}_{Y'}(\mathcal{I}').
\end{align*}
Next, (\ref{4501}) says $\SD^{\leq 1}_Y(\mathcal{I})\SO_{Y'} \subset \SE^{-n_1} \cdot \SD^{\leq 1}_{Y'}(\mathcal{I}\SO_{Y'})$. Combining this with the above computation, $\SD^{\leq 1}_Y(\mathcal{I})\SO_{Y'} \subset \SE^{N-n_1} \cdot \SD^{\leq 1}_{Y'}(\mathcal{I}') = \SE^{(a_1-1)n_1} \cdot \SD^{\leq 1}_{Y'}(\mathcal{I}')$. The fractional ideal $\SE^{-(a_1-1)n_1} \cdot \SD^{\leq 1}_Y(\mathcal{I})\SO_{Y'}$ is contained in $\SD^{\leq 1}_{Y'}(\mathcal{I}')$, and hence, is an ideal $\mathcal{I}_1$ on $Y'$. Since $\SE$ is a principal ideal, we get the desired factorization $\SD^{\leq 1}_Y(\mathcal{I})\SO_{Y'} = \SE^{(a_1-1)n_1} \cdot \mathcal{I}_1$, with $\mathcal{I}_1 \subset \SD^{\leq 1}_{Y'}(\mathcal{I}')$. 
\end{proof}

For part (ii) of (\ref{4502}), adopt the set-up at the beginning of this proof. Then $x_1$ factors as $u^{n_1} \cdot x_1'$ in $\SO_{Y'}$ (cf. (\ref{4401}) and (\ref{4402})), whence \begin{align*}
    (x_1\mathcal{I})\SO_{Y'} = \SE^{(a_1+1)n_1} \cdot (x_1'\mathcal{I}').
\end{align*}
Once again, an application of (\ref{4404}(ii)) completes the proof.
\end{proof}

We can now prove (\ref{4503}) in general:

\begin{proof}[Proof of \emph{(\ref{4503})}]
We have already shown the case $j=1$. In general, induct on $j$. Assume (\ref{4503}) is known for some $1 \leq j < a_1$, and we prove the lemma for $j+1$. Applying (\ref{4502}(i)) repeatedly, $\SJ^{\frac{a_1-j}{a_1}}$ is $\SD^{\leq j}(\mathcal{I})$-admissible, and induction hypothesis says $\SD^{\leq j}(\mathcal{I})\SO_{Y'} = \SE^{(a_1-j)n_1} \cdot \mathcal{I}_j$ for an ideal $\mathcal{I}_j$ on $Y'$, with $\mathcal{I}_j \subset \SD^{\leq j}(\mathcal{I}')$. Applying the case $j=1$ with \begin{enumerate}
    \item[(a)] $\SJ$ replaced by $\SJ^{\frac{a_1-j}{a_1}} = \overline{\SJ}^{n_1(a_1-j)}$;
    \item[(b)] $\mathcal{I}$ replaced by $\SD^{\leq j}(\mathcal{I})$, 
\end{enumerate}
we see that $\SD^{\leq j+1}(\mathcal{I})\SO_{Y'} = \SD^{\leq 1}(\SD^{\leq j}(\mathcal{I}))\SO_{Y'}$ factors as $\SE^{(a_1-j-1)n_1} \cdot \mathcal{I}_{j+1}$ for some ideal $\mathcal{I}_{j+1}$ on $Y'$, with $\mathcal{I}_{j+1} \subset \SD^{\leq 1}(\mathcal{I}_j) \subset \SD^{\leq 1}(\SD^{\leq j}(\mathcal{I}')) = \SD^{\leq j+1}(\mathcal{I}')$, as desired.
\end{proof}

(\ref{4502}(i)) provides us with the first piece of information about the invariant of an $\mathcal{I}$-admissible toroidal center at a point $y \in Y$:

\begin{corollary}\label{4504}
Let $\mathcal{I}$ be an ideal on $Y$, and fix $y \in Y$ such that $\mathcal{I}_y \neq 0$. If $\logord_y(\mathcal{I}) = b_1 < \infty$, and $\SJ^{(y)} = (x_1^{a_1},\dotsc,x_k^{a_k},(Q \subset M)^r)$ is a $\mathcal{I}$-admissible toroidal center at $y$ with $k \geq 1$, then $a_1 \leq b_1$.
\end{corollary}

\begin{proof}
Suppose for a contradiction that $a_1 > b_1$. We apply (\ref{4502}(i)) repeatedly to conclude that $ (\SJ^{(y)})^{\frac{a_1-b_1}{a_1}}$ is $\SD^{\leq b_1}(\mathcal{I})$-admissible. By restricting to a smaller open affine neighbourhood $U_y$ of $y$ on which $\SJ^{(y)}$ is defined, we may arrange for $\SD^{\leq b_1}(\mathcal{I})$ to be $\SO_Y$ when restricted to $U_y$. Replacing $Y$ with $U_y$, $(\SJ^{(y)})^{\frac{a_1-b_1}{a_1}}$ is $\SO_Y$-admissible, but \begin{align*}
    (\SJ^{(y)})^{\frac{a_1-b_1}{a_1}} = \big(x_1^{a_1-b_1},\dotsc,x_k^{\frac{(a_1-b_1)a_k}{a_1}},(Q \subset M)^{\frac{(a_1-b_1)r}{a_1}}\big),
\end{align*}
with $a_1-b_1 > 0$, a contradiction.
\end{proof}

\subsection{The invariant of a toroidal center is well-defined}\label{4.6}

We prove (\ref{3202}) in this section. The main ingredient of the proof is (\ref{4504}), but we will need two lemmas.

\begin{lemma}\label{4601}
Let $y \in Y$, let $\kappa(y)$ be the residue field at $y$, and let $\fs_y$ be the logarithmic stratum of $Y$ at $y$ \emph{(\ref{B103})}. Let $\SJ^{(y)} = (x_1^{a_1},\dotsc,x_k^{a_k},(Q \subset M)^r)$ be a toroidal center at $y$, with $k \geq 1$. For a homogeneous section $fT^\ell \in \SJ^{(y)}$, write the image of $f$ under $\SO_{Y,y} \twoheadrightarrow \SO_{\fs_y,y} \to \widehat{\SO}_{\fs_y,y} \simeq \kappa(y)\llbracket x_1,\dotsc,x_n \rrbracket$ as $\sum_{\vec{\alpha}}{c_{\vec{\alpha}} \cdot x_1^{\alpha_1}\dotsm x_n^{\alpha_n}}$ for some $c_{\vec{\alpha}} \in \kappa(y)$. Then $\sum_{i=1}^k{\frac{\alpha_i}{a_i}} \geq \ell$ whenever $c_{\vec{\alpha}} \neq 0$.
\end{lemma}

\begin{proof}
We may replace $Y$ by $\fs_y$, and reduce to the case where $Y$ is a smooth $\kk$-variety with trivial logarithmic structure, and $\SJ^{(y)} = (x_1^{a_1},\dotsc,x_k^{a_k})$ with $k \geq 1$. Replacing $\SJ^{(y)}$ by $(\SJ^{(y)})^\ell$, we may assume $\ell = 1$. This is the same exact situation as in \cite[5.2.1]{ATW20a}. We recall its proof. Consider the following valuation in $\ZR\left(\widehat{\SO}_{Y,y}\right)$: \[
    \nu_\SJ\left(\sum_{\vec{\alpha}}{c_{\vec{\alpha}} \cdot x_1^{\alpha_1}\dotsm x_n^{\alpha_n}}\right) := \min_{\substack{\vec{\alpha} \\ c_{\vec{\alpha}} \neq 0}}\left(\sum_{i=1}^k{\frac{\alpha_i}{a_i}}\right).
\]
The hypothesis implies $\SJ^{(y)}$ is $(f)$-admissible, and hence, by (\ref{3207}), $\gamma_{\SJ^{(y)}} \leq \gamma_{(f)}$. Therefore, \[
    \min_{\substack{\vec{\alpha} \\ c_{\vec{\alpha}} \neq 0}}\left(\sum_{i=1}^k{\frac{\alpha_i}{a_i}}\right) = \nu_{\SJ}(f) = \gamma_{(f),\nu_{\SJ}} \geq \gamma_{\SJ^{(y)},\nu_{\SJ}} = \min\lbrace a_i \cdot \nu_{\SJ}(x_i) \colon 1 \leq i \leq k \rbrace = 1.
\]
This completes the proof.
\end{proof}

\begin{lemma}[Exchange]\label{4602}
Let $y \in Y$, and let $\SJ^{(y)} = (x_1^{a_1},\dotsc,x_k^{a_k},(Q \subset M)^r)$ be a toroidal center at $y$, with $k \geq 1$. Suppose that $x_1',x_2,\dotsc,x_n$ is also a system of ordinary parameters at $y$, and suppose $\SJ^{(y)}$ is $((x_1')^{a_1})$-admissible. After possibly passing to a smaller affine neighbourhood of $y$ on which $\SJ^{(y)}$ is defined, we have $\SJ^{(y)} = ((x_1')^{a_1},x_2^{a_2},\dotsc,x_k^{a_k},(Q \subset M)^r)$.
\end{lemma}

\begin{proof}
The hypothesis says that $\SJ^{(y)}$ contains $(\SJ')^{(y)} = ((x_1')^{a_1},x_2^{a_2},\dotsc,x_n^{a_n},(Q \subset M)^r)$. This is necessarily an equality near $y$, by passing to completion at $y$, and seeing that the $\kappa(y)$-dimensions of each $T^N$-graded piece on both sides match.
\end{proof}

We can now prove (\ref{3202}):

\begin{proof}[Proof of \emph{(\ref{3202})}]
Suppose that $\SJ^{(y)}$ admits the following presentations: \begin{align*}
    (x_1^{a_1},\dotsc,x_k^{a_k},(Q \subset M)^r) = \SJ^{(y)} = ((x_1')^{b_1},\dotsc,(x_\ell')^{b_\ell},(Q' \subset M)^s).
\end{align*}
Note that $k = 0$ if and only if $\ell = 0$, in which case $\inv(\SJ^{(y)}) = (\infty)$. Henceforth, assume $k \geq 1$, and hence $\ell \geq 1$. By replacing $\SJ^{(y)}$ by some power of itself, we may assume $a_1$ and $b_1$ are integers. Observe that in particular, $(x_1^{a_1},\dotsc,x_k^{a_k},(Q \subset M)^r)$ is $((x_1')^{b_1})$-admissible. Using (\ref{4504}), we see that $a_1 \leq b_1$. Reversing the roles, we get $b_1 \leq a_1$, whence $a_1 = b_1$. Applying (\ref{4502}(i)) repeatedly, we see that $(\SJ^{(y)})^{1/a_1} = (x_1,x_2^{a_2/a_1},\dotsc,x_k^{a_k/a_1},(Q \subset M)^{r/a_1})$ is $(x_1')$-admissible. Extending $x_1,\dotsc,x_k$ to a system of ordinary parameters $x_1,\dotsc,x_n$ at $y$, and passing to completion at $y$, write the image of $x_1'$ under $\SO_{Y,y} \twoheadrightarrow \SO_{\fs_y,y} \to \widehat{\SO}_{\fs_y,y} \simeq \kappa(y)\llbracket x_1,\dotsc,x_n \rrbracket$ as $\sum_{\vec{\alpha}}{c_{\vec{\alpha}}x_1^{\alpha_1}\dotsm x_n^{\alpha_n}}$ for some $c_{\vec{\alpha}} \in \kappa(y)$. Applying (\ref{4601}), we see that whenever $c_{\vec{\alpha}} \neq 0$, then $\alpha_1+\sum_{i=2}^k{\frac{\alpha_i}{a_i/a_1}} \geq 1$. Consequently, if we let $k_0 = \max\lbrace 1 \leq i \leq k \colon a_i = a_1 \rbrace \geq 1$, the image of $x_1'$ in $\SO_{\fs_y,y}$ lies in $(x_1,\dotsc,x_{k_0}) + \mathfrak{m}_{\fs_y,y}^2$, where $\mathfrak{m}_{\fs_y,y}$ is the maximal ideal of $\SO_{\fs_y,y}$. Therefore, after possibly reordering $x_1,\dotsc,x_{k_0}$, we may replace $x_1$ by $x_1'$ so that $(x_1',x_2,\dotsc,x_n)$ is a system of ordinary parameters at $y$. It is essential to note that the reordering does not mess up the presentation of $\SJ^{(y)} = (x_1^{a_1},\dotsc,x_k^{a_k},(Q \subset M)^r)$, since $a_1 = \dotsb = a_{k_0}$. Applying (\ref{4602}), we obtain: \begin{align*}
    ((x_1')^{a_1},x_2^{a_2},\dotsc,x_k^{a_k},(Q \subset M)^r) = \SJ^{(y)} = ((x_1')^{a_1},(x_2')^{b_2},\dotsc,(x_\ell')^{b_\ell},(Q' \subset M)^s).
\end{align*}
Now restrict to the hypersurface $H$ containing $y$ given by $x_1' = 0$, and we get: \begin{align*}
    (x_2^{a_2},\dotsc,x_k^{a_k},(Q \subset M)^r) = \SJ^{(y)}_H = ((x_2')^{b_2},\dotsc,(x_\ell')^{b_\ell},(Q' \subset M)^s).
\end{align*}
By induction hypothesis, we obtain $k=\ell$, and $a_i=b_i$ for $2 \leq i \leq k = \ell$, as desired. Moreover, $Q \neq \emptyset$ if and only if $Q' \neq \emptyset$. However, we remind the reader (\ref{3203}(ii)) that $Q$ may be different from $Q'$, and $r$ may be not equal to $s$.
\end{proof}

\section{Coefficient ideals}\label{chpt5}

As always, assume $Y$ is a strict toroidal $\kk$-scheme. Fix an ideal $\mathcal{I} \subset \SO_Y$. 

\subsection{Maximal contact element}\label{5.1}
In this section, we assume $1 \leq a = \max\logord(\mathcal{I}) < \infty$. Following \cite[Definition 3.79]{Kol07}, the \emph{maximal contact ideal} of $\mathcal{I}$ is defined as: \begin{align*}
    \MC(\mathcal{I}) = \SD^{\leq a-1}(\mathcal{I}).
\end{align*}
For a point $y \in Y$ with $a = \logord_y(\mathcal{I})$, a \emph{maximal contact element} of $\mathcal{I}$ at $y$ is a section of $\MC(\mathcal{I})$ over a neighbourhood of $y$ in $Y$, which can be extended to a system of ordinary parameters at $y$ (or equivalently, has logarithmic order $1$ at $y$). Maximal contact elements at such points $y \in Y$ always exist, because we are in characteristic zero. The vanishing locus of a maximal contact element of $\mathcal{I}$ at $y$ is called a \emph{hypersurface of maximal contact} for $\mathcal{I}$ through $y$. It is well-known that hypersurfaces of maximal contact play a crucial role in resolution of singularities in characteristic zero, in the sense that they allow for induction on dimension: namely, one passes to a hypersurface of maximal contact in the induction step.

Following \cite[Definition 3.53]{Kol07}, we say $\mathcal{I} \subset \SO_Y$ is \emph{MC-invariant}, if \begin{align*}
    \MC(\mathcal{I}) \cdot \SD^{\leq 1}(\mathcal{I}) \subset \mathcal{I}.
\end{align*}
The reason why we care about such a property is reflected in the theorem below:

\begin{theorem}[Invariance of maximal contact for MC-invariant ideals]\label{5101}
Assume $\mathcal{I}$ is MC-invariant. For every $y \in Y$ such that $\logord_y(\mathcal{I}) = a \geq 1$, and every pair of maximal contact elements $x$ and $x'$ of $\mathcal{I}$, there exist strict and \'etale morphisms \begin{align*}
    \widetilde{U} \xrightrightarrows[\phi_{x'}]{\phi_x} Y
\end{align*}
from a strict toroidal $\kk$-scheme $\widetilde{U}$ into $Y$, and a point $\widetilde{y}$ of $\widetilde{U}$ such that $\phi_x(\widetilde{y}) = y = \phi_{x'}(\widetilde{y})$, satisfying the following properties: \begin{enumerate}
    \item $\phi_x^\ast(\mathcal{I}) = \phi_{x'}^\ast(\mathcal{I})$;
    \item $\phi_x^\ast(x) = \phi_{x'}^\ast(x')$ in $\MC(\widetilde{\mathcal{I}})$, where $\widetilde{\mathcal{I}}$ denotes the ideal in \textup{(i)}.
\end{enumerate}
\end{theorem}

The statement (and proof) of (\ref{5101}) follows \cite[Lemma 3.5.5]{Wlo05}, \cite[Lemma 5.3.3]{ATW20a} and \cite[Theorem 3.92]{Kol07} closely. See Appendix \ref{appC} for a proof. 

\subsection{Coefficient ideals}\label{5.2}

In this section, we recall the method of taking coefficient ideals. This originates from Hironaka \cite{Hir64}, and has been studied extensively in the papers of Bierstone-Milman (\cite{BM08}, etc), Encinas-Villamayor (\cite{EV00}, etc), W{\l}odarczyk \cite{Wlo05}, and many others. Our treatment closely follows \cite{ATW19}, which studies coefficient ideals from the Rees algebra approach of \cite{EV07}.

For an integer $a \geq 1$, consider the graded $\SO_Y$-subalgebra $\mathscr{G}_\bullet(\mathcal{I},a) \subset \SO_Y[T]$ generated by $\SO_Y$ and $\SD^{\leq j}(\mathcal{I}) \cdot T^{a-j}$ for every $0 \leq j < a$. Its graded pieces are \begin{align*}
    \mathscr{G}_s(\mathcal{I},a) := \left(\prod_{j=0}^{a-1}{(\SD^{\leq j}(\mathcal{I}))^{c_j}} \colon c_j \in \NN,\; \sum_{j=0}^{a-1}{(a-j)c_j} \geq s \right) \subset \SO_Y \qquad \textup{for } s \geq 1.
\end{align*}
The main reason for putting $\SD^{\leq j}(\mathcal{I})$ in degree $a-j$ is the following lemma: 

\begin{lemma}\label{5201}
Let $y \in Y$, and let $a \geq 1$ be an integer. If $\logord_y(\mathcal{I}) \geq a$, then we have $\logord_y(\mathscr{G}_s(\mathcal{I},a)) \geq s$ for every $s \geq 1$.
\end{lemma}

\begin{proof}
Each term $\prod_{j=0}^{a-1}{(\SD^{\leq j}(\mathcal{I}))^{c_j}}$ in $\mathscr{G}_s(\mathcal{I},a)$ has logarithmic order at $y$ given by \begin{align*}
    \sum_{j=0}^{a-1}{c_j \cdot \logord_y(\SD^{\leq j}(\mathcal{I}))} = \sum_{j=0}^{a-1}{c_j(\logord_y(\mathcal{I})-j)} \geq \sum_{j=0}^{a-1}{c_j(a-j)} \geq s
\end{align*}
whence $\logord_y(\mathscr{G}_s(\mathcal{I},a)) \geq s$.
\end{proof}

\begin{remark}\label{5202}
Since the formation of $\SD^{\leq 1}$, as well as taking products and sums of ideals, are functorial for logarithmically smooth morphisms, it follows that the formation of $\mathscr{G}_s(-,a)$ is also functorial for logarithmically smooth morphisms, i.e. if $\widetilde{Y} \rightarrow Y$ is a logarithmically smooth morphism of toroidal $\kk$-schemes, then $\mathscr{G}_s(\mathcal{I},a)\SO_{\widetilde{Y}} = \mathscr{G}_s(\mathcal{I}\SO_{\widetilde{Y}},a)$.
\end{remark}

The graded pieces satisfy the following standard properties:

\begin{lemma}[cf. \textrm{\cite[Proposition 3.99]{Kol07}}]\label{5203}
Let $\mathscr{G}_\bullet = \mathscr{G}_\bullet(\mathcal{I},a)$ as above, and assume $1 \leq a = \max\logord(\mathcal{I}) < \infty$. Then: \begin{enumerate}
    \item $\mathscr{G}_{s+1} \subset \mathscr{G}_s$ for every $s$.
    \item $\mathscr{G}_s \cdot \mathscr{G}_t \subset \mathscr{G}_{s+t}$ for every $s,t$.
    \item $\SD^{\leq 1}(\mathscr{G}_{s+1}) = \mathscr{G}_s$ for every $s$.
    \item $\SD^{\leq s-1}(\mathscr{G}_s) = \mathscr{G}_1 = \MC(\mathcal{I})$ for every $s$. In particular, $s = \max\logord(\mathscr{G}_s)$.
    \item For every $s$, $\mathscr{G}_s$ is MC-invariant.
    \item $\mathscr{G}_s \cdot \mathscr{G}_t = \mathscr{G}_{s+t}$ whenever $t \geq (a-1) \cdot \lcm(2,\dotsc,a)$ and $s$ is a multiple of $\lcm(2,\dotsc,m)$. In fact, the same holds if $t \geq a!$.
    \item $(\mathscr{G}_s)^j = \mathscr{G}_{js}$ whenever $s = r \cdot \lcm(2,\dotsc,a)$ for some $r \geq a-1$. In fact, the same holds for $s = a!$.
    \item $(\SD^{\leq i}(\mathscr{G}_s))^s \subset \mathscr{G}_s^{s-i}$ whenever $s = r \cdot \lcm(2,\dotsc,a)$ for some $r \geq a-1$, and $0 \leq i < s$. In fact, the same holds for $s = a!$.
\end{enumerate}
\end{lemma}

\begin{proof}
Even though we are in the logarithmic case, the proof for \cite[Proposition 3.99]{Kol07} works verbatim, but one should be aware of an inconsequential but noteworthy difference: for the inclusion $\mathscr{G}_s \subset \SD^{\leq 1}(\mathscr{G}_{s+1})$ in (iii), the proof utilizes a maximal contact element $x$ of $\mathcal{I}$ at a point, which in the logarithmic case is an ordinary parameter, and hence the corresponding logarithmic derivation is still $\frac{\partial}{\partial x}$.
\end{proof}

\begin{corollary}\label{5204}
Assume $1 \leq a = \max\logord(\mathcal{I}) < \infty$, and let $y \in Y$. If $\logord_y(\mathcal{I}) = a$, then $\logord_y(\mathscr{G}_s(\mathcal{I},a)) = s$ for every $s \geq 1$. Moreover, if $x$ is a maximal contact element for $\mathcal{I}$ at $y$, then $x$ is also a maximal contact element for $\mathscr{G}_s(\mathcal{I},a)$ at $y$.
\end{corollary}

\begin{proof}
This is a consequence of (\ref{5201}) and (\ref{5203}(iv)).
\end{proof}

With the exception of (viii), all the properties in (\ref{5203}) are self-explanatory. For example, (\ref{5203}(vii)) says that the $(a!)$-Veronese subalgebra $\mathscr{G}_{a!\bullet}(\mathcal{I},a)$ of $\mathscr{G}_\bullet(\mathcal{I},a)$ is generated in degree $1$, i.e. it is the Rees algebra of:

\begin{definition}[Coefficient ideal]\label{5205}
Let $\mathcal{I}$ be an ideal on $Y$, and assume that $1 \leq a = \max\logord(\mathcal{I}) < \infty$. The \emph{coefficient ideal} of the marked ideal $(\mathcal{I},a)$ is\begin{align*}
    \mathscr{C}(\mathcal{I},a) := \mathscr{G}_{a!}(\mathcal{I},a) \subset \SO_Y.
\end{align*}
\end{definition}

Historically, the coefficient ideal provides a method to enrich an ideal with its higher derivatives, which retains information that would otherwise be lost when one restricts the original ideal (as opposed to the coefficient ideal) to a hypersurface of maximal contact.

Finally, let us explicate the property in (\ref{5203}(viii)). Following \cite[Definition 3.83]{Kol07}, we say an ideal $\mathcal{I}$ on $Y$, with $1 \leq a = \max\logord(\mathcal{I}) < \infty$, is \emph{$\SD$-balanced} (in the logarithmic sense) if \begin{align*}
    \SD^{\leq i}(\mathcal{I})^a \subset \mathcal{I}^{a-i} \qquad \textrm{for} \quad 0 \leq i < a.
\end{align*}
In particular, (\ref{5203}(viii)) says that if $a = \max\logord(\mathcal{I}) < \infty$, the coefficient ideal $\mathscr{C}(\mathcal{I},a)$ is $\SD$-balanced.
    
The ``$\SD$-balanced'' property plays a subtle role in our paper. Namely, let $x$ be a maximal contact element of $\mathcal{I}$ at some point $y \in Y$, and denote the corresponding hypersurface of maximal contact by $H$. If one extends $x$ to a system of ordinary parameters at $y$, one easily sees that $\SD^{\leq 1}(\mathcal{I}|_H) \subset \SD^{\leq 1}(\mathcal{I})|_H$. Note, however, that the reverse inclusion does not hold in general. As noted in the paragraph before \cite[Definition 3.83]{Kol07}, the ``$\SD$-balanced'' property provides a partial remedy to this issue. 

Let us be more precise about this by stating the issue in terms of admissibility of toroidal centers. Namely, let $\SJ^{(y)}$ be a toroidal center at $y$, and assume the restriction of $\SJ^{(y)}$ to $H$, denoted $\SJ^{(y)}_H$, is $\mathcal{I}|_H$-admissible. Then a repeated application of (\ref{4502}(i)) tells us that after relacing $\SJ_H^{(y)}$ by some power of itself, $\SJ^{(y)}_H$ is $\SD^{\leq i}(\mathcal{I}|_H)$-admissible. Unfortunately, $\SD^{\leq i}(\mathcal{I}|_H)$-admissibility does not imply $\SD^{\leq i}(\mathcal{I})|_H$-admissibility. However, if one assumes that $\mathcal{I}$ is $\SD$-balanced (with $a = \max\logord(\mathcal{I})$), then $(\SD^{\leq i}(\mathcal{I})|_H)^a \subset (\mathcal{I}|_H)^{a-i}$, so that applying (\ref{3210}(iii)) twice gives the following chain of implications: \begin{align*}
    \textup{$\SJ^{(y)}_H$ is $\mathcal{I}|_H$-admissible} \; &\Rightarrow \; \textup{$(\SJ^{(y)}_H)^{a-i}$ is $(\mathcal{I}|_H)^{a-i}$-admissible} \\
    &\Rightarrow \; \textup{$(\SJ^{(y)}_H)^{a-i}$ is $(\SD^{\leq i}(\mathcal{I})|_H)^a$-admissible} \\
    &\Rightarrow \; \textup{$(\SJ^{(y)}_H)^{\frac{a-i}{a}}$ is $\SD^{\leq i}(\mathcal{I})|_H$-admissible}.
\end{align*}
It turns out that this strategy works out very well in \S\ref{6.2} below (see the proof of (\ref{6201}(i))).

\subsection{Formal decomposition}\label{5.3}

Let $y \in Y$, and assume $\logord_y(\mathcal{I}) = a$ (where $a \geq 1$ is an integer). Let $x_1$ be a maximal contact element of $\mathcal{I}$ at a point $y \in Y$. Extending it to a system of ordinary parameters $x_1,\dotsc,x_n$ at $y$, we have \begin{align*}
    \widehat{\SO}_{Y,y} \simeq \kappa\llbracket x_1,x_2\dotsc,x_n,M \rrbracket, \qquad \textup{where } \kappa = \kappa(y) \; \textup{ and } \; M = \overline{\SM}_{Y,y}.
\end{align*}
For integers $s \geq 1$, \begin{enumerate}
    \item let $\widehat{\mathscr{G}}_s(\mathcal{I},a) = \mathscr{G}_s(\mathcal{I},a)\widehat{\SO}_{Y,y}$,
    \item let $\overline{\mathscr{C}}_s(\mathcal{I},a)$ denote the ideal generated by the image of $\widehat{\mathscr{G}}_s(\mathcal{I},a)$ under the reduction homomorphism $\widehat{\SO}_{Y,y} = \kappa\llbracket x_1,x_2,\dotsc,x_n,M \rrbracket \twoheadrightarrow \kappa\llbracket x_2,\dotsc,x_n,M \rrbracket$,
    \item and let $\widetilde{\mathscr{C}}_s(\mathcal{I},a) = \overline{\mathscr{C}}_s(\mathcal{I},a)\kappa\llbracket x_1,x_2,\dotsc,x_n,M \rrbracket = \overline{\mathscr{C}}_s(\mathcal{I},a)\widehat{\SO}_{Y,y}$.
\end{enumerate}

\begin{proposition}[Formal decomposition, cf. \textup{\cite[Proposition 4.4.1]{ATW19}}]\label{5301}
After passing to completion at $y$, we have \begin{align*}
    \widehat{\mathscr{G}}_s(\mathcal{I},a) = (x_1^s) + (x_1^{s-1})\widetilde{\mathscr{C}}_1(\mathcal{I},a) + \dotsb + (x_1)\widetilde{\mathscr{C}}_{s-1}(\mathcal{I},a) + \widetilde{\mathscr{C}}_s(\mathcal{I},a), \qquad \textup{where } s \geq 1.
\end{align*}
In particular, \begin{align*}
    \widehat{\mathscr{C}}(\mathcal{I},a) = (x_1^{a!}) + (x_1^{a!-1})\widetilde{\mathscr{C}}_1(\mathcal{I},a)+\dotsb+(x_1)\widetilde{\mathscr{C}}_{a!-1}(\mathcal{I},a)+\widetilde{\mathscr{C}}_{a!}(\mathcal{I},a).
\end{align*}
\end{proposition}

\begin{proof}
We shall prove by induction on $s$, with the case $s=1$ being clear. For integers $N \geq s$, we have the ideals $(x_1^{N+1}) \subset \widehat{\mathscr{G}}_s(\mathcal{I},a)$, which are stable under the linear operator $x_1\frac{\partial}{\partial x_1}$. Thus, $x_1\frac{\partial}{\partial x_1}$ descends to a linear operator on $\widehat{\mathscr{G}}_s(\mathcal{I},a)/(x_1^{N+1})$, and decomposes it into a direct sum of $m$-eigenspaces for integers $0 \leq m \leq N$. These $m$-eigenspaces are independent of choice of $N \geq m$. Therefore, we can write the $m$-eigenspace as $x_1^m \cdot \widehat{\mathscr{G}}_s^{(m)}(\mathcal{I},a)$ for subspaces $\widehat{\mathscr{G}}_s^{(m)}(\mathcal{I},a) \subset \kappa\llbracket x_2,\dotsc,x_n,M \rrbracket$, so that \begin{align*}
   \widehat{\mathscr{G}}_s(\mathcal{I},a)/(x_1^{N+1}) = \bigoplus_{m=0}^N{x_1^m \cdot \widehat{\mathscr{G}}_s^{(m)}(\mathcal{I},a)}.
\end{align*}
This implies that \begin{align}\label{eq5.1}
    \widehat{\mathscr{G}}_s(\mathcal{I},a) = \left(x_1^m \cdot \widehat{\mathscr{G}}_s^{(m)}(\mathcal{I},a) \colon 0 \leq m \leq N \right) + \left(x_1^{N+1}\right).
\end{align}
Next, we explicate the terms in the above equation. The simple terms are $\widehat{\mathscr{G}}_s^{(0)}(\mathcal{I},a) = \overline{\mathscr{C}}_s(\mathcal{I},a)$, and $\widehat{\mathscr{G}}_s^{(m)}(\mathcal{I},a) = \kappa\llbracket x_2,\dotsc,x_n,M \rrbracket$ for $m \geq s$. For integers $0 < m < s$, we have \begin{align*}
    \widehat{\mathscr{G}}_s^{(m)}(\mathcal{I},a) = \frac{\partial^m}{\partial x_1^m}\big(x_1^m \cdot \widehat{\mathscr{G}}_s^{(m)}(\mathcal{I},a)\big) &\subset \SD^{\leq m}(\widehat{\mathscr{G}}_s(\mathcal{I},a)) \cap \kappa\llbracket x_2,\dotsc,x_n,M \rrbracket \\
    &= \widehat{\mathscr{G}}_{s-m}(\mathcal{I},a) \cap \kappa\llbracket x_2,\dotsc,x_n,M \rrbracket \subset \overline{\mathscr{C}}_{s-m}(\mathcal{I},a)
\end{align*}
where the equality in the second line follows from (\ref{5203}(iii)). Substituting these into equation \eqref{eq5.1} with $N=s$, we get: \begin{align*}
    \widehat{\mathscr{G}}_s(\mathcal{I},a) \subset \widetilde{\mathscr{C}}_s(\mathcal{I},a) + (x_1)\widetilde{\mathscr{C}}_{s-1}(\mathcal{I},a) + \dotsb + (x_1^{s-1})\widetilde{\mathscr{C}}_1(\mathcal{I},a) + (x_1^s).
\end{align*}
The induction hypothesis gives: \begin{align*}
    (x_1)\widetilde{\mathscr{C}}_{s-1}(\mathcal{I},a) + \dotsb + (x_1^{s-1})\widetilde{\mathscr{C}}_1(\mathcal{I},a) + (x_1^s) = (x_1)\widehat{\mathscr{G}}_{s-1}(\mathcal{I},a) \subset \widehat{\mathscr{G}}_s(\mathcal{I},a).
\end{align*}
Since $\widetilde{\mathscr{C}}_s(\mathcal{I},a) \subset \widehat{\mathscr{G}}_s(\mathcal{I},a)$ as well, the proposition follows.
\end{proof}

\section{Invariants and toroidal centers associated to ideals}\label{chpt6}

\subsection{Defining invariants and toroidal centers at points}\label{6.1}

For an ideal $\mathcal{I}$ on a strict toroidal $\kk$-scheme $Y$ and $y \in Y$, we shall first associate some preliminary data, namely: \begin{enumerate}
    \item a finite sequence of natural numbers $(b_1,\dotsc,b_k) \in \NN^k$,
    \item a finite sequence of ordinary parameters $x_1,\dotsc,x_k$ at $y$,
    \item and an ideal $Q$ of $M = \overline{\SM}_{Y,y}$.
\end{enumerate}
We do this by induction, which terminates only once $Q$ is defined. For the base case, we consider: \begin{enumerate}
    \item[\framebox{\textbf{Case 1a}}] If $\logord_y(\mathcal{I}) = 0$ (i.e. $\mathcal{I}_y = (1)$), then set $k := 1$, $b_1 := 0$, and $Q := \emptyset$. Let $x_1$ be any ordinary parameter at $p$.
    \item[\framebox{\textbf{Case 1b}}] If $\logord_y(\mathcal{I}) = \infty$ (i.e. $\SM(\mathcal{I})_y \neq (1)$), we do not define any $b_i$ or $x_i$ (i.e. set $k := 0$), and define $Q$ by passing the stalk of $\alpha_Y^{-1}(\SM(\mathcal{I}))$ at $y$ to $\overline{\SM}_{Y,y}$, which we denote by $\overline{\alpha_Y^{-1}(\SM(\mathcal{I}))_y}$. Note that it may happen that $\mathcal{I}_y = \SM(\mathcal{I})_y = 0$, in which case $Q := \emptyset$.
    \item[\framebox{\textbf{Case 2}}] If not, set $b_1:=\logord_y(\mathcal{I}) \in \NN_{\geq 1}$ and let $x_1$ be a maximal contact element of $\mathcal{I}$ at $y$. 
\end{enumerate}
In \framebox{\textbf{Case 2}}, set $\mathcal{I}[1] = \mathcal{I}$, and we shall define the remaining $b_i$, $x_i$ and $Q$ by means of induction. Assuming that $\mathcal{I}[i],b_i,x_i$ are defined for $i \leq \ell$, we set \begin{align*}
    \mathcal{I}[\ell+1] := \mathscr{C}(\mathcal{I}[\ell],b_\ell)|_{V(x_1,\dotsc,x_\ell)}.
\end{align*}
In what follows, we pull back the logarithmic structure $\SM_Y$ on $Y$ back to define a logarithmic structure $\alpha_{V(x_1,\dotsc,x_\ell)} \colon \SM_{V(x_1,\dotsc,x_\ell)} \rightarrow \SO_{V(x_1,\dotsc,x_k)}$ on $V(x_1,\dotsc,x_\ell)$. Note that since $x_1,\dotsc,x_\ell$ are ordinary parameters at $y$, $V(x_1,\dotsc,x_\ell)$ is a strict toroidal $\mathds{k}$-scheme under this logarithmic structure.
\begin{enumerate}
    \item[\framebox{\textbf{Case A}}]  If $\SM(\mathcal{I}[\ell+1])_y \neq (1)$ (i.e. $\logord_y(\mathcal{I}[\ell+1]) = \infty$), no further $b_i$ or $x_i$ are defined. Define $Q$ to be the preimage of $\overline{\alpha_{V(x_1,\dotsc,x_\ell)}^{-1}(\SM(\mathcal{I}[\ell+1]))_y}$ under the canonical isomorphism $\overline{\SM}_{Y,y} \xrightarrow{\simeq} \overline{\SM}_{V(x_1,\dotsc,x_k),y}$.
    \item[\framebox{\textbf{Case B}}] If not, set $b_{\ell+1} := \logord_y(\mathcal{I}[\ell+1]) \in \NN_{\geq 1}$, and define $x_{\ell+1}$ to be a lifting to $\SO_Y$ of the maximal contact element of $\mathcal{I}[\ell+1]$ at $y$.
\end{enumerate}
This concludes the induction. Although different choices of ordinary parameters $x_i$ can be made above, the next lemma shows that the $b_i$ and $Q$ are well-defined:

\begin{lemma}\label{6101}The $b_i$ and $Q$ are independent of the choices of ordinary parameters $x_i$ above.
\end{lemma}

\begin{proof}
We proceed by induction on $k=$ the number of $b_i$. The case $k = 0$ occurs if and only if $\SM(\mathcal{I})_y \neq (1)$, in which case there are no $b_i$ and the definition of $Q$ does not require choices. Henceforth, consider $k \geq 1$ (i.e. $\logord_y(\mathcal{I}) < \infty$). Evidently the integer $b_1 = \logord_y(\mathcal{I})$ requires no choices. Next, suppose that we are presented with two choices of maximal contact elements $x,x'$ of $\mathcal{I}$ at $y$. We can replace $Y$ by a neighbourhood of $y$ so that $\max\logord(\mathcal{I}) = b_1$: then $\mathscr{C}(\mathcal{I},b_1)$ is MC-invariant (\ref{5203}(v)), and $x,x'$ are still maximal contact elements of $\mathscr{C}(\mathcal{I},b_1)$ at $y$ (\ref{5204}). Therefore, we can apply (\ref{5101}) to $\mathscr{C}(\mathcal{I},b_1)$: we get strict and \'etale morphisms $\phi_{x,x'} \colon \widetilde{U} \rightrightarrows Y$, and a point $\widetilde{y} \in \widetilde{U}$ such that $\phi_x(\widetilde{y}) = y = \phi_{x'}(\widetilde{y})$. Moreover, $\phi_x^\ast(\mathscr{C}(\mathcal{I},b_1)) = \phi_{x'}^\ast(\mathscr{C}(\mathcal{I},b_1))$ (call this ideal $\widetilde{\mathcal{I}}$) and $z = \phi_x^\ast(x) = \phi_{x'}^\ast(x') \in \widetilde{I}$. Letting $\mathcal{I}[2] = \mathscr{C}(\mathcal{I},b_1)|_{V(x)}$ and $\mathcal{I}[2'] = \mathscr{C}(\mathcal{I},b_1)|_{V(x')}$, we have: \begin{align}\label{eq6.1}
    \phi_x^\ast(\mathcal{I}[2]) = \widetilde{\mathcal{I}}|_{V(z)} = \phi_{x'}^\ast(\mathcal{I}[2']).
\end{align}

If $k = 1$, we are in \framebox{\textbf{Case A}} above. By \cite[Proposition IV.3.1.6]{Ogu18} and (\ref{B113}(iii)), \begin{align}\label{eq6.2}
    \phi_x^\ast\big(\SM(\mathcal{I}[2])\big) = \SM\big(\widetilde{\mathcal{I}}|_{V(z)}\big) = \phi_{x'}^\ast\big(\SM(\mathcal{I}[2'])\big).
\end{align}
Since $\phi_x$ is strict, $\phi_x^\flat \colon \phi_x^\ast(\SM_Y) \rightarrow \SM_{\widetilde{U}}$ is an isomorphism. We therefore get isomorphisms \[
    \overline{\SM}_{V(x),y} \xleftarrow{\simeq} \overline{\SM}_{Y,y} \xrightarrow{\simeq} \overline{\phi_x^\ast(\SM_Y)}_{\widetilde{y}} \xrightarrow{\simeq} \overline{\SM}_{\widetilde{U},\widetilde{y}}
\]
which maps $\overline{\alpha_{V(x)}^{-1}\big(\SM(\mathcal{I}[2])\big)_y}$ on the left, isomorphically, onto $\overline{\alpha_{\widetilde{U}}^{-1}\big( \phi_x^\ast\big(\SM(\mathcal{I}[2])\big)\big)_{\widetilde{y}}}$ on the right. The same statement holds with $V(x)$ replaced by $V(x')$, $\phi_x$ replaced by $\phi_{x'}$, and $\mathcal{I}[2]$ replaced by $\mathcal{I}[2']$. Combining this and (\ref{eq6.2}), one concludes that $Q$ is also independent of choices. 

On the other hand, if $k \geq 2$, we are in \framebox{\textbf{Case B}} above. Then (\ref{eq6.1}) implies \begin{align*}
    \logord_y(\mathcal{I}[2]) = \logord_{\widetilde{y}}(\widetilde{\mathcal{I}}|_{V(z)}) = \logord_y(\mathcal{I}[2']).
\end{align*}
Thus, $b_2$ is independent of choices. By induction hypothesis, the remaining $b_3,b_4,\dotsc$ and $Q$ are independent of choices.
\end{proof}

We are now ready to define the key invariant associated to an ideal at a point:

\begin{definition}[Invariant of an ideal at a point]\label{6102}
Let $\mathcal{I}$ be an ideal on $Y$, and fix $y \in Y$. The (logarithmic) \emph{invariant} of $\mathcal{I}$ at $y$ is defined as: \begin{align*}
    \inv_y(\mathcal{I}) := \begin{cases}
        \Big(b_1,\frac{b_2}{(b_1-1)!},\frac{b_3}{(b_1-1)! \cdot (b_2-1)!},\dotsb,\frac{b_k}{\prod_{i=1}^{k-1}{(b_i-1)!}}\Big) & \textup{if } Q = \emptyset \\
       \Big(b_1,\frac{b_2}{(b_1-1)!},\frac{b_3}{(b_1-1)! \cdot (b_2-1)!},\dotsb,\frac{b_k}{\prod_{i=1}^{k-1}{(b_i-1)!}},\infty\Big) \qquad & \textup{if } Q \neq \emptyset
    \end{cases}
\end{align*}
where $(b_1,\dotsc,b_k)$ and $Q$ are defined for $\mathcal{I}$ at $y$ as before. We will denote the finite entries of $\inv_y(\mathcal{I})$ by $a_i$, so in particular, $a_1=b_1$. We also set $\max\inv(\mathcal{I}) := \max_{y \in Y}{\inv_y(\mathcal{I})}$.
\end{definition}

Observe that $\inv_y(\mathcal{I})$ is the empty sequence $()$ if and only if $\mathcal{I}_y = 0$ (i.e. $y \notin \Supp(\mathcal{I})$). Moreover, $\inv_y(\mathcal{I}) = (0)$ if and only if $\mathcal{I}_y = (1)$ (i.e. $y \notin V(\mathcal{I})$), while $\inv_y(\mathcal{I}) = (a_1)$ for an integer $a_1 \geq 1$ if and only if $\mathcal{I}_y = (x_1^{a_1})$. Finally, $\inv_y(\mathcal{I}) = (\infty)$ if and only if $\SM(\mathcal{I}_y) \neq (1)$ (i.e. $y \in V(\SM(\mathcal{I}))$).

\begin{lemma}\label{6103}
$\inv_y$ satisfies the following properties: \begin{enumerate}
    \item If $\logord_y(\mathcal{I}) = a_1 < \infty$, and $x_1$ is a maximal contact element of $\mathcal{I}$ at $y$, then $\inv_y(\mathcal{I})$ is the concatenation $\Big(a_1,\frac{\inv_y(\mathscr{C}(\mathcal{I},a_1)|_{x_1 = 0})}{(a_1-1)!}\Big)$.
    \item $\inv_y(\mathcal{I})$ is upper semi-continuous on $Y$ (with respect to the lexicographic order which was described in \S\ref{1.1}).
    \item If $\widetilde{Y} \rightarrow Y$ is a logarithmically smooth morphism of strict toroidal $\mathds{k}$-schemes which maps $\widetilde{y} \in \widetilde{Y}$ to $y \in Y$, then $\inv_{\widetilde{y}}(\mathcal{I}\SO_{\widetilde{Y}}) = \inv_y(\mathcal{I})$. If $\widetilde{Y} \to Y$ is moreover surjective, then $\max\inv(\mathcal{I}\mathscr{O}_{\widetilde{Y}}) = \max\inv(\mathcal{I})$.
\end{enumerate}
\end{lemma}

\begin{proof}
Part (i) is evident from (\ref{6102}), while part (iii) follows from (\ref{B115}(iv)) and (\ref{5202}). For part (ii), fix some nondecreasing truncated sequence of nonnegative rational numbers $(a_1,\dotsc,a_k)$ whose last entry could possibly be $\infty$. We need to show the locus $Z$ of points $y \in Y$ such that $\inv_y(\mathcal{I}) \geq (a_1,\dotsc,a_k)$ is closed in $Y$. We do so by induction on $k$. If $k=0$, $Z = Y \setminus \Supp(\mathcal{I})$. Since $Y$ is a disjoint union of its irreducible components (\ref{B105}(iii)), $\Supp(\mathcal{I})$ is a union of some of the irreducible components of $Y$, whence it is open (and closed) in $Y$, so $Z$ is closed in $Y$. Now assume $k \geq 1$. If $a_1 = 0$, $Z = Y$. If $a_1 \in \QQ_{\geq 0} \setminus \ZZ_{\geq 0}$, $Z = V(\SD^{\leq \lceil a_1 \rceil -1}(\mathcal{I}))$ by (\ref{B115}(i)). If $a_1 = \infty$, then $k=1$ and $Z = V(\SM(\mathcal{I}))$ by (\ref{B115}(ii)). Finally, consider $a_1 \in \ZZ_{>0}$. By (\ref{B115}(i)), the locus $W$ of points $y \in Y$ with $\logord_y(\mathcal{I}) > a_1$ is $V(\SD^{\leq a_1}(\mathcal{I}))$. Using part (i) of this lemma and induction hypothesis, the locus $W'$ of points $y \in V(x_1)$ such that $\inv_y(\mathscr{C}(\mathcal{I},a_1)|_{x_1=0}) \geq (a_1-1)! \cdot (a_2,\dotsc,a_k)$ is closed in $V(x_1)$ (and hence, in $Y$). Note that if $y \in W'$, then $\logord_y(\mathcal{I}) \geq a_1$ (if not, the stalk of $\mathscr{C}(\mathcal{I},a_1)$ at $y$ is $(1)$, whence $\inv_y(\mathscr{C}(\mathcal{I},a_1)|_{x_1=0}) = (0) < (a_1-1)! \cdot (a_2,\dotsc,a_k)$). Using part (i) of this lemma again, $Z = W \cup W'$, so $Z$ is closed in $Y$, as desired.
\end{proof}

\begin{definition}[Toroidal center associated to an ideal at a point]\label{6104}
Let $\mathcal{I}$ be a ideal on $Y$, and fix $y \in Y$ such that $\mathcal{I}_y \neq 0$. For a choice of ordinary parameters $x_1,\dotsc,x_k$ associated to $\mathcal{I}$ at $y$ as above, the corresponding toroidal center $\SJ^{(y)}(\mathcal{I})$ at $y$ associated to $\mathcal{I}$ is defined as: \begin{align*}
    \SJ^{(y)}(\mathcal{I}) := \begin{cases}
        \left(x_1^{b_1},x_2^{\frac{b_2}{(b_1-1)!}},x_3^{\frac{b_3}{(b_1-1)! \cdot (b_2-1)!}},\dotsc,x_k^{\frac{b_k}{\prod_{i=1}^{k-1}{(b_i-1)!}}}\right) & \textup{if } Q = \emptyset \\
       \left(x_1^{b_1},x_2^{\frac{b_2}{(b_1-1)!}},x_3^{\frac{b_3}{(b_1-1)! \cdot (b_2-1)!}},\dotsc,x_k^{\frac{b_k}{\prod_{i=1}^{k-1}{(b_i-1)!}}},(Q \subset M)^{\frac{1}{\prod_{i=1}^k{(b_i-1)!}}}\right) \; & \textup{if } Q \neq \emptyset
    \end{cases}
\end{align*}
where $(b_1,\dotsc,b_k)$ and $Q$ are defined for $\mathcal{I}$ at $y$ as before. (We use the convention that $x_1^0 := 1$.) Observe it has invariant equal to $\inv_y(\mathcal{I})$. For the remainder of this paper, we denote $\SJ^{(y)}(\mathcal{I})$ by $(x_1^{a_1},\dotsc,x_k^{a_k},(Q \subset M)^{1/d})$, where $Q$ could be $\emptyset$, and $d$ is always the positive integer $\prod_{i=1}^k{(b_i-1)!}$. 
\end{definition}

We will show later in (\ref{6202}) that $\SJ^{(y)}(\mathcal{I})$ does not actually depend on the choice of ordinary parameters $x_1,\dotsc,x_k$ associated to $\mathcal{I}$ at $y$, which justifies the notation.

\subsection{The associated toroidal center is uniquely admissible}\label{6.2}

The goal of this subsection is to show:

\begin{theorem}[Unique Admissibility]\label{6201}
Let $\mathcal{I}$ be an ideal on $Y$, and fix $y \in Y$ such that $\mathcal{I}_y \neq 0$. \begin{enumerate}
    \item For any choice of ordinary parameters $x_i$ as in \emph{\S\ref{6.1}}, the toroidal center $\SJ^{(y)}(\mathcal{I})$ at $y$ in \emph{(\ref{6104})} is $\mathcal{I}$-admissible.
    \item Every $\mathcal{I}$-admissible toroidal center $\SJ^{(y)}$ at $y$ has invariant $\inv(\SJ^{(y)}) \leq \inv_y(\mathcal{I})$ (where $<$ refers to the lexicographic order which was described in \emph{\S\ref{1.1})}. Consequently, we have the characterization: \begin{align*}
        \inv_y(\mathcal{I}) = \max_{\SJ^{(y)} \textup{ $\mathcal{I}$-admissible}}{\inv(\SJ^{(y)})}.
    \end{align*}
    \item Let $\SJ^{(y)}  = ((x_1')^{a_1},\dotsc,(x_k')^{a_k},(Q' \subset M)^r)$ be a $\mathcal{I}$-admissible toroidal center at $y$, with invariant $\inv(\SJ^{(y)}) = \inv_y(\mathcal{I})$. For any choice of ordinary parameters $x_1,\dotsc,x_k$ associated to $\mathcal{I}$ at $y$ as in \emph{\S\ref{6.1}}, we have $\SJ^{(y)} = (x_1^{a_1},\dotsc,x_k^{a_k},(Q' \subset M)^r)$, after possibly passing to a smaller affine neighbourhood of $y$ on which $\SJ^{(y)}$ is defined.
\end{enumerate}
\end{theorem}

Before proving the theorem, let us note an immediate consequence of (\ref{6201}(iii)):

\begin{corollary}\label{6202}Let $\mathcal{I}$ be an ideal on $Y$, and fix $y \in Y$ such that $\mathcal{I}_y \neq 0$. Then the stalk of $\SJ^{(y)}(\mathcal{I})$ \emph{(\ref{6104})} at $y$ does not depend on the choice of ordinary parameters $x_i$ associated to $\mathcal{I}$ at $y$. \qed
\end{corollary}

We shall divide the proof of (\ref{6201}) into two parts. In the proof of both parts, we will need the following lemma for the induction step:

\begin{lemma}\label{6203}
Let $\mathcal{I}$ be an ideal on $Y$, and let $y \in Y$ be such that $\mathcal{I}_y \neq 0$. Let $\SJ^{(y)} = (x_1^{a_1},\dotsc,x_k^{a_k},(Q \subset M)^r)$ be a toroidal center at $y$, where $k \geq 1$, and $a_1 \geq 1$ is an integer. \begin{enumerate}
    \item Suppose $\SJ^{(y)}$ is $\mathcal{I}$-admissible. Then for any integer $1 \leq m \leq a_1$ and $s \geq 1$, $(\SJ^{(y)})^{s/m}$ is $\mathscr{G}_s(\mathcal{I},m)$-admissible. 
    \item Conversely, if $(\SJ^{(y)})^{(m-1)!}$ is $\mathscr{C}(\mathcal{I},m)$-admissible for some integer $1 \leq m \leq a_1$, then $\SJ^{(y)}$ is $\mathcal{I}$-admissible.
\end{enumerate}
In particular, for any integer $1 \leq m \leq a_1$, $\SJ^{(y)}$ is $\mathcal{I}$-admissible if and only if $(\SJ^{(y)})^{(m-1)!}$ is $\mathscr{C}(\mathcal{I},m)$-admissible.
\end{lemma}

\begin{proof}
If $\SJ^{(y)}$ is $\mathcal{I}$-admissible, iterating (\ref{4502}(i)) tells us that for all $0 \leq j \leq m-1$, $(\SJ^{(y)})^{\frac{a_1-j}{a_1}}$ is $\SD^{\leq j}(\mathcal{I})$-admissible. For natural numbers $c_0,\dotsc,c_{m-1}$, (\ref{3210}(ii)) implies that $(\SJ^{(y)})^{\sum_{j=0}^{m-1}{\frac{a_1-j}{a_1}c_j}}$ is $\big(\prod_{j=0}^{m-1}{(\SD^{\leq j}(\mathcal{I}))^{c_j}}\big)$-admissible. Since $m \leq a_1$, we have $\frac{m-j}{m} \leq \frac{a_1-j}{a_1}$, whence $(\SJ^{(y)})^{\sum_{j=0}^{m-1}{\frac{m-j}{m}c_j}}$ is $\big(\prod_{j=0}^{m-1}{(\SD^{\leq j}(\mathcal{I}))^{c_j}}\big)$-admissible. For $(c_0,\dotsc,c_{m-1}) \in \NN^m$ satisfying $\sum_{j=0}^{m-1}{(m-j)c_j} \geq s$, we have $\sum_{j=0}^{m-1}{\frac{m-j}{m}c_j} \geq \frac{s}{m}$, and hence, $(\SJ^{(y)})^{s/m}$ is $\big(\prod_{j=0}^{m-1}{(\SD^{\leq j}(\mathcal{I}))^{c_j}}\big)$-admissible. By (\ref{3210}(i)), $(\SJ^{(y)})^{s/m}$ is $\mathscr{G}_s(\mathcal{I},m)$-admissible. This proves (i).

Conversely, if $(\SJ^{(y)})^{(m-1)!}$ is $\mathscr{C}(\mathcal{I},m)$-admissible, then $(\SJ^{(y)})^{(m-1)!}$ is $\mathcal{I}^{(m-1)!}$-admissible. By (\ref{3210}(iii)), $\SJ^{(y)}$ is $\mathcal{I}$-admissible. This proves (ii).
\end{proof}

We can now prove (\ref{6201}(i)):

\begin{proof}[Proof of \emph{(\ref{6201}(i))}]
Write $\SJ^{(y)} = \SJ^{(y)}(\mathcal{I})$ in this proof. We proceed by induction on the length $L$ of $\inv_y(\mathcal{I}) = \inv(\SJ^{(y)})$. The base case is $L=1$. The case $\inv(\SJ^{(y)}) = (a_1)$, with $a_1 < \infty$, is evident. If $\inv(\SJ^{(y)}) = (\infty)$, then $\SJ^{(y)}$ is $\mathcal{I}$-admissible because $\SM(\mathcal{I})_y \supset \mathcal{I}_y$. Henceforth, assume $L \geq 2$, so in particular, the first entry in $\inv(\SJ^{(y)})$ is an integer $a_1 \geq 1$. By (\ref{6203}), we may replace $\mathcal{I}$ by $\mathscr{C} = \mathscr{C}(\mathcal{I},a_1)$ and replace $\SJ^{(y)}$ by $(\SJ^{(y)})^{(a_1-1)!}$. By (\ref{3209}), we may pass to completion at $y$, and instead show that $(\widehat{\SJ}^{(y)})^{(a_1-1)!}$ is $\widehat{\mathscr{C}}$-admissible. By (\ref{5301}) we can decompose \begin{align*}
    \widehat{\mathscr{C}} = (x_1^{a_1!}) + (x_1^{a_1!-1})\widetilde{\mathscr{C}}_1 + \dotsb + (x_1)\widetilde{\mathscr{C}}_{a_1!-1} + \widetilde{\mathscr{C}}_{a_1!}, \qquad \textup{where } \widetilde{\mathscr{C}}_{a_1!-i} = \widetilde{\mathscr{C}}_{a_1!-i}(\mathcal{I},a_1),
\end{align*}
and therefore by (\ref{3210}(i)), it remains to show $(\widehat{\SJ}^{(y)})^{(a_1-1)!}$ is $\big((x_1^i)\widetilde{\mathscr{C}}_{a_1!-i}\big)$-admissible for $0 \leq i \leq a_1!$. The case $i = a_1!$ is straightforward.

For the remaining $0 \leq i < a_1!$, let $H$ denote the hypersurface of maximal contact $x_1 = 0$, and let $\SJ_H^{(y)}$ be the restricted toroidal center $\SJ^{(y)}|_{x_1=0}$. By (\ref{6103}(i)), as well as the induction hypothesis (applied to $\mathscr{C}(\mathcal{I},a_1)|_{x_1=0}$), $(\SJ_H^{(y)})^{(a_1-1)!}$ is $\mathscr{C}|_H$-admissible. Since $\mathscr{C}$ is $\SD$-balanced\footnote{(more precisely, after restricting a neighbourhood $U$ of $y$ on which $\max\logord(\mathcal{I}|_U) = a_1$)} (by (\ref{5203}(viii))), we have $(\SD^{\leq i}(\mathscr{C})|_H)^{a_1!} \subset (\mathscr{C}|_H)^{a_1!-i}$. By (\ref{3210}(iii)), we see that $(\SJ_H^{(y)})^{(a_1-1)! \cdot (a_1!-i)}$ is $(\mathscr{C}|_H)^{(a_1!-i)}$-admissible, and hence, $(\SD^{\leq i}(\mathscr{C})|_H)^{a_1!}$-admissible. Consequently, (\ref{3211}) implies that $(\SJ^{(y)})^{(a_1-1)! \cdot (a_1!-i)}$ is $(\SD^{\leq i}(\mathscr{C})|_H\SO_Y)^{a_1!}$-admissible. By a repeated application of (\ref{4502}(ii)), we see that $(\SJ^{(y)})^{(a_1-1)! \cdot a_1!}$ is $\big((x_1^{ia_1!})(\SD^{\leq i}(\mathscr{C})|_H\SO_Y)^{a_1!}\big)$-admissible. By another application of (\ref{3210}(iii)), $(\SJ^{(y)})^{(a_1-1)!}$ is $\big((x_1^i)(\SD^{\leq i}(\mathscr{C})|_H\SO_Y)\big)$-admissible. Recall that $\SD^{\leq i}(\mathscr{C}) = \mathscr{G}_{a_1!-i}$ (\ref{5203}(iii)). Therefore, passing to completion at $y$, $(\widehat{\SJ}^{(y)})^{(a_1-1)!}$ is $\big((x_1^i)\widetilde{\mathscr{C}}_{a_1!-i}\big)$-admissible. This completes the proof.
\end{proof}

Next, we prove the remaining two parts of (\ref{6201}). The proof of these two parts should be compared to the proof of (\ref{3202}) in \S\ref{4.6}.

\begin{proof}[Proof of \emph{(\ref{6201}(ii)$+$(iii))}]
We prove both parts by induction on the length $L$ of $\inv_y(\mathcal{I})$. Consider the base case $L=1$. If $\inv_y(\mathcal{I}) = (\infty)$, there is nothing to show. On the other hand, if $\inv_y(\mathcal{I}) = (a_1)$ with $a_1 < \infty$, then $\mathcal{I}_y = (x_1^{a_1})$ for some ordinary parameter $x_1$ at $y$, and both parts are immediate. Henceforth, assume $L \geq 2$. Let $\SJ^{(y)} = ((x_1')^{b_1},\dotsc,(x_\ell')^{b_\ell},(Q' \subset M)^r)$ be a $\mathcal{I}$-admissible toroidal center at $y$. Since $L \geq 2$, the first entry in $\inv_y(\mathcal{I})$ is an integer $a_1 \geq 1$, where $a_1 = \logord_y(\mathcal{I}) < \infty$. Consequently, $\ell \geq 1$. Applying (\ref{4504}), $b_1 \leq a_1$. If $b_1 < a_1$, $\inv(\SJ^{(y)}) \leq \inv_y(\mathcal{I})$ follows. Thus, assume $b_1 = a_1 < \infty$ for the remainder of this proof. 

Let $x_1$ be a maximal contact element for $\mathcal{I}$ at $y$. Applying (\ref{4502}(i)) repeatedly, $(\SJ^{(y)})^{1/a_1} = (x_1',(x_2')^{b_2/a_1},\dotsc,(x_\ell')^{b_\ell/a_1},(Q' \subset M)^{r/a_1})$ is $\SD^{\leq a_1-1}(\mathcal{I})$-admissible, and hence, $(x_1)$-admissible. Extending $x_1',\dotsc,x_\ell'$ to a system of ordinary parameters $x_1',\dotsc,x_n'$ at $y$, and passing to completion at $y$, we can write the image of $x_1$ under $\SO_{Y,y} \twoheadrightarrow \SO_{\fs_y,y} \to \widehat{\SO}_{\fs_y,y} \simeq \kappa(y)\llbracket x_1',\dotsc,x_n' \rrbracket$ as $\sum_{\vec{\alpha}}{c_{\vec{\alpha}}(x_1')^{\alpha_1}\dotsm (x_n')^{\alpha_n}}$ for some $c_{\vec{\alpha}} \in \kappa(y)$. By (\ref{4601}), $\sum_{i=1}^k{\frac{\alpha_i}{b_i/a_1}} \geq 1$ whenever $c_{\vec{\alpha}} \neq 0$. Consequently, if we let $\ell_0 = \max\lbrace 1 \leq i \leq \ell \colon b_i = a_1 \rbrace \geq 1$, then the image of $x_1$ in $\SO_{\fs_y,y}$ lies in $(x_1',\dotsc,x_{\ell_0}') + \mathfrak{m}_{\fs_y,y}^2$, where $\mathfrak{m}_{\fs_y,y}^2$ is the maximal ideal of $\SO_{\fs_y,y}$. Therefore, after possibly reordering $x_1',\dotsc,x_{\ell_0}'$, we may replace $x_1'$ by $x_1$ so that $(x_1,x_2',\dotsc,x_n')$ is a system of ordinary parameters at $y$. Note that any such reordering does not mess up the presentation of $\SJ^{(y)} = ((x_1')^{b_1},\dotsc,(x_\ell')^{b_\ell},(Q \subset M)^r)$, since $a_1 = b_1 = \dotsb = b_{\ell_0}$. Applying (\ref{4602}): \begin{align*}
    \SJ^{(y)} = (x_1^{a_1},(x_2')^{b_2},\dotsc,(x_k')^{b_\ell},(Q' \subset M)^r).
\end{align*}
The next natural step is to pass to the induction step.

Let $\mathscr{C} = \mathscr{C}(\mathcal{I},a_1)$. By (\ref{6203}), $(\SJ^{(y)})^{(a_1-1)!}$ is $\mathscr{C}$-admissible. Let $H$ denote the hypersurface of maximal contact given by $x_1=0$, and let $\SJ^{(y)}_H$ denote the restricted toroidal center $\SJ^{(y)}|_{x_1=0}$. Then $(\SJ^{(y)}_H)^{(a_1-1)!}$ is $\mathscr{C}|_H$-admissible. By induction hypothesis (for (\ref{6201}(ii))) applied to $\mathscr{C}|_H$, we see that $\inv(\SJ_H^{(y)})^{(a_1-1)!} \leq \inv_y(\mathscr{C}|_H)$, so $\inv(\SJ_H^{(y)}) \leq \frac{1}{(a_1-1)!} \cdot \inv_y(\mathscr{C}|_H)$. Applying (\ref{6103}(i)), we obtain $\inv(\SJ^{(y)}) = (a_1,\inv(\SJ_H^{(y)}))\leq \big(a_1,\frac{\inv_y(\mathscr{C}|_H)}{(a_1-1)!}\big) = \inv_y(\mathcal{I})$, proving (\ref{6201}(ii)). 

In the event that $\inv(\SJ^{(y)}) = \inv_y(\mathcal{I})$, then $\inv(\SJ^{(y)}_H)^{(a_1-1)!} = \inv_y(\mathscr{C}|_H)$, so that $\ell = k$ and $b_i = a_i$ for $1 \leq i \leq k = \ell$. Let $x_1,x_2,\dotsc,x_k$ be ordinary parameters associated to $\mathcal{I}$ at $y$ (where $x_1$ was arbtirarily chosen earlier), as in \S\ref{6.1}. By induction hypothesis (for (\ref{6201}(iii))) applied to the $\mathscr{C}|_H$-admissible toroidal center $(\SJ_H^{(y)})^{(a_1-1)!}$ at $y$, we have \begin{align*}
    \SJ^{(y)}_H = ((x_2')^{a_2},\dotsc,(x_k')^{a_k},(Q' \subset M)^r) = (x_2^{a_2},\dotsc,x_k^{a_k},(Q' \subset M)^r).
\end{align*}
In the above expression, $x_i'$ is more precisely the reduction of $x_i'$ modulo $x_1=0$, and similarly $x_i$ is the reduction of $x_i$ modulo $x_1=0$. We claim that this implies \begin{align*}
    \SJ^{(y)} = (x_1^{a_1},(x_2')^{a_2},\dotsc,(x_k')^{a_k},(Q' \subset M)^r) = (x_1^{a_1},x_2^{a_2},\dotsc,x_k^{a_k},(Q' \subset M)^r).
\end{align*}
This follows by the same method illustrated in (\ref{3205}), i.e. by checking both sides are equal as idealistic $\QQ$-exponents, and hence, as integrally closed Rees algebras. This proves (\ref{6201}(iii)).
\end{proof}

\begin{corollary}\label{6204}
Let $\mathcal{I}$ be an ideal on $Y$, and fix $y \in Y$ such that $\mathcal{I}_y \neq 0$. Then: \begin{enumerate}
    \item $\inv_y(\mathcal{I}^m) = m \cdot \inv_y(\mathcal{I})$.
    \item Assume $\logord_y(\mathcal{I}) = a_1 < \infty$. For any integer $1 \leq m \leq a_1$, $\inv_y(\mathscr{C}(\mathcal{I},m)) = (m-1)! \cdot \inv_y(\mathcal{I})$.
\end{enumerate}
\end{corollary}

\begin{proof}
Apply (\ref{6201}(ii)), in conjunction with (\ref{3210}(iii)) and (\ref{6203}).
\end{proof}

\subsection{Compatibility of associated toroidal centers}\label{6.3}

In \S\ref{6.1}, we defined the toroidal center associated to an ideal $\mathcal{I} \subset \SO_Y$ at a point $y \in Y$. The next theorem glues toroidal centers at points $y \in Y$ with invariant $\inv_y(\mathcal{I}) = \max\inv(\mathcal{I})$:

\begin{theorem}[Gluing]\label{6301}Let $\mathcal{I}$ be a nowhere zero ideal on $Y$, and define $\max\inv(\mathcal{I}) := \max_{y \in Y}{\inv_y(\mathcal{I})}$. There exists a unique $\mathcal{I}$-admissible toroidal center $\SJ = \SJ(\mathcal{I})$ on $Y$ such that for all $y \in Y$, there exists an open affine neighbourhood $U_y$ of $y$ on which: \begin{enumerate}
    \item If $\inv_y(\mathcal{I}) = \max\inv_y(\mathcal{I})$, then $\SJ|_{U_y}$ is the toroidal center $\SJ^{(y)}(\mathcal{I})$ at $y$ \emph{(\ref{6104})}.
    \item If $\inv_y(\mathcal{I}) < \max\inv(\mathcal{I})$, then $\SJ|_{U_y} = \SO_Y[T]|_{U_y}$.
\end{enumerate}
\end{theorem}

\begin{definition}[Toroidal center associated to an ideal]\label{6302}
Let $\mathcal{I}$ be a nowhere zero ideal on $Y$. The toroidal center associated to $\mathcal{I}$ is $\SJ(\mathcal{I})$ in (\ref{6301}).
\end{definition}

\begin{proof}
Since $\inv_y(\mathcal{I})$ is upper semi-continuous (\ref{6103}(ii)), the locus $V$ of points $y \in Y$ where $\inv_y(\mathcal{I}) < \max\inv(\mathcal{I})$ is open. We claim that we can glue the following: \begin{enumerate}
    \item[(I)] $\SO_Y[T]|_V$;
    \item[(II)] for each $y \in Y$ with $\inv_y(\mathcal{I}) = \max\inv(\mathcal{I})$, the toroidal center $\SJ^{(y)}(\mathcal{I})$ restricted to an appropriately chosen open affine neighbourhood $U_y$ of $y$,
\end{enumerate}
to obtain a toroidal center $\SJ$ on $Y$. This toroidal center would have the desired properties.

First fix $y \in Y$ with $\inv_y(\mathcal{I}) = \max\inv(\mathcal{I})$. Let $\SJ^{(y)}(\mathcal{I}) = (x_1^{a_1},\dotsc,x_k^{a_k},(Q \subset M)^{1/d})$ be defined on an open affine neighbourhood $U_y$ of $y$ in $Y$. Recall that $x_i$ are choices of ordinary parameters associated to $\mathcal{I}$ at $y$, and $Q$ is the preimage of $\overline{\alpha_{V(x_1,\dotsc,x_k)}^{-1}(\SM(\mathcal{I}[k+1]))_y}$ under the canonical isomorphism $M=\overline{\SM}_{Y,y} \xrightarrow{\simeq} \overline{\SM}_{V(x_1,\dotsc,x_k),y}$, where $\mathcal{I}[k+1]$ is the ideal on $V(x_1,\dotsc,x_k)$ which was defined inductively in \S\ref{6.1}. Moreover, we have a chart $\beta \colon M \rightarrow H^0(U_y,\SM_Y|_{U_y})$ which is neat at $y$. Then our claim in the preceding paragraph amounts to showing that after possibly shrinking $U_y$, one has, for each $y' \in U_y$, the following statements: \begin{enumerate}
    \item[(a)] If $\inv_{y'}(\mathcal{I}) = \max\inv(\mathcal{I})$, then the stalks of $\SJ^{(y)}(\mathcal{I})$ and $\SJ^{(y')}(\mathcal{I})$ at $y'$ coincide.
    \item[(b)] If $\inv_{y'}(\mathcal{I}) < \max\inv(\mathcal{I})$, then the stalk of $\SJ^{(y)}(\mathcal{I})$ at $y'$ is $\SO_{Y,y'}[T]$.
\end{enumerate}

For part (a), $x_1,\dotsc,x_k$ are also ordinary parameters associated to $\mathcal{I}$ at $y'$. By unique admissibility (\ref{6201}(iii)), we have $\SJ^{(y')}(\mathcal{I}) = (x_1^{a_1},\dotsc,x_k^{a_k},(Q' \subset M')^{1/d})$, and (\ref{6101}) says that $Q'$ is equal to the preimage of  $\overline{\alpha_{V(x_1,\dotsc,x_k)}^{-1}(\SM(\mathcal{I}[k+1]))_{y'}}$ under the canonical isomorphism $M' = \overline{\SM}_{Y,y'} \xrightarrow{\simeq} \overline{\SM}_{V(x_1,\dotsc,x_k),y'}$. On the other hand, (\ref{6303}(ii)) says the ideal of $\SM_{V(x_1,\dotsc,x_k)}|_{U_y \cap V(x_1,\dotsc,x_k)}$ generated by the image of $\overline{\alpha_{V(x_1,\dotsc,x_k)}^{-1}(\SM(\mathcal{I}[k+1]))_y}$ under the chart $\overline{\beta}$: \begin{equation*}
    \begin{tikzcd}
    M = \overline{\SM}_{Y,y} \arrow[to=1-2, "\beta"] \arrow[to=2-1, "\simeq"] & H^0(U_y,\SM_Y|_{U_y}) \arrow[to=2-2] \\
    \overline{\SM}_{V(x_1,\dotsc,x_k),y} \arrow[to=2-2,"\overline{\beta}"] & H^0(U_y \cap V(x_1,\dotsc,x_k),\SM_{V(x_1,\dotsc,x_k)}|_{U_y \cap V(x_1,\dotsc,x_k)})
    \end{tikzcd}
\end{equation*}
is equal to $\alpha_{V(x_1,\dotsc,x_k)}^{-1}(\SM(\mathcal{I}[k+1]))|_{U_y \cap V(x_1,\dotsc,x_k)}$, from which part (a) follows.

For part (b), let $\inv_{y'}(\mathcal{I}) = (a_1',\dotsc,a_\ell') < \max\inv(\mathcal{I})$. First consider the case when there exists $1 \leq j \leq k$ such that $a_j' < a_j$. Let $j_0 = \min\lbrace 1 \leq j \leq k \colon a_j' < a_j \rbrace$. By (\ref{B115}(i)) and unique admissibility (\ref{6201}(iii)), we may adjust $x_{j_0}$, if needed, so that $x_{j_0}$ is a unit in $\SO_{Y,y'}$, which yields part (b) for this first case. 

The other case occurs when $a_i' = a_i$ for all $1 \leq i \leq k$. In this case, $x_1,\dotsc,x_k$ are also ordinary parameters associated to $\mathcal{I}$ at $y'$ (as in part (a)). Let us always rule out the case $Q = \emptyset$ (which occurs if and only if $\mathcal{I}[k+1]_y = 0$, or equivalently, $\inv_y(\mathcal{I}) = (a_1,\dotsc,a_k)$), by shrinking $U_y$ so that $\mathcal{I}[k+1]_{y'} = 0$ (and hence $\inv_{y'}(\mathcal{I}) = \inv_y(\mathcal{I})$) for every $y' \in U_y \cap V(x_1,\dotsc,x_k)$. On the other hand, if $Q' \neq \emptyset$, then $a_{k+1}' < \infty$, and (\ref{B115}(i)) implies $\SM(\mathcal{I}[k+1])_{y'} = (1)$. Combining that with (\ref{6303}(ii)) as in part (a) completes the proof for part (b) in the second case.
\end{proof}

\begin{lemma}\label{6303}
Let $Y$ be a fs Zariski logarithmic scheme, and let $\beta \colon M \rightarrow H^0(Y,\SM_Y)$ be a chart for $\SM_Y$ which is neat at some $y \in Y$ (so we shall identify $M = \overline{\SM}_{Y,y}$ in the statements below). For an ideal $\mathcal{Q} \subset H^0(Y,\SM_Y)$, we have: \begin{enumerate}
    \item $\beta^{-1}(\mathcal{Q}) = \overline{\mathcal{Q}_y}$ (where the latter denotes the passage of the stalk of $\mathcal{Q}$ at $y$ to $\overline{\SM}_{Y,y}$).
    \item The image of $\overline{\mathcal{Q}_y} \subset \overline{\SM}_{Y,y} = M$ under $\beta$ generates $\mathcal{Q}$.
\end{enumerate}
In particular, $\mathcal{Q} = \emptyset$ if and only if $\overline{\mathcal{Q}}_y = \emptyset$.
\end{lemma}

\begin{proof}
For (i), $\gamma \colon \overline{\SM}_{Y,y} \xrightarrow{\beta} H^0(Y,\SM_Y) \rightarrow \SM_{Y,y}$ defines a splitting $\SM_{Y,y} \simeq \overline{\SM}_{Y,y} \oplus \SO_{Y,y}^\ast$. The splitting allows us to write every $m \in \SM_{Y,y}$ as $(\overline{m},u_m)$ for unique $\overline{m} \in \overline{\SM}_{Y,y}$ and $u_m \in \SO_{Y,y}^\ast$. Under this notation, we have $\beta^{-1}(\mathcal{Q}) = \gamma^{-1}(\mathcal{Q}_y) = \lbrace \overline{m} \in \overline{\SM}_{Y,y} \colon (\overline{m},0) \in \mathcal{Q}_y \rbrace = \lbrace \overline{m} \in \overline{\SM}_{Y,y} \colon m = (\overline{m},u_m) \in \mathcal{Q}_y \rbrace = \overline{\mathcal{Q}_y}$, as desired. 

For (ii), $\beta$ factors as $M \stackrel{\iota}{\hookrightarrow} M \oplus \SO_Y^\ast \stackrel{\pi}{\twoheadrightarrow} \SM_Y$. It then suffices to show that $\iota(\overline{\mathcal{Q}_y})$ generates $\pi^{-1}(\mathcal{Q})$. But (i) implies $\overline{\mathcal{Q}_y} = \iota^{-1}(\pi^{-1}(\mathcal{Q})) = \lbrace m \in M \colon (m,0) = \iota(m) \in  \pi^{-1}(\mathcal{Q}) \rbrace$, whose image under $\iota$ evidently generates $\pi^{-1}(\mathcal{Q})$.
\end{proof}

\subsection{The case of toroidal Deligne-Mumford stacks over \texorpdfstring{$\mathds{k}$}{k}}\label{6.4}

The goal of this section is to extend the definition of associated toroidal centers and associated invariants to toroidal $\mathds{k}$-schemes, or more generally, \emph{toroidal Deligne-Mumford stacks over $\mathds{k}$}.

\begin{lemma}[Functoriality of associated toroidal centers]\label{6401}
Let $f \colon \widetilde{Y} \rightarrow Y$ be a logarithmically smooth morphism of strict toroidal $\mathds{k}$-schemes, which maps $\widetilde{y} \in Y$ to $y \in Y$. For an ideal $\mathcal{I}$ on $Y$ satisfying $\mathcal{I}_y \neq 0$, we have $\SJ^{(y)}(\mathcal{I})\SO_{\widetilde{Y}} = \SJ^{(\widetilde{y})}(\mathcal{I}\SO_{\widetilde{Y}})$ on an affine open neighbourhood of $\widetilde{y}$. If $\mathcal{I}$ is a nowhere zero ideal on $Y$ and $f$ is moreover surjective, then $\SJ(\mathcal{I})\SO_{\widetilde{Y}} = \SJ(\mathcal{I}\SO_{\widetilde{Y}})$.
\end{lemma}

\begin{proof}
For the first assertion, we may replace $Y$ with an affine open neighbourhood of $y$ on which $\SJ^{(y)}(\mathcal{I})$ is defined. Firstly observe that if $\logord_y(\mathcal{I}) = \infty$, then $\SM(\mathcal{I})\SO_{\widetilde{Y}} = \SM(\mathcal{I}\SO_{\widetilde{Y}})$ (\ref{B113}(iii)), and the lemma is immediate. On the other hand, if $\logord_y(\mathcal{I}) = b_1 < \infty$, then any maximal contact element $x_1$ of $\mathcal{I}$ at $y$ is also a maximal contact element of $\mathcal{I}\SO_{\widetilde{Y}}$ at $\widetilde{y}$. Let $V_Y(x_1)$ (resp. $V_{\widetilde{Y}}(x_1)$) be the hypersurface on $Y$ (resp. $\widetilde{Y}$) given by $x_1 = 0$. We restrict to the logarithmically smooth morphism $V_{\widetilde{Y}}(x_1) \rightarrow V_Y(x_1)$. By (\ref{5202}), we have $\mathscr{C}(\mathcal{I}\SO_{\widetilde{Y}},b_1)|_{V_{\widetilde{Y}}(x_1)} = \mathscr{C}(\mathcal{I},b_1)\SO_{\widetilde{Y}}|_{V_{\widetilde{Y}}(x_1)} = \mathscr{C}(\mathcal{I},b_1)|_{V_Y(x_1)}\SO_{V_{\widetilde{Y}}(x_1)}$. The first assertion is then proven by applying the induction hypothesis to the ideal $\mathscr{C}(\mathcal{I},b_1)|_{V_Y(x_1)}$ on $V_Y(x_1)$. The second assertion follows from the first, since surjectivity of $f$ implies that  $\max\inv(\mathcal{I}\mathscr{O}_{\widetilde{Y}}) = \max\inv(\mathcal{I})$ (\ref{6103}(iii)).
\end{proof}

\begin{corollary}\label{6402}
Let $Y$ be a toroidal Deligne-Mumford stack over $\mathds{k}$, and fix an atlas $p_{1,2} \colon Y_1 \rightrightarrows Y_0$ of $Y$ by schemes such that $Y_0$ is a strict toroidal $\mathds{k}$-scheme\footnote{See paragraph after (\ref{B201}).}. Let $y \in \abs{Y}$, and let $\mathcal{I}$ be an ideal on $Y$ such that $\mathcal{I}_y \neq 0$. \begin{enumerate}
    \item If $y_1,y_2 \in Y_0$ are points over $y$, then $\inv_{y_1}(\mathcal{I}\SO_{Y_0}) = \inv_{y_2}(\mathcal{I}\SO_{Y_0})$.
    \item If $y_1$ is a point over $y$, the toroidal center $\SJ^{(y_1)}(\mathcal{I}\SO_{Y_0})$ descends to a toroidal center\footnote{One can extend the definition of toroidal centers to toroidal Deligne-Mumford stacks over $\mathds{k}$, which we have opted not to state explicitly.} $\SJ^{(y)}(\mathcal{I})$ on an open substack of $Y$ containing $y$. 
\end{enumerate}
If $\mathcal{I}$ is a nowhere zero ideal on $Y$, then the toroidal center $\SJ(\mathcal{I}\SO_{Y_0})$ descends to a toroidal center $\SJ(\mathcal{I})$ on $Y$.
\end{corollary}

Because of (\ref{6402}(i)), we can define the invariant $\inv_y(\mathcal{I})$ of $\mathcal{I}$ at $y$ to be $\inv_{y_1}(\mathcal{I}\SO_{Y_0})$ for any point $y_1 \in Y_0$ above $y$.

\begin{proof}
Let $(y_1,y_2) \in Y_1$ denote the point mapping to $y_i$ via $p_i$ for $i=1,2$. Since $p_{1,2}$ are both strict and \'etale, (\ref{6103}(iii)) implies $\inv_{y_1}(\mathcal{I}\SO_{Y_0}) = \inv_{(y_1,y_2)}(\mathcal{I}\SO_{Y_1}) = \inv_{y_2}(\mathcal{I}\SO_{Y_0})$, so part (i) follows. If that invariant is equal to $\max\inv(\mathcal{I}\SO_{Y_0})$, then (\ref{6401}) implies $p_1^\ast\SJ^{(y_1)}(\mathcal{I}\SO_{Y_0}) = \SJ^{(y_1,y_2)}(\mathcal{I}\SO_{Y_1}) = p_2^\ast\SJ^{(y_2)}(\mathcal{I}\SO_{Y_0})$. If not, evidently the same equality holds. Therefore, we obtain the desired descent in the final statement. Part (ii) is a consequence of the final statement by replacing $Y$ with an invariant open affine neighbourhood of $y_1$ on which $\SJ^{(y_1)}(\mathcal{I}\SO_{Y_0})$ is defined.
\end{proof}

\section{Logarithmic principalization}\label{chpt7}

\subsection{Statement of theorem}\label{7.1}

The goal of this section is to prove:

\begin{theorem}[Logarithmic Principalization]\label{7101}There is a functor $F_{\textup{log-pr}}$ associating to \begin{quote}
    a nowhere zero, proper ideal $\mathcal{I}$ on a toroidal Deligne-Mumford stack $Y$ over a field $\mathds{k}$ of characteristic zero
\end{quote}
an $\mathcal{I}$-admissible toroidal center $\SJ = \SJ(\mathcal{I})$ with reduced toroidal center $\overline{\SJ}$, weighted toroidal blow-up $Y' = \Bl_Y(\overline{\SJ}) \rightarrow Y$, and the weak transform $F_{\textup{log-pr}}(\mathcal{I} \subsetneq \SO_Y) = (\mathcal{I}' \subset \SO_{Y'})$ such that $\max\inv(\mathcal{I}') < \max\inv(\mathcal{I})$. Functoriality here is with respect to logarithmically smooth, surjective morphisms. 

In particular, there is an integer $N \geq 1$ so that the iterated application $(\mathcal{I}_N \subset \SO_{Y_N}) = F_{\textup{log-pr}}^{\circ N}(\mathcal{I} \subsetneq \SO_Y)$ of $F_{\textup{log-pr}}$ has $\mathcal{I}_N = (1)$. This stabilized functor $F_{\textup{log-pr}}^{\circ \infty}$ is functorial for all logarithmically smooth morphisms, whether or not surjective.
\end{theorem}

We prove the principalization theorem as a consequence of the results in \S\ref{7.2}. We remind the reader that the notion of weak transform was introduced in (\ref{4405}).

\subsection{The invariant drops}\label{7.2}

Let $\mathcal{I}$ be an ideal on a strict toroidal $\mathds{k}$-scheme $Y$, and let $y \in Y$ such that $\mathcal{I}_y \neq 0$. Let $\SJ = \SJ^{(y)}(\mathcal{I}) = (x_1^{a_1},\dotsc,x_k^{a_k},(Q \subset M)^{1/d})$ be the toroidal center associated to $\mathcal{I}$ at $y$ (with the notation following (\ref{6104})), which has invariant $\inv(\SJ) = \inv_y(\mathcal{I})$. Let $\overline{\SJ}$ be the reduced toroidal center associated to $\SJ$, so that $\overline{\SJ} = \SJ^{1/(a_1n_1)} = (x_1^{1/n_1},\dotsc,x_k^{1/n_k},(Q \subset M)^{1/(a_1n_1d)})$ if $k \geq 1$, and $\overline{\SJ} = \SJ$ if $k=0$.

For this section only, let us work locally at $y$ and replace $Y$ with the open affine neighbourhood of $y$ on which $\SJ$ is defined. For any integer $c \geq 1$, we shall write $Y_c' \rightarrow Y$ for the weighted toroidal blow-up along $\overline{\SJ}^{1/c}$, with exceptional ideal $\SE_c$. By (\ref{4404}), there exists an ideal $\mathcal{I}_c'$ on $Y_c'$ such that $\mathcal{I}\SO_{Y_c'}$ factors as $\SE_c^{a_1n_1c} \cdot \mathcal{I}_c'$ if $k \geq 1$; on the other hand, if $k = 0$, $\mathcal{I}\SO_{Y_c'}$ factors as $\SE_c^{c} \cdot \mathcal{I}_c'$. The goal of this section is to show:

\begin{theorem}[The invariant drops]\label{7201}
Let the notation be as above, and assume $\mathcal{I}_y \neq (1)$. For every integer $c \geq 1$, and every point $y' \in \abs{Y_c'}$ over $y$, we have $\inv_{y'}(\mathcal{I}'_c) < \inv_y(\mathcal{I})$.
\end{theorem}

Of course, we are only interested in the theorem for the case $c=1$. We prove it for all integers $c \geq 1$, so that induction can take place. Let us first deal with the special case $k = 0$:

\begin{lemma}[Cleaning up, cf. \textup{\cite[Proposition 2.2.1]{ATW20a}}]\label{7202}
Let the notation be as above. Assume $k=0$, i.e. $\inv_y(\mathcal{I}) = (\infty)$, and write $\SJ = (Q \subset M)$, where $M = \overline{\SM}_{Y,y}$, and $Q = \overline{\alpha_Y^{-1}(\SM(\mathcal{I}))_y}$. Then \emph{(\ref{7201})} holds.
\end{lemma}

\begin{proof}
We use a $(mT^c)$-chart as in (\ref{4402}), where $m$ belongs to a fixed finite set of generators for $Q$. In this case, $Y_c'$ is the pullback of the logarithmically smooth morphism $\Spec(M_m \rightarrow \mathds{k}[x_1,\dotsc,x_n,M_m]) \rightarrow \Spec(M \rightarrow \mathds{k}[x_1,\dotsc,x_n,M])$ back to $Y$, via the strict and smooth morphism $Y \rightarrow \Spec(M \rightarrow \mathds{k}[x_1,\dotsc,x_n,M])$, where $M_m$ is the saturation of the submonoid of $M[m^{1/c}]^{\gp}$ generated by $M[m^{1/c}]$ and $\lbrace q' = q/m = q/u^c \colon q \in Q \rbrace$. Since $Y_c' \rightarrow Y$ is logarithmically smooth, (\ref{B113}(iii)) implies that $\SM(\mathcal{I}\SO_{Y_c'}) = \SM(\mathcal{I})\SO_{Y_c'}$, so $\SM(\mathcal{I}\SO_{Y_c'})_{y'} = \SM(\mathcal{I})_y\SO_{Y_c',y'}$. Since every $q \in Q$ factors as $q' \cdot u^c$ in $M_m$, we have $\SM(\mathcal{I})_y\SO_{Y_c',y'} = (\SE_c)_{y'}^c$. Therefore, $(\mathcal{I}\SO_{Y_c'})_{y'} = (\SE_c)_{y'}^c \cdot (\mathcal{I}_c')_{y'} = \SM(\mathcal{I}\SO_{Y_c'})_{y'} \cdot (\mathcal{I}_c')_{y'}$. Applying (\ref{B113}(iv)), one sees that $\SM(\mathcal{I}_c')_{y'} = (1)$, whence $\logord_{y'}(\mathcal{I}_c') < \infty$. Thus, $\inv_y(\mathcal{I}_c') < (\infty) = \inv_y(\mathcal{I})$.
\end{proof}

For the case $k \geq 1$, the next lemma (and its corollary) shows we can replace $\mathcal{I}$ by the coefficient ideal $\mathscr{C}(\mathcal{I},a_1)$:

\begin{lemma}[cf. \textup{\cite[Lemma 3.3]{BM08}}]\label{7203}
Let the notation be as above. Assume $k \geq 1$, so that $a_1 = \logord_y(\mathcal{I}) < \infty$, and let $\mathscr{C} = \mathscr{C}(\mathcal{I},a_1)$. For every integer $c \geq 1$, factorize $\mathscr{C}\SO_{Y'_c} = \SE_c^{a_1!n_1c} \cdot \mathscr{C}'_c$ for some ideal $\mathscr{C}_c'$ on $Y_c'$, as in \emph{(\ref{4404})}. Then we have the inclusions: $(\mathcal{I}_c')^{(a_1-1)!} \subset \mathscr{C}_c' \subset \mathscr{C}(\mathcal{I}_c',a_1)$.
\end{lemma}

\begin{proof}
We have: \begin{align*}
    \mathscr{C}\SO_{Y_c'} = \left(\prod_{j=0}^{a_1-1}{(\SD^{\leq j}_Y(\mathcal{I})\SO_{Y_c'})^{c_j}} \colon c_j \in \NN,\; \sum_{j=0}^{a_1-1}{(a_1-j)c_j \geq a_1!}\right).
\end{align*}
Applying (\ref{4503}), we see that for every $1 \leq j < a_1$, \begin{align*}
    \SD^{\leq j}(\mathcal{I})\SO_{Y_c'} \subset \SE_c^{(a_1-j)n_1c} \cdot \SD^{\leq j}(\mathcal{I}_c').
\end{align*}
We also have $\mathcal{I}\SO_{Y_c'} = \SE_c^{a_1n_1c} \cdot \mathcal{I}_c'$. Plugging this into the first equation yields: \begin{align*}
    \mathscr{C}\SO_{Y_c'} \subset \SE_c^{a_1!n_1c} \cdot \left(\prod_{j=0}^{a_1-1}{(\SD^{\leq j}(\mathcal{I}_c'))^{c_j}} \colon c_j \in \NN,\; \sum_{j=0}^{a_1-1}{(a_1-j)c_j} \geq a_1!\right) = \SE_c^{a_1!n_1c} \cdot \mathscr{C}(\mathcal{I}'_c,a_1).
\end{align*}
Thus, we get the second inclusion $\mathscr{C}_c' \subset \mathscr{C}(\mathcal{I}_c',a_1)$. The first inclusion follows from the inclusion $\mathscr{C}\SO_{Y_c'} \supset (\mathcal{I}\SO_{Y_c'})^{(a_1-1)!} = \SE_c^{a_1!n_1c} \cdot (\mathcal{I}_c')^{(a_1-1)!}$.
\end{proof}

\begin{corollary}\label{7204}
Assume the hypotheses of \emph{(\ref{6203})}. For every point $y' \in \abs{Y_c'}$ over $y$, we have: \begin{enumerate}
    \item $\inv_{y'}(\mathscr{C}_c') = (a_1-1)! \cdot \inv_{y'}(\mathcal{I}_c')$.
    \item $\inv_{y'}(\mathcal{I}_c') < \inv_y(\mathcal{I})$ if and only if $\inv_{y'}(\mathscr{C}_c') < \inv_y(\mathscr{C})$.
\end{enumerate}
\end{corollary}

\begin{proof}
By (\ref{7203}), we have: \begin{align*}
    \inv_{y'}\big((\mathcal{I}_c')^{(a_1-1)!}\big) \geq \inv_{y'}(\mathscr{C}_c') \geq \inv_{y'}(\mathscr{C}(\mathcal{I}_c',a_1-1)),
\end{align*}
but (\ref{6204}(ii)) implies $\inv_{y'}(\mathscr{C}(\mathcal{I}_c',a_1)) = (a_1-1)! \cdot \inv_{y'}(\mathcal{I}_c') = \inv_{y'}\big((\mathcal{I}_c')^{(a_1-1)!}\big)$. This forces equality throughout, yielding (i). Part (ii) follows from (i) and (\ref{6204}(ii)).
\end{proof}

\begin{proof}[Proof of \emph{(\ref{7201})}]
We induct on the length $L$ of $\inv_y(\mathcal{I})$. First consider the base case $L=1$. The sub-case $\inv_y(\mathcal{I}) = (\infty)$ is settled in (\ref{7202}). On the other hand, if $\inv_y(\mathcal{I}) = (a_1)$ with $a_1 < \infty$, then $\mathcal{I}_y = (x_1^{a_1})$, with weak transform $(\mathcal{I}_c')_{y'} = (1)$. Henceforth, assume $L \geq 2$. In particular, $k \geq 1$, so (\ref{7204}) says that we can replace\footnote{Note that $ \SJ^{(y)}(\mathscr{C}) = \SJ^{(y)}(\mathcal{I})^{(a_1-1)!} = \overline{\SJ}^{a_1!n_1}$, so $\overline{\SJ}$ is also the reduced toroidal center associated to $\SJ^{(y)}(\mathscr{C})$.} $\mathcal{I}$ with $\mathscr{C} = \mathscr{C}(\mathcal{I},a_1)$, and show that the invariant drops for $\mathscr{C}$. 

Let us first outline the set-up for induction. Let $H$ be the hypersurface of maximal contact for $\mathcal{I}$ through $y$ given by $x_1=0$, and let $\SJ_H = (x_2^{a_2},\dotsc,x_k^{a_k},(Q \subset M)^{1/d})$ be the restriction of $\SJ$ to $H$. Let $\overline{\SJ}_H$ denote the reduced toroidal center associated to $\SJ_H$, so $\overline{\SJ}_H = (\SJ_H)^{c'/(a_1n_1)}$ where $c' = \gcd(n_2,\dotsc,n_k)$. Note that $\SJ^{(y)}(\mathscr{C}|_H) = \SJ_H^{(a_1-1)!} = \overline{\SJ}_H^{(a_1!n_1)/c'}$, so $\overline{\SJ}_H$ is the reduced toroidal center associated to $\SJ^{(y)}(\mathscr{C}|_H)$. Since the length of $\inv_y(\mathscr{C}|_H)$ is $<L$, the induction hypothesis implies in particular that the invariant of $\mathscr{C}|_H$ at $y$ drops after the weighted toroidal blow-up along $\overline{\SJ}_H^{1/(cc')}$. But the weighted toroidal blow-up along $\overline{\SJ}_H^{1/(cc')}$ coincides with the the proper transform $H_c' \rightarrow H$ of $H$ via the weighted toroidal blow-up along $\overline{\SJ}^{1/c}$ (\ref{4406}(ii)).

Therefore, to leverage on the preceding paragraph, we consider the following two cases: \begin{enumerate}
    \item[(a)] $y'$ is in the $(x_1T^{n_1c})$-chart of $Y_c'$;
    \item[(b)] $y'$ is in the proper transform $H_c'$, in which case $y'$ is in the other charts of $Y_c'$.
\end{enumerate}

For case (a), the local section $x_1^{a_1!}$ of $\mathscr{C}$ factors as $x_1^{a_1!} = u^{a_1!n_1c} \cdot 1$ in $\mathscr{C}\SO_{Y_c'} = u^{a_1!n_1c} \cdot \mathscr{C}_c'$, where $u$ is the equation for $\SE_c$. Therefore, $(\mathscr{C}'_c)_{y'} = (1)$, i.e. $\inv_{y'}(\mathscr{C}'_c) = (0) < \inv_y(\mathscr{C})$, as desired. 

For case (b), we saw earlier that the induction hypothesis implies: \begin{align}\label{eq7.1}
    \inv_{y'}(\mathscr{C}_c'|_{H_c'}) < \inv_y(\mathscr{C}|_H).
\end{align}
Moreover, the local section $x_1^{a_1!}$ of $\mathscr{C}$ now factors as $x_1^{a_1!} = u^{a_1!n_1c} \cdot (x_1')^{a_1!}$ in $\mathscr{C}\SO_{Y_c'} = u^{a_1!n_1c} \cdot \mathscr{C}_c'$, where $u$ is the equation for $\SE_c$, and $x_1'$ is the equation for $H_c'$. Thus, $(x_1')^{a_1!} \subset (\mathscr{C}_c')_{y'}$, so that $\logord_{y'}(\mathscr{C}_c') \leq a_1!$. Let us now consider two sub-cases of (b): \begin{enumerate}
    \item[(bi)] If $\logord_{y'}(\mathscr{C}_c') < a_1!$, then a fortiori $\inv_{y'}(\mathscr{C}_c') < \inv_y(\mathscr{C})$.
    \item[(bii)] On the other hand, if $\logord_{y'}(\mathscr{C}_c') = a_1!$, then $x_1'$ is a maximal contact element for $\mathscr{C}_c'$ at $y'$, so $H_c'$ is a hypersurface of maximal contact for $\mathscr{C}_c'$ through $y'$. Therefore: \begin{align*}
    \inv_{y'}(\mathscr{C}_c') &= \left(a_1!,\frac{\inv_{y'}(\mathscr{C}(\mathscr{C}_c',a_1!)|_{H_c'})}{(a_1!-1)!}\right) \qquad \textrm{by (\ref{6103}(i))} \\
    &\leq \left(a_1!,\frac{\inv_{y'}(\mathscr{C}(\mathscr{C}_c'|_{H_c'},a_1!))}{(a_1!-1)!}\right) \qquad \textrm{since } \mathscr{C}(\mathscr{C}_c',a_1!)|_{H_c'} \supset \mathscr{C}(\mathscr{C}_c'|_{H_c'},a_1!) \\
    &= \big(a_1!,\inv_{y'}(\mathscr{C}_c'|_{H_c'})\big) \qquad \qquad \qquad \; \textrm{by (\ref{6204}(ii))} \\
    &< \big(a_1!,\inv_y(\mathscr{C}|_H)\big) \qquad \qquad \qquad \;\;\;\; \textrm{by (\ref{eq7.1})} \\
    &= (a_1-1)! \cdot \inv_y(\mathcal{I}) \qquad \qquad \qquad \;\; \textrm{by (\ref{6103}(i))} \\
    &= \inv_y(\mathscr{C}) \qquad \qquad \qquad \qquad \qquad \;\;\; \textrm{by (\ref{6204}(ii))}
    \end{align*}
    as desired.
\end{enumerate}
This completes the proof of the induction step.
\end{proof}

\subsection{Proof of logarithmic principalization}\label{7.3}

\begin{proof}[Proof of \emph{(\ref{7101})}]
For the first paragraph of the theorem, take $\SJ = \SJ(\mathcal{I})$ as in \S\ref{chpt6}. Following the notation in (\ref{6104}), write $\SJ = (x_1^{a_1},\dotsc,x_k^{a_k},(Q \subset M)^{1/d})$, and write $\overline{\SJ} = \SJ^{1/(a_1n_1)} = (x_1^{1/n_1},\dotsc,x_k^{1/n_k},(Q \subset M)^{1/(a_1n_1d)})$. Let $Y' = \Bl_Y(\overline{\SJ}) \rightarrow Y$ as in the theorem. By (\ref{4404}(i)), $\mathcal{I}\SO_{Y'}$ factors as $\SE^{a_1n_1} \cdot \mathcal{I}'$. By (\ref{6201}(ii)), $\mathcal{I}'$ is the weak transform of $\mathcal{I}$. By (\ref{7201}), $\max\inv(\mathcal{I}') < \max\inv(\mathcal{I})$. Functoriality with respect to logarithmically smooth surjective morphisms follows from (\ref{6401}).

The second paragraph of the theorem is now immediate by a standard argument. Namely, if $Y_n \rightarrow \dotsb \rightarrow Y$ is the logarithmic principalization of $\mathcal{I} \subsetneq \SO_Y$ and $\widetilde{Y} \rightarrow Y$ is a logarithmically smooth morphism of toroidal Deligne-Mumford stacks over $\mathds{k}$ with $\widetilde{\mathcal{I}} = \mathcal{I}\SO_{\widetilde{Y}}$, then the logarithmic principalization of $\widetilde{\mathcal{I}} \subsetneq \SO_{\widetilde{Y}}$ is obtained from the pullback of $Y_n \rightarrow \dotsb \rightarrow Y$ to $\widetilde{Y}$ by removing empty blow-ups.
\end{proof}

\subsection{Proof of logarithmic embedded resolution}\label{7.4}

\begin{proof}[Proof of \emph{(\ref{1101})}]
This proceeds in the same way as the proof of \cite[Theorem 1.1.1]{ATW19}. For the first paragraph in the theorem, one applies (\ref{7101}) to the ideal $\mathcal{I} = \mathcal{I}_X$ defining $X$ in $Y$, and replaces the weak transform $\mathcal{I}'$ by the proper transform $\mathcal{I}_{X'} \supset \mathcal{I}'$. This implies part (ii) of the theorem, i.e. $\max\inv(\mathcal{I}_{X'}) \leq \max\inv(\mathcal{I}') < \max\inv(\mathcal{I}_X)$. Parts (i) and (iv) of the theorem were noted in the paragraphs between (\ref{4402}) and (\ref{4403}), while (iii) follows from the fact that the chosen toroidal center $\SJ = \SJ(\mathcal{I})$ is $\mathcal{I}$-admissible (\ref{6201}(i)).

The second paragraph of the theorem is just a repeated application of the first paragraph. One stops at the point where $\max\inv(\mathcal{I}_{X_N})$ is the sequence $(1,\dotsc,1)$ of length $c$ (where $c$ is the codimension of $X$ in $Y$): at this point, the toroidal center $\SJ_N$, whose support is contained in $X_N$, is everywhere of the form $(x_1,\dotsc,x_c)$ for ordinary parameters $x_i$, and hence the support of $\SJ_N$ is in particular toroidal. Since $\inv_p(\mathcal{I}_{X_N}) = (1,\dotsc,1)$ at a point $p$ at which $X_N$ is logarithmically smooth, we have that the support of $\SJ_N$ contains a dense open in $X_N$, whence they coincide, and $X_N$ is toroidal. This gives (a), while (b) and (c) are immediate from (iii) and (iv).
\end{proof}

\subsection{Proof of re-embedding principle}\label{7.5}

\begin{proof}[Proof of \emph{(\ref{1202})}]
We may assume $Y$ is a strict toroidal $\kk$-scheme. Write $\mathbb{A}^1_\mathds{k} = \Spec(\mathds{k}[x_0])$, and let $\mathcal{I}_{X \subset Y}$ denote the ideal of $X$ in $Y$. Then the ideal $\mathcal{I}_{X \subset Y_1}$ of $X$ in $Y_1$ is $(x_0)+\mathcal{I}_{X \subset Y}$. Then $\SD_Y^{\leq 1}(\mathcal{I}_{X \subset Y_1})_p = (1)$ with maximal contact element $x_0$ everywhere, so that $\mathcal{I}[2] = \mathcal{I}_{X \subset Y_1}|_{V(x_0)=Y} = \mathcal{I}_{X \subset Y}$. Therefore, part (i) follows by definition of the invariant (\S\ref{6.1}).

For (ii), first note that if $\SJ(\mathcal{I}_{X \subset Y}) = (x_1^{a_1},\dotsc,x_k^{a_k},(Q \subset M)^{1/d})$, then $\SJ(\mathcal{I}_{X \subset Y_1}) = (x_0,x_1^{a_1},\dotsc,x_k^{a_k},(Q \subset M)^{1/d})$. Then the fact that $Y'$ is identified with the proper transform $V(x_0')$ of $Y = V(x_0) \subset Y_1$ in $Y_1'$ follows from (\ref{4406}(ii)). Moreover, if $\mathcal{I}_{X \subset Y}'$ (resp. $\mathcal{I}_{X \subset Y_1}'$) denotes the underlying ideal of $X' \subset Y'$ (resp. $X_1' \subset Y_1'$), then $\mathcal{I}_{X \subset Y_1}' = (x_0') + \mathcal{I}_{X \subset Y}'$, and hence part (ii) follows.
\end{proof}

\section{An example}\label{chpt8}

Consider the set-up in \S\ref{1.3}. We show, by way of example, that the toroidal Deligne-Mumford stacks $Y_i$ obtained in our logarithmic embedded resolution algorithm $Y_N \rightarrow \dotsb \rightarrow Y_1 \rightarrow Y_0 = Y$ is not necessarily smooth over $\mathds{k}$, and the proper transform $X_N = Y_N \times_Y X_N$ is not necessarily smooth over $\mathds{k}$. This necessitates the need for \emph{resolution of toroidal singularities}, as outlined in (\ref{1203}).

\subsection{Resolution of toroidal singularities is necessary}\label{8.1}

We revisit the following singular surface in \cite[8.3]{ATW19}: \begin{align*}
    X = V(\mathcal{I}) = V(x^2yz+y^4z) \subset Y = \mathbb{A}_{\mathds{k}}^3.
\end{align*}
While $Y_1$ and $Y_2$ for this example are smooth over $\mathds{k}$, we will see below that $Y_3$ is not. We do this by focusing on a particular chart at each step of our logarithmic resolution algorithm. \begin{enumerate}
    \item[(I)] Since $\SD^{\leq 4}(\mathcal{I}) = (x,y,z)$, we have $\max\inv(X \subset Y) = \inv_{(0,0,0)}(X \subset Y) = (4,4,4)$, and $\SJ(\mathcal{I}) = (x^4,y^4,z^4)$. Rescaling, the first step in our logarithmic resolution algorithm involves the blow-up $Y_1 \rightarrow Y$ along $\overline{\SJ}(\mathcal{I}) = (x,y,z)$. Here $Y_1$ is a priori a strict toroidal $\mathds{k}$-scheme, but in fact it is also smooth over $\mathds{k}$. This can be seen by examining the $x$, $y$ and $z$-charts. For example, the $z$-chart $Y_1^{(z)}$ of $Y_1$ is given by the strict toroidal $\mathds{k}$-scheme which is also smooth over $\mathds{k}$: \begin{align*}
        \qquad Y_1^{(z)} = \Spec(\NN^1 \rightarrow \mathds{k}[x_1,y_1,\underline{z}_1])
    \end{align*}
    where $x=x_1\underline{z}_1$, $y = y_1\underline{z}_1$, and $z=\underline{z}_1$ is the equation of the exceptional divisor. Here we underline $\underline{z}_1$ to indicate that it is the image of the standard basis vector $e_1$ of $\NN^1$ under the logarithmic structure $\NN^1 \rightarrow \mathds{k}[x_1,y_1,\underline{z}_1]$ (as given by (\ref{4401})). In this chart, the equation $(x^2yz+y^4z)$ of $X \subset Y$ becomes $\underline{z}_1^4(x_1^2y_1+y_1^4\underline{z}_1)$, with proper transform \begin{align*}
        \qquad X_1^{(z)} = V(\mathcal{I}_1^{(z)}) = V(x_1^2y_1+y_1^4\underline{z}_1) \subset Y_1^{(z)} = \Spec(\NN^1 \rightarrow \mathds{k}[x_1,y_1,\underline{z}_1]).
    \end{align*}
    \item[(II)] Next, we have $\SD^{\leq 1}(\mathcal{I}_1^{(z)}) = (x_1y_1,x_1^2+4y_1^3\underline{z}_1,y_1^4\underline{z}_1)$ and $\SD^{\leq 2}(\mathcal{I}_1^{(z)}) = (x_1,y_1)$, whence $\mathscr{C}(\mathcal{I}_1^{(z)},3)|_{y_1 = 0} = (x_1^6)$. Therefore, $\max\inv(X_1^{(z)} \subset Y_1^{(z)}) = (3,3) < (4,4,4) = \max\inv(X \subset Y)$, and $\SJ(\mathcal{I}_1^{(z)}) = (x_1^3,y_1^3)$. Rescaling, the second step in our logarithmic resolution algorithm for the $z$-chart involves the blow-up $Y_2^{(z)} \rightarrow Y_1^{(z)}$ along $\overline{\SJ}(\mathcal{I}_1^{(z)}) = (x_1,y_1)$. Similar to (I), $Y_2^{(z)}$ is a strict toroidal $\mathds{k}$-scheme which is smooth over $\mathds{k}$. For example, the $y_1$-chart $Y_2^{(z,y_1)}$ of $Y_2^{(z)}$ is given by the strict toroidal $\mathds{k}$-scheme which is also smooth over $\mathds{k}$: \begin{align*}
        \qquad Y_2^{(z,y_1)} = \Spec(\NN^2 \rightarrow \mathds{k}[x_2,\underline{\smash{y}}_2,\underline{z}_2])
    \end{align*}
    where $x_1 = \underline{\smash{y}}_2x_2$, $\underline{z}_1 = \underline{z}_2$, and $y_1 = \underline{\smash{y}}_2$ is the equation of the exceptional divisor. Once again we underline $\underline{\smash{y}}_2$ and $\underline{z}_2$ to indicate that they are the respective images of the standard basis vectors $e_1$ and $e_2$ of $\NN^2$ under the logarithmic structure $\NN^2 \rightarrow \mathds{k}[x_2,\underline{\smash{y}}_2,\underline{z}_2]$ on $Y_2^{(z,y_1)}$ (as given by (\ref{4401})). In this chart, the equation $(x_1^2y_1 + y_1^4\underline{z}_1)$ of $X_1^{(z)} \subset Y_1^{(z)}$ becomes $\underline{\smash{y}}_2^3(x_2^2+\underline{\smash{y}}_2\underline{z}_2)$, with proper transform \begin{align*}
        \qquad X_2^{(z,y_1)} = V(\mathcal{I}_2^{(z,y_1)}) = V(x_2^2+\underline{\smash{y}}_2\underline{z}_2) \subset Y_2^{(z,y_2)} = \Spec(\NN^2 \rightarrow \mathds{k}[x_2,\underline{\smash{y}}_2,\underline{z}_2]).
    \end{align*} 
    \item[(III)] Finally, we have $\SD^{\leq 1}(\mathcal{I}_2^{(z,y_1)}) = (x_2,\underline{\smash{y}}_2\underline{z}_2)$, whence $\max\inv(X_2^{(z,y_1)} \subset Y_2^{(z,y_1)}) = (2,\infty) < (3,3) = \max\inv(X_1^{(z)} \subset Y_1^{(z)})$. Since $\mathscr{C}(\mathcal{I}_2^{(z,y_1)},2)|_{x_2=0} = (\underline{\smash{y}}_2\underline{z}_2)$, we also have $\SJ(\mathcal{I}_2^{(z,y_1)}) = (x_2^2,\underline{\smash{y}}_2\underline{z}_2)$. Rescaling, the third step in our logarithmic resolution algorithm for the $y_1$-chart involves the blow-up $Y_3^{(z,y_1)} \rightarrow Y_2^{(z,y_1)}$ along $\overline{\SJ}(\mathcal{I}_2^{(z,y_1)}) = (x_2,(\underline{\smash{y}}_2\underline{z}_2)^{1/2})$. A priori $Y_3^{(z,y_1)}$ is a toroidal Deligne-Mumford stack over $\mathds{k}$, but this time the $x_2$-chart $Y_3^{(z,y_1,x_2)}$ of $Y_3^{(z,y_1)}$ is no longer smooth over $\mathds{k}$: \begin{align*}
        \qquad Y_3^{(z,y_1,x_2)} = \Spec\left(\frac{\NN^4}{\langle e_2+e_3 \thicksim e_4+2e_1 \rangle} \rightarrow \frac{\mathds{k}[\underline{x}_3,\underline{\smash{y}}_3,\underline{z}_3,\underline{w}_3]}{(\underline{\smash{y}}_3\underline{z}_3-\underline{w}_3\underline{x}_3^2)}\right)
    \end{align*}
    where $\underline{\smash{y}}_2 = \underline{\smash{y}}_3$, $\underline{z}_2 = \underline{z}_3$, $\underline{\smash{y}}_2\underline{z}_2 = \underline{w}_3x_2^2$, and $x_2 = \underline{x}_3$ is the equation of the exceptional divisor. As before, we underline $\underline{x}_3$, $\underline{\smash{y}}_3$, $\underline{z}_3$ and $\underline{w}_3$ to indicate that they are the respective images of the standard basis vectors $e_1$, $e_2$, $e_3$ and $e_4$ of $\NN^4$ under the logarithmic structure $\NN^4/\langle e_2+e_3 \thicksim e_4+2e_1 \rangle \rightarrow \mathds{k}[\underline{x}_3,\underline{\smash{y}}_3,\underline{z}_3,\underline{w}_3]/(\underline{\smash{y}}_3\underline{z}_3-\underline{w}_3\underline{x}_3^2)$ on $Y_3^{(z,y_1,x_2)}$ (as given by (\ref{4401})). In this chart, the equation $(x_2^2+\underline{\smash{y}}_2\underline{z}_2)$ of $X_2^{(z,y_1)} \subset Y_2^{(z,y_1)}$ becomes $\underline{x}_3^2(1+\underline{w}_3)$, with proper transform \begin{align*}
        \qquad \qquad X_3^{(z,y_1,x_2)} = V(1+\underline{w}_3) \subset \Spec\left(\frac{\NN^4}{\langle e_2+e_3 \thicksim e_4+2e_1 \rangle} \rightarrow \frac{\mathds{k}[\underline{x}_3,\underline{\smash{y}}_3,\underline{z}_3,\underline{w}_3]}{(\underline{\smash{y}}_3\underline{z}_3-\underline{w}_3\underline{x}_3^2)}\right).
    \end{align*}
    Note that $\max\inv(X_3^{(z,y_1,x_2)} \subset Y_3^{(z,y_1,x_2)}) = (1) < (2,\infty) = \max\inv(X_2^{(z,y_1)} \subset Y_2^{(z,y_1)})$, so our logarithmic embedded resolution algorithm stops here (for this chart). In other words, $X_3^{(z,y_1,x_2)}$ is toroidal. However, as a scheme, \begin{align*}
        \qquad X_3^{(z,y_1,x_2)} \simeq \Spec\left(\frac{\mathds{k}[x_3,y_3,z_3]}{(x_3^2+y_3z_3)}\right)
    \end{align*}
    is not smooth over $\mathds{k}$.
\end{enumerate}

\appendix
\section{Zariski-Riemann space}\label{appA}

In this appendix, fix an algebraic function field $K$, over a ground field $\kk$.

\subsection{Zariski-Riemann space of \texorpdfstring{$K/\kk$}{K/k}}\label{A.1}

The Zariski-Riemann space of $K/\kk$, which we shall describe shortly, was originally called the Riemann manifold of $K/\kk$ by Zariski, in his proof of resolution of singularities of $\kk$-varieties\footnote{See footnote \#2 at the start of \S 2.} of dimensions $2$ and $3$. This notion is implicit in Hironaka's work on resolution of singularities for all dimensions in characteristic zero. It also plays an essential role in \cite{ATW19}, as well as this paper. We shall describe this space in steps: \begin{enumerate}
    \item As a set, \begin{align*}
        \ZR(K,\kk) := \lbrace\textup{valuation rings $R$ of $K$ containing $\kk$} \rbrace.
    \end{align*}
    We usually denote an element $R$ of $\ZR(K,\kk)$ by its corresponding valuation $\nu \colon K^\ast \twoheadrightarrow G$ instead, where $G = \lbrace xR \colon x \in K^\ast \rbrace$ is the value group of $\nu$. In that case, we write $R_\nu$ for $R$, and $G_\nu$ for $G$. We denote the unique maximal ideal of $R_\nu$ by $\mathfrak{m}_\nu$, and its residue field by $\kappa_\nu = R_\nu/\mathfrak{m}_\nu$.
    \item As a topological space, $\ZR(K,\kk)$ has a basis of open sets given by $\mathscr{F} = \lbrace U(x_1,\dotsc,x_n) \colon$ $n \geq 0$ and $x_i \in K^\ast \rbrace$, where $U(x_1,\dotsc,x_n) = \lbrace \nu \in \ZR(K,\kk) \colon R_\nu \supset \kk[x_1,\dotsc,x_n] \rbrace$.
    \item Finally, as a locally ringed space, $\ZR(K,\kk)$ is equipped with a sheaf of rings $\SO = \SO_{\ZR(K,\kk)}$ described by \begin{align*}
        \SO(U) := \bigcap_{\nu \in U}{R_\nu} \qquad \textup{where } U \subset \ZR(K,\kk) \textup{ is open}.
    \end{align*}
    In particular, $\SO(U(x_1,\dotsc,x_n))$ is the integral closure of $\kk[x_1,\dotsc,x_n]$ in $K$ \cite[Theorem 10.4]{Mat89}. Then $\SO$ is a subsheaf of the constant sheaf $K$ on $\ZR(K,\kk)$, and the stalk of $\SO$ at $\nu$ is $R_\nu$. Note that $\ZR(K,\kk)$ also carries a sheaf of ordered groups $\Gamma = K^\ast/\SO^\ast$, whose sections over an open set $U$ are: \begin{align*}
        \textup{\footnotesize{$\Big\lbrace (s_\nu)_{\nu \in U} \in \prod_{\nu \in U}{G_\nu} \colon \forall \nu \in U$, $\exists$ open set $\nu \in V \subset U$ and $\exists x \in K^\ast$ such that $\forall \nu' \in V$, $s_{\nu'} = \nu'(x)$}} \Big\rbrace
    \end{align*}
    and whose stalk at $\nu$ is $G_\nu$, with a morphism of sheaves of ordered groups $\val \colon K^\ast \twoheadrightarrow \Gamma$. The image $\val(\SO \setminus \lbrace 0 \rbrace) \subset \Gamma$ is the sheaf of monoids consisting of non-negative sections of $\Gamma$, denoted $\Gamma_+$. Explicitly, its sections over an open set $U$ are: \begin{align*}
        \; \textup{\footnotesize{$\Big\lbrace (s_\nu)_{\nu \in U} \in \prod_{\nu \in U}{G_\nu} \colon \forall \nu \in U$, $\exists$ open set $\nu \in V \subset U$ and $\exists\ 0 \neq x \in \SO(V)$ such that $\forall \nu' \in V$, $s_{\nu'} = \nu'(x)$}} \Big\rbrace.
    \end{align*}
\end{enumerate}

Two remarks are in order. Firstly, $\ZR(K,\kk)$ is quasi-compact. For a proof, see \cite[Theorem 10.5]{Mat89}. Secondly, $\ZR(K,\kk)$ can be characterized as an inverse limit of projective models of $K/\kk$ in the category of locally ringed spaces. Let us expound on this further. By a projective model $Y$ of $K/\kk$, we mean $Y$ is a projective $\kk$-variety, whose field of functions $K(Y)$ is isomorphic to $K$. For every $\nu \in \ZR(K,\kk)$, there exists a unique dotted arrow making the triangles in the diagram below commute: \begin{equation*}
    \begin{tikzcd}
    \Spec(K) \arrow[to=1-3, "\textup{generic pt.}"] \arrow[to=2-1] & & Y \arrow[to=2-3] \\
    \Spec(R_\nu) \arrow[to=1-3, dotted, "f_\nu"] \arrow[to=2-3] & & \Spec(\kk)
    \end{tikzcd}
\end{equation*}
The composition $\Spec(\kappa_\nu) = \Spec(R_\nu/\mathfrak{m}_\nu) \rightarrow \Spec(R_\nu) \xrightarrow{f_\nu} Y$ demarcates a point $y_\nu$ on $Y$, which is called the center of $R_\nu$ on $Y$ \cite[Exercise II.4.5]{Har77}. This gives an injective local $\kk$-homomorphism $f_\nu^\# \colon \SO_{Y,y_\nu}$ $\rightarrow R_\nu$ of local rings whose field of fractions is $K$, in which case we say $R_\nu$ dominates $\SO_{Y,y_\nu}$ via $f_\nu$ (cf. \cite[Lemma II.4.4]{Har77}).

The projective models of $K/\kk$ form an inverse system as follows: an arrow from a projective model $Y_a$ to another $Y_b$ is a birational morphism $\varphi_{a \rightarrow b} \colon Y_b \rightarrow Y_a$. For every $\nu$ in $\ZR(K,\kk)$, $\varphi_{a \rightarrow b}$ necessarily maps the center $y_{b,\nu}$ of $R_\nu$ on $Y_b$ to the center $y_{a,\nu}$ of $R_\nu$ on $Y_a$. In other words, $\varphi^b_a$ induces a local homomorphism $\varphi_{a \rightarrow b}^\#$ of local rings with field of fractions $K$, which makes the diagram below commute: \begin{equation*}
    \begin{tikzcd}
    \SO_{Y_a,y_{a,\nu}} \arrow[to=1-3, "\varphi_{a \rightarrow b}^\#"] \arrow[to=2-2, swap, "f_{a,\nu}^\#"] & & \SO_{Y_b,y_{b,\nu}} \arrow[to=2-2, "f_{b,\nu}"] \\
    & R_\nu &
    \end{tikzcd}
\end{equation*}
The join $Y_c$ of two projective models $Y_a$ and $Y_b$ admits birational morphisms $Y_c \rightarrow Y_a$ and $Y_c \rightarrow Y_b$, whence this is indeed an inverse system.

\begin{proof}[Proof that $\ZR(K,\kk)$ is the set-theoretic inverse limit]
As shown above, a point $\nu \in \ZR(K,\kk)$ determines a collection of points $\lbrace y_\nu \in Y \colon Y$ is a projective model of $K/\kk \rbrace$ --- which is (by definition) preserved by arrows in the inverse system --- and hence determines a point in the inverse limit. 

Conversely, a point in the inverse limit is a collection of points $\Sigma = \lbrace y_\Sigma \in Y \colon Y$ is a projective model of $K/\kk \rbrace$ which is preserved by arrows in the inverse system. Let $R$ be the direct limit of the system whose objects are the local rings $\SO_{Y,y_\Sigma}$, and whose arrows are given by the local $\kk$-homomorphisms $\varphi_{a \rightarrow b}^\#$ of local rings with field of fractions $K$ (where $a \rightarrow b$ is an arrow in the inverse system of projective models of $K/\kk$). 

Since $R$ is the direct limit of a system of local rings with local homomorphisms, $R$ is a local ring with maximal ideal $(\mathfrak{m}_{Y,y_\Sigma} \colon Y$ is a projective model for $K/\kk)$. By \cite[Theorem I.6.1A]{Har77}, $R$ is a valuation ring of $K$ containing $\kk$, and hence determines a point $\nu \in \ZR(K,\kk)$. For each projective model $Y$ of $K/\kk$, $\SO_{Y,y_\Sigma}$ must be the unique local ring of $Y$ dominated by $R$, whence the center of $R$ on $Y$ is $y_\Sigma$. One can also use \cite[Theorem I.6.1A]{Har77} to show that given a point $\nu \in \ZR(K,\kk)$, $R_\nu$ is the direct limit of the system of local rings $\SO_{Y,y_\nu}$. This establishes the desired one-to-one correspondence of points.
\end{proof}

\begin{proof}[Proof that $\ZR(K,\kk)$ is the topological inverse limit]
The inverse limit topology on $\ZR(K,\kk)$ is the coarsest topology such that the projection maps $\pi_Y \colon \ZR(K,\kk) \rightarrow Y$ (where $Y$ is a projective model for $K/\kk$), which sends $\nu \in \ZR(K,\kk)$ to the center $y_\nu$ of $R_\nu$ on $Y$, are continuous.

Let $Y$ be a projective model for $K/\kk$. Let $U = \Spec(A) \subset Y$ be an open affine subset. Then $\pi_Y^{-1}(\Spec(A))$ consists of $\nu \in \ZR(K,\kk)$ such that there exists a dotted arrow filling in the diagram below: \begin{equation*}
    \begin{tikzcd}
    K & A \arrow[to=1-1] \arrow[to=2-1, dotted] \\
    R_\nu \arrow[to=1-1] & \kk \arrow[to=1-2] \arrow[to=2-1] 
    \end{tikzcd}
\end{equation*}
Since $A$ is a finitely generated $\kk$-algebra, we can write $A = \kk[x_1,\dotsc,x_n]$ for $x_i \in K^\ast$. Then $U(x_1,\dotsc,x_n) = \pi_Y^{-1}(\Spec(A))$. Conversely, given $x_1,\dotsc,x_n \in K^\ast$, we can find $x_{n+1},\dotsc,x_m \in K^\ast$ such that $A = \kk[x_1,\dotsc,x_m]$ has fraction field $K$. The projection $T_i \mapsto x_i$ gives a presentation of $A$ as $A \simeq \kk[T_1,\dotsc,T_m]/\mathfrak{p}$, where $\mathfrak{p}$ is a prime ideal of the polynomial ring $\kk[T_1,\dotsc,T_m]$. We can homogenize the prime ideal $\mathfrak{p} \subset \kk[T_1,\dotsc,T_m]$ to a homogeneous prime ideal $\mathfrak{P} \subset \kk[T_1,\dotsc,T_m,T_{m+1}]$. Then $U = \Spec(A)$ is an open affine subset of $Y = \Proj(\kk[T_1,\dotsc,T_{m+1}]/\mathfrak{P})$, which is a projective model of $K/\kk$ with $\pi_Y^{-1}(\Spec(A)) = U(x_1,\dotsc,x_m) \subset U(x_1,\dotsc,x_n)$. Since open affines form a basis for Zariski topology, we are done.
\end{proof}

\begin{proof}[Proof that $\ZR(K,\kk)$ is the inverse limit in the category of locally ringed spaces]
Set \[
    \SO' := \varinjlim_Y{\pi_Y^{-1}\SO_Y}
\]
where the direct limit is taken over projective models $Y$ of $K/\kk$. This is the correct sheaf of rings on $\ZR(K,\kk)$ as the inverse limit in the category of locally ringed spaces (see for example \cite[Theorem 4 and Corollary 5]{Gil11}). It remains to note that $\SO' = \SO_{\ZR(K,\kk)}$. For this, observe there are morphisms $\pi_Y^{-1}\SO_Y \rightarrow \SO_{\ZR(K,\kk)}$ (adjoint to the canonical morphisms $\SO_Y \rightarrow (\pi_Y)_\ast\SO_{\ZR(K,\kk)}$) for each projective model $Y$ of $K/\kk$, culminating in a morphism $\SO' \rightarrow \SO_{\ZR(K,\kk)}$ which is an isomorphism by checking it on stalks.
\end{proof}

Note that $\ZR(K,\kk)$ is also the inverse limit of a similar system of proper models of $K/\kk$ (proper $\kk$-varieties, whose field of fractions is isomorphic to $K$), in which the projective models of $K/\kk$ form a cofinal subsystem (by Chow's Lemma \cite[Exercise II.4.10]{Har77}). 

\subsection{Zariski-Riemann space of a \texorpdfstring{$\kk$}{k}-variety}\label{A.2}

More generally, we can define the Zariski-Riemann space for a $\kk$-variety $Y$. Let $K$ be the field of fractions $K(Y)$ of $Y$. Since $Y$ is separated but not necessarily proper, not every $\nu \in \ZR(K,\kk)$ possesses a center $y_\nu$ on $Y$, but if it does, the center $y_\nu$ is unique. Therefore, we set \begin{align*}
    \ZR(Y) := \lbrace \nu \in \ZR(K,\kk) \colon \nu \textup{ has a center on } Y \rbrace \subset \ZR(K,\kk).
\end{align*}
This agrees with the notation in \cite{ATW19}. If $Y$ is a proper model of $K/\kk$, $\ZR(Y)$ is simply $\ZR(K,\kk)$ defined in Appendix \ref{A.1}. We let $\ZR(Y)$ inherit its topology, sheaf of rings $\SO_{\ZR(Y)}$, and sheaf of ordered groups $\Gamma_Y$ from $\ZR(K,\kk)$. As before, $\ZR(Y)$ is the inverse limit of the system of modifications $Y' \rightarrow Y$ in the category of morphisms of locally ringed spaces into $Y$. We write $\pi_Y$ for the morphism $\ZR(Y) \rightarrow Y$ sending $\nu$ to the center of $\nu$ on $Y$. 

Note that $\ZR(Y)$ is quasi-compact and open in $\ZR(K,\kk)$. This can be seen as follows. First suppose $Y = \Spec(A)$ is an affine $\kk$-variety, with $A$ generated as a $\kk$-algebra by $x_1,\dotsc,x_n \in K$. In this case, we have seen earlier that $\ZR(Y)$ is $U(x_1,\dotsc,x_n)$, and is quasi-compact (by \cite[Theorem 10.5]{Mat89}). In general, since $Y$ is covered by finitely many affine opens, one deduces that $\ZR(Y)$ is quasi-compact and open in $\ZR(K,\kk)$. We conclude this section with a noteworthy fact:

\begin{lemma}\label{A201}
Let $Y$ be a $\kk$-variety, with morphism $\pi_Y \colon \ZR(Y) \rightarrow Y$. If $Y$ is normal, then the morphism $\pi_Y^\# \colon \SO_Y \rightarrow (\pi_Y)_\ast\SO_{\ZR(Y)}$ is an isomorphism of sheaves on $Y$.
\end{lemma}

\begin{proof}
Since open affines form a basis for the Zariski topology on $Y$, it suffices to check this isomorphism on open affines $U = \Spec(A) \subset Y$. Since $Y$ is normal, $\SO_Y(U) = A$ is normal, whence by \cite[Theorem 10.4]{Mat89}, \begin{align*}
    \SO_Y(U) = \bigcap_{\substack{\nu \in \ZR(K,\kk) \\ R_\nu \supseteq A}}{R_\nu}
\end{align*}
But the set of $\nu \in \ZR(K,\kk)$ such that $R_\nu \supseteq A$ is precisely the set of $\nu \in \ZR(Y)$ which has a center on $U = \Spec(A) \subset Y$. Therefore, $\SO_Y(U) = \bigcap_{\nu \in \pi_Y^{-1}(U)}{R_\nu} = \SO_{\ZR(K,\kk)}(\pi_Y^{-1}(U))$.
\end{proof}

\subsection{Functoriality with respect to dominant morphisms}\label{A.3}

If $f \colon Y' \rightarrow Y$ is a dominant morphism of $\kk$-varieties, $f$ induces a morphism $\ZR(f) \colon \ZR(Y') \rightarrow \ZR(Y)$ of locally ringed spaces, which maps $R_\nu \mapsto R_\nu \cap K(Y)$. The morphism $\SO_{\ZR(Y)} \rightarrow \ZR(f)_\ast\SO_{\ZR(Y')}$ is given by the inclusion $\bigcap_{\nu \in U}{R_\nu} \hookrightarrow \bigcap_{\eta \in \ZR(f)^{-1}(U)}{R_\eta}$ over an open set $U$, and is stalk-wise given by the local homomorphism $R_\nu \cap K(Y) \hookrightarrow R_\nu$. This morphism $\SO_{\ZR(Y)} \rightarrow \ZR(f)_\ast\SO_{\ZR(Y')}$ descends to a morphism of sheaves of ordered groups $\Gamma_Y \rightarrow \ZR(f)_\ast\Gamma_{Y'}$, as well as a morphism of sheaves of monoids $\Gamma_{Y,+} \rightarrow \ZR(f)_\ast\Gamma_{Y',+}$.

\section{Toroidal geometry}\label{appB}

In this appendix, we briefly mention some preliminaries on logarithmic geometry pertinent to this paper. Most of the notation and language here follows \cite{Ogu18} closely. Other relevant references include \cite{Kat89}, \cite{Kat94}, \cite{Niz06}, \cite{AT17}, \cite{ATW20a}, \cite{ATW20b}, and \cite{GR19}. 

\subsection{Toroidal \texorpdfstring{$\kk$}{k}-schemes}\label{B.1}

In this section, $Y$ denotes a logarithmic scheme, and we denote its underlying scheme by $\underline{Y}$, and its underlying logarithmic structure by $\alpha_Y \colon \SM_Y \rightarrow \SO_Y$. Occasionally we also use the letter $\underline{Y}$ to denote the logarithmic scheme given by the scheme $\underline{Y}$ equipped with the trivial logarithmic structure.

\begin{definition}\label{B101}
We say $Y$ is \emph{fs}, if: \begin{quote}
    $Y$ admits a covering $\mathcal{U}$ (in the Zariski or \'etale topology, depending if $\SM_Y$ is Zariski or not) such that the pullback of $\SM_Y$ to each $U$ in $\mathcal{U}$ admits a chart subordinate to a fs (= \emph{fine} and \emph{saturated}) monoid $M$ --- or equivalently, $U$ admits a strict morphism $U \rightarrow \Spec(M \rightarrow \ZZ[M])$ for a fs monoid $M$.
\end{quote}
\end{definition}

\begin{xpar}\label{B102}
Let $Y$ be a fs logarithmic scheme. In what follows, $y$ always denotes a point in $Y$, while $\overline{y}$ denotes a geometric point over $y$. Then one can show: \begin{enumerate}
    \item $\overline{\SM}_{Y,\overline{y}}^{\gp}$ is a free abelian group of finite rank $r(y)$. Note $r(y)$ is independent of choice of $\overline{y}$ over $y$. See \cite[Proposition I.1.3.5(2)]{Ogu18}.
    \item $r(y) = \rk(\overline{\SM}_{Y,\overline{y}}^{\gp})$ is upper semi-continuous on $Y$, i.e. for each $n \in \NN$, \begin{align*}
        Y^{(n)} := \lbrace y \in Y \colon \rk(\overline{\SM}_{Y,\overline{y}}^{\gp}) \leq n \rbrace
    \end{align*}
    is Zariski open in $Y$ \cite[Corollary II.2.16]{Ogu18}. In particular, $Y^\ast := \lbrace y \in Y \colon \SM_{Y,\overline{y}} := \SO_{Y,\overline{y}}^\ast \rbrace$ is a Zariski open in $Y$, called the \emph{locus of triviality} of $Y$.
    \item For each $n \in \NN$, \begin{align*}
        Y(n) := \lbrace y \in Y \colon \rk(\overline{\SM}_{Y,\overline{y}}) = n \rbrace \subset Y^{(n)}
    \end{align*}
    is a Zariski closed subscheme of $Y^{(n)}$, and has the following \'etale-local description: for all $\overline{y} \in Y(n)$, $\SO_{Y(n),\overline{y}} = \SO_{Y,\overline{y}}/I(\overline{y})$, where $I(\overline{y})$ is the ideal of $\SO_{Y,\overline{y}}$ generated by the image of the unique maximal ideal $\SM_{Y,\overline{y}}^+$ of $\SM_{Y,\overline{y}}$ under $\SM_{Y,\overline{y}} \xrightarrow{\alpha_{Y,\overline{y}}} \SO_{Y,\overline{y}}$. See \cite[2.2.10]{AT17}.
    \item After replacing $Y$ by an \'etale neighbourhood of $\overline{y}$, $Y$ admits a fine chart $M \rightarrow H^0(Y,\SM_Y)$ which is neat at $\overline{y}$, i.e. the composition $M \rightarrow H^0(Y,\SM_Y) \rightarrow \SM_{Y,\overline{y}} \rightarrow \overline{\SM}_{Y,\overline{y}}$ is an isomorphism. See \cite[Proposition III.1.2.7]{Ogu18}. In particular, \'etale locally every logarithmic scheme $Y$ is a Zariski logarithmic scheme.
\end{enumerate}
If $\SM_Y$ is Zariski, all statements apply with $\overline{y}$ replaced by the scheme-theoretic point $y \in Y$, and (iv) holds after replacing $Y$ by a Zariski neighbourhood of $y$.
\end{xpar}

\begin{definition}[Logarithmic stratification]\label{B103}
Let $Y$ be a fs logarithmic scheme. The \emph{logarithmic stratification} of $Y$ is the stratification given by $\lbrace Y(n) \colon n \in \NN \rbrace$ in (\ref{B102}(iii)). For each $y \in Y$, we set $\fs_y = Y(n)$ for $n=\rk(\overline{\SM}_{Y,\overline{y}}^{\gp})$, and $\fs_y$ is called the logarithmic stratum through $y$. 
\end{definition}

\begin{definition}\label{B104}
We say that a fs logarithmic scheme $Y$ is \emph{logarithmically regular} at a point $y \in Y$, if for some (and hence any) geometric point $\overline{y}$ over $y$, \begin{center}
    $\fs_y$ is regular at $\overline{y}$ and the equality $\dim(\SO_{Y,\overline{y}}) = \rk(\overline{\SM}_{Y,\overline{y}}^{\gp}) + \dim(\SO_{\fs_y,\overline{y}})$ holds.
\end{center}
If $Y$ is a fs Zariski logarithmic scheme, we say $Y$ is logarithmically regular at $y \in Y$, if the same statement holds with $\overline{y}$ replaced by the scheme-theoretic point $y$ throughout. We say $Y$ is logarithmically regular if $Y$ is logarithmically regular at every point $y \in Y$.
\end{definition}

\begin{xpar}\label{B105}
Let $Y$ be a fs logarithmic scheme. \begin{enumerate}
    \item In general, for every $y \in Y$, $\dim(\SO_{Y,\overline{y}}) \leq \rk(\overline{\SM}_{Y,\overline{y}}^{\gp}) + \dim(\SO_{\fs_y,\overline{y}})$ \cite[Lemma 2.3]{Kat94}. 
    \item Let $U = Y^\ast$ be the triviality locus of $Y$, with open embedding into $Y$ denoted by $j$. If $Y$ is logarithmically regular, then $\alpha_Y \colon \SM_Y \rightarrow \SO_Y$ is injective, and the image of $\alpha_Y$ is $j_\ast(\SO_U^\ast) \cap \SO_Y$. If $D = Y \setminus U$ is nonempty, then $D$ is a divisor on $Y$, called the toroidal divisor of $Y$. See \cite[Theorem 3.2.4]{Kat94} and \cite[Proposition 2.6]{Niz06}.
\end{enumerate}
If $\SM_Y$ is Zariski, then the above statements hold with $\overline{y}$ replaced by the scheme-theoretic points $y \in Y$. In addition:
\begin{enumerate} 
    \item[(iii)] If $Y$ is logarithmically regular, $\underline{Y}$ is Cohen-Macaulay and normal \cite[Theorem 4.1]{Kat94}. In particular, $\underline{Y}$ is reduced, and if $Y$ is locally Noetherian, $\underline{Y}$ is a disjoint union of its irreducible components. Moreover, $\underline{Y}$ is catenary, so each non-empty logarithmic stratum $Y(n)$ of $Y$ has pure codimension $n$. 
    \item[(iv)] It is a fact that if $Y$ is logarithmically regular at all closed points in $Y$, then $Y$ is logarithmically regular \cite[Proposition 7.1]{Kat94}.
\end{enumerate}
\end{xpar}

We collate the aforementioned properties in the following definition:

\begin{definition}[Toroidal $\kk$-schemes \textrm{\cite[2.3.4]{AT17}}]\label{B106}
Let $\kk$ be a field of characteristic zero. A \emph{toroidal $\kk$-scheme} is a fs logarithmic $\kk$-scheme $Y$ which is logarithmically regular, such that $\underline{Y}$ is of finite type over $\kk$. If moreover $\SM_Y$ is a Zariski logarithmic structure, then we say $Y$ is a strict toroidal $\kk$-scheme.
\end{definition}

Note that every regular $\kk$-scheme is a toroidal $\kk$-scheme when we equip it with the trivial logarithmic structure. The remark in \cite[Remark 2.3.5]{AT17} deserves mention here: if $\kk = \overline{\kk}$, strict toroidal $\kk$-varieties correspond to the toroidal embeddings without self-intersections in \cite{KKMSD73}. More generally, toroidal $\kk$-varieties correspond to general toroidal embeddings, possibly with self-intersections.

\begin{xpar}\label{B107}
Let $Y$ be a strict toroidal $\kk$-scheme. \begin{enumerate}
    \item For every $y \in Y$, fix $x_1,\dotsc,x_n \in \SO_{Y,y}$ which reduce to a regular system of parameters $x_1,\dotsc,x_n$ of $\SO_{\fs_y,y}$, fix a local fs chart $\beta \colon M \to H^0(U,\SM_Y|_U)$ at $y$ which is neat at $y$, and fix a coefficient field $\kappa$ for $\widehat{\SO}_{Y,y}$. Then the induced surjective homomorphism \begin{align*}
        \kappa\llbracket X_1,\dotsc,X_n,M=\overline{\SM}_{Y,y} \rrbracket \to \widehat{\SO}_{Y,y}, \qquad X_i \mapsto x_i
    \end{align*}
    is an isomorphism \cite[Theorem 3.2(1)]{Kat94}.
    \item Endow $\Spec(\kk)$ with the trivial logarithmic structure. Then $Y$ is logarithmically smooth over $\kk$. Moreover, if $\widetilde{Y}$ is a fs Zariski logarithmic $\kk$-scheme which admits a logarithmically smooth morphism $f \colon \widetilde{Y} \rightarrow Y$ to a strict toroidal $\kk$-scheme $Y$, then $\widetilde{Y}$ is also a strict toroidal $\kk$-scheme. See \cite[Proposition 8.3]{Kat94}.
\end{enumerate}
\end{xpar}

\'Etale locally every toroidal $\kk$-scheme is a strict toroidal $\kk$-scheme (\ref{B102}(iv)). Therefore, if we want to understand the \'etale-local structure of toroidal $\kk$-schemes, it suffices to explicate the local structure of strict toroidal $\kk$-schemes. We shall do this via a choice of:

\begin{definition}
[Logarithmic coordinates and parameters \textrm{\cite[3.1.2]{ATW20a}}]\label{B108}
Let $Y$ be a strict toroidal $\kk$-scheme, and let $y \in Y$. Set $n = \codim_{\fs_y}\overline{\lbrace y \rbrace}$, $N = \dim(\fs_y)$ and $M = \overline{\SM}_{Y,y}$. By a system of \emph{logarithmic coordinates} at $y$, we mean the following data: \begin{enumerate}
    \item sections $x_1,\dotsc,x_N$ of $\SO_{Y,y}$ whose images under $\SO_{Y,y} \twoheadrightarrow \SO_{\fs_y,y} \xrightarrow{d} \Omega^1_{\underline{\fs_y},y}$ reduce to a $\kappa(y)$-basis for $\Omega^1_{\underline{\fs_y}}(y)$, and such that the images of the first $n$ sections $x_1,\dotsc,x_n$ in $\SO_{\fs_y,y}$ form a regular system of parameters of $\SO_{\fs_y,y}$,
    \item and a local fs chart $\beta \colon M \rightarrow H^0(U,\SM_Y|_U)$ at $y$, which is neat at $y$.
\end{enumerate}
We usually denote this data by \[
\left((x_1,\dotsc,x_{N}),\ M=\overline{\SM}_{Y,y} \xrightarrow{\beta} H^0(U,\SM_Y|_U)\right).
\]
We call $\lbrace x_1,\dotsc,x_{N} \rbrace$ a system of \emph{ordinary coordinates} at $y$, and we call the subset $\lbrace x_1,\dotsc,x_{n} \rbrace$ a system of \emph{ordinary parameters} at $y$. The sub-data \[
    \left((x_1,\dotsc,x_n),M=\overline{\SM}_{Y,y} \xrightarrow{\beta} H^0(U,\SM_Y|_U)\right)
\]
is called a system of logarithmic parameters at $y$.

The elements of $\alpha_Y(\beta(M \setminus \lbrace 0 \rbrace))$ are called \emph{monomial parameters} at $y$. For an element $m \in M \setminus \lbrace 0 \rbrace$, we usually use the same letter $m$ for $\beta(m)$, and write $\exp(m)$ for the monomial parameter $\alpha_Y(\beta(m))$. If $(d,D)$ denotes the universal logarithmic derivation $\SO_Y \oplus \SM_Y \rightarrow \Omega^1_Y$, this notation should remind you of the ``exponential rule'' in calculus: $d(\exp(m)) = \exp(m) \cdot Dm$.
\end{definition}

In what follows, we denote the logarithmic tangent sheaf of a strict toroidal $\kk$-scheme $Y$ over $\kk$ by $\SD^1_Y$ (instead of the usual $T^1_Y$ or $T_{Y/\kk}^1$).

\begin{lemma}\label{B109}
Let $Y$ be a strict toroidal $\kk$-scheme. Let $y \in Y$, and fix a system of logarithmic parameters $\left((x_1,\dotsc,x_{N}),\ M=\overline{\SM}_{Y,y} \xrightarrow{\beta} H^0(U,\SM_Y|_U)\right)$ at $y$. Then: \begin{enumerate}
    \item The universal logarithmic derivation $(d,D) \colon \SO_Y \oplus \SM_Y \to \Omega^1_Y$ induces a natural isomorphism $\Omega^1_{Y,y} \xleftarrow{\simeq} \left(\bigoplus_{i=1}^{N}{\SO_{Y,y} \cdot dx_i}\right) \oplus \left(\SO_{Y,y} \otimes \overline{\SM}_{Y,y}^{\gp}\right)$. In particular, $\Omega^1_{Y,y}$ is generated as a $\SO_{Y,y}$-module by $dx_i$ for $1 \leq i \leq N$, as well as $D(M)$.
    
    \item For every element $L$ of $\Hom\left(\overline{\SM}_{Y,y}^{\gp},\SO_{Y,y}\right)$, there exists a unique derivation $(\mathbb{D}_L,L) \in \SD^1_{Y,y}$ such that $\mathbb{D}_L(\exp(m)) = \exp(m) \cdot L(m)$ for every monomial parameter $\exp(m)$, and $\mathbb{D}_L(x_i) = 0$ for every $1 \leq i \leq N$. This defines an isomorphism $\SD^1_{Y,y} \xleftarrow{\simeq} \left(\bigoplus_{i=1}^{N}{\SO_{Y,y} \cdot \frac{\partial}{\partial x_i}}\right) \oplus \Hom\left(\overline{\SM}_{Y,y}^{\gp},\SO_{Y,y}\right)$, where $\frac{\partial}{\partial x_i}$ is the derivation dual to $x_i$.
    
    \item Fix a basis $m_1,\dotsc,m_r \in M$ for $M^{\gp} = \overline{\SM}_{Y,y}^{\gp}$, and write $u_i = \exp(m_i)$ for $1 \leq i \leq r$. Then $\Omega^1_{Y,y}$ is a free $\SO_{Y,y}$-module with basis $dx_1,\dotsc,dx_N,\frac{du_1}{u_1},\dotsc,\frac{du_r}{u_r}$, and $\SD^1_{Y,y}$ is a free $\SO_{Y,y}$-module with dual basis $\frac{\partial}{\partial x_1},\dotsc,\frac{\partial}{\partial x_{N}},u_1\frac{\partial}{\partial u_1},\dotsc,u_r\frac{\partial}{\partial u_r}$.
\end{enumerate}
\end{lemma}

\begin{proof}[Sketch of proof]
Let $\mathbb{A}_M = \Spec(M \rightarrow \kk[M])$. Adapting the diagram in the proof of \cite[Theorem IV.3.3.3]{Ogu18}, one can deduce the following split short exact sequence \begin{align*}
    0 \rightarrow \kappa(y) \otimes \overline{\SM}_{Y,y}^{\gp} \simeq \Omega^1_{\mathbb{A}_M/\kk}(y) \rightarrow \Omega^1_{Y/\kk}(y) \rightarrow \Omega^1_{Y/\underline{Y}}(y) \simeq \Omega^1_{\underline{\fs_y}}(y) \rightarrow 0,
\end{align*}
from which (i) follows by Nakayama's Lemma. Part (ii) is the dual of (i), and part (iii) follows from (i) and (ii). (An alternative proof can be found in \cite[Lemma 3.34]{ATW20a}.)
\end{proof} 

We can now explicate the local structure of strict toroidal $\kk$-schemes:

\begin{theorem}\label{B110}
Let $Y$ be a strict toroidal $\kk$-scheme. Fix $y \in Y$, and set $M = \overline{\SM}_{Y,y}$. Then the following statements hold: \begin{enumerate}
    \item After replacing $Y$ with a Zariski neighbourhood of $y$, $Y$ admits a strict morphism $f \colon Y \rightarrow \Spec(M \rightarrow \kk[M])$.
    \item After replacing $Y$ by a Zariski neighbourhood of $y$, $f$ admits a factorization \begin{align*}
        U \xrightarrow{f_1} \Spec(M \rightarrow \kk[M \oplus \NN^n]) \xrightarrow{g_1} \Spec(M \rightarrow \kk[M])
    \end{align*}
    where $n = \codim_{\fs_y}\overline{\lbrace y \rbrace}$, $f_1$ is strict and smooth of relative dimension $\dim\overline{\lbrace y \rbrace}$, $f_1$ maps $y$ to the vertex of $\Spec(M \rightarrow \kk[M \oplus \NN^n])$, and $g_1$ is induced by the inclusion $M \hookrightarrow M \oplus \NN^n$.
    \item After replacing $Y$ by a Zariski neighbourhood of $y$, $f_1$ admits a factorization \begin{align*}
        U \xrightarrow{f_2} \Spec(M \rightarrow \kk[M \oplus \NN^N]) \xrightarrow{g_2} \Spec(M \rightarrow \kk[M \oplus \NN^n])
    \end{align*}
    where $N = \dim(\fs_y)$, $f_2$ is strict and \'etale, and $g_2$ is induced by the inclusion $\NN^n \hookrightarrow \NN^N$ into the first $n$ coordinates.
\end{enumerate}
\end{theorem}

\begin{proof}[Sketch of proof]
Part (i) follows from (\ref{B102}(iv)). The remaining parts follow from (\ref{B109}) and \cite[Theorem IV.3.2.3(2) and Proposition IV.3.16]{Ogu18}.
\end{proof}

The remainder of this section reviews some notions developed in \cite[\S 3]{ATW20a}, which are pertinent to this paper.

\begin{definition}[Logarithmic differential operators]\label{B111}Let $Y$ be a strict toroidal $\kk$-scheme. \begin{enumerate}
    \item For each natural number $n \geq 1$, let $\SD_Y^{\leq n}$ be the $\SO_Y$-submodule of the total sheaf $\SD_{\underline{Y}}^\infty$ of differential operators on $\underline{Y}$ generated by $\SO_Y$ and the images of $(\SD^1_Y)^{\otimes i}$ for $1 \leq i \leq n$. $\SD_Y^{\leq n}$ is called the sheaf of \emph{logarithmic differential operators on $Y$ of order $\leq n$}.
    \item The direct limit $\bigcup_{n \in \NN}{\SD_Y^{(\leq n)}} \subset \SD_{\underline{Y}}^\infty$ is called the \emph{total sheaf of logarithmic differential operators} of $Y$, and is denoted by $\SD^\infty_Y$. 
    \item Given an ideal $\mathcal{I}$ on $Y$, let $\SD_Y^{\leq n}(\mathcal{I})$ (resp. $\SD_Y^\infty(\mathcal{I})$) denote the ideal on $Y$ generated by the image of $\mathcal{I}$ under $\SD_Y^{\leq n}$ (resp. under $\SD_Y^\infty$). 
\end{enumerate}

When $Y$ is clear from context, we usually write $\SD_Y^{\leq n}$ as $\SD^{\leq n}$ (likewise for $\SD_Y^\infty$). We caution the reader that the definition in (\ref{B111}) only makes sense for $\characteristic(\kk) = 0$.
\end{definition}

\begin{definition}[Monomial ideals and saturation]\label{B112}
Let $Y$ be a strict toroidal $\kk$-scheme, and let $\mathcal{I}$ be an ideal on $Y$. \begin{enumerate}
    \item We say that $\mathcal{I}$ is a \emph{monomial ideal} if is generated by the image of an ideal $\mathcal{Q} \subset \SM_Y$ under $\alpha_Y \colon \SM_Y \rightarrow \SO_Y$.
    \item The \emph{monomial saturation} of $\mathcal{I}$, denoted $\SM(\mathcal{I})$, is defined to be the intersection of the collection of all monomial ideals on $Y$ containing $\mathcal{I}$.
\end{enumerate}
\end{definition}

Evidently, $\SM(\mathcal{I})$ contains $\mathcal{I}$, and $\mathcal{I}$ is monomial if and only if $\mathcal{I} = \SM(\mathcal{I})$. 

\begin{lemma}\label{B113}
Let $Y$ be a strict toroidal $\kk$-scheme. The following statements hold for an ideal $\mathcal{I}$ on $Y$: \begin{enumerate}
    \item $\mathcal{I}$ is monomial if and only if $\SD^{\leq 1}_Y(\mathcal{I}) = \mathcal{I}$.
    \item $\SD^{\infty}_Y(\mathcal{I}) = \SM(\mathcal{I})$.
    \item If $f \colon \widetilde{Y} \rightarrow Y$ is a logarithmically smooth morphism of strict toroidal $\kk$-schemes, then $\SD^{\leq n}_{\widetilde{Y}}(\mathcal{I}\SO_{\widetilde{Y}}) = \SD^{\leq n}_Y(\mathcal{I})\SO_{\widetilde{Y}}$ for all natural numbers $n \geq 1$, and $\SM(\mathcal{I}\SO_{\widetilde{Y}}) = \SM(\mathcal{I})\SO_{\widetilde{Y}}$.
    \item If $\SQ$ is a monomial ideal on $Y$, then $\SD_Y^{\leq n}(\SQ \cdot \mathcal{I}) = \SQ \cdot \SD_Y^{\leq n}(\mathcal{I})$ for all natural numbers $n \geq 1$.
\end{enumerate}
\end{lemma}

\begin{proof}
This is \cite[Corollary 3.3.12, Theorem 3.4.2 and Lemma 3.5.2]{ATW20a}.
\end{proof}

\begin{definition}[Logarithmic order]\label{B114}
Let $Y$ be a strict toroidal $\kk$-scheme. If $\mathcal{I}$ is an ideal on $Y$, the \emph{logarithmic order} of $\mathcal{I}$ at a point $y \in Y$ is defined as \begin{align*}
    \logord_y(\mathcal{I}) = \ord_y(\mathcal{I}|_{\fs_y}) \in \NN \cup \lbrace \infty \rbrace,
\end{align*}
where $\ord_y$ refers to the usual order of an ideal at a point (see for example, in \cite[Definition 3.47]{Kol07}). The maximal logarithmic order of $\mathcal{I}$ is $\max\logord(\mathcal{I}) = \max_{y \in Y}{\logord_y(\mathcal{I})}$.
\end{definition}

\begin{lemma}\label{B115}
Let $Y$ be a strict toroidal $\kk$-scheme. The following statements hold for an ideal $\mathcal{I}$ on $Y$, and a point $y \in Y$: \begin{enumerate}
    \item $\logord_y(\mathcal{I}) = \min\lbrace n \in \NN \colon \SD_Y^{\leq n}(\mathcal{I})_y = \SO_{Y,y} \rbrace$, where we take $\min(\emptyset) = \infty$ by convention.
    \item  $\logord_y(\mathcal{I}) = \infty$ if and only if $y \in V(\SM(\mathcal{I}))$.
    \item $\SM(\mathcal{I}) = (1)$ if and only if $\max\logord(\mathcal{I}) < \infty$.
    \item If $f \colon \widetilde{Y} \rightarrow Y$ is a logarithmically smooth morphism of strict toroidal $\kk$-schemes, and $\widetilde{y} \in \widetilde{Y}$ maps to $y \in Y$, then $\logord_{\widetilde{y}}(\mathcal{I}\SO_{\widetilde{Y}}) = \logord_y(\mathcal{I})$.
\end{enumerate}
\end{lemma}

Note that parts (i) and (ii) say that $\logord_y(\mathcal{I})$ is upper semi-continuous on $Y$: (i) for a natural number $n$, $V(\SD^{\leq n}_Y(\mathcal{I}))$ is the locus of points $y \in Y$ satisfying $\logord_y(\mathcal{I}) > n$; (ii) $V(\SM(\mathcal{I}))$ is the locus of points $y \in Y$ satisfying $\logord_y(\mathcal{I}) = \infty$.

\begin{proof}
This is \cite[Lemma 3.6.3, Lemma 3.6.5, Corollary 3.66 and Lemma 3.6.8]{ATW20a}.
\end{proof}

\subsection{Toroidal Deligne-Mumford stacks over \texorpdfstring{$\kk$}{k}}\label{B.2}

Let $\kk$ be a field of characteristic zero. Before defining the notion of a toroidal Deligne-Mumford stack over $\kk$, we recall some preliminaries from \cite[3.3]{ATW20b}. A logarithmic structure $\SM_Y$ on a Deligne-Mumford stack $Y$ is a sheaf of monoids on the \'etale site $Y_{\textup{\'et}}$, and a homomorphism $\alpha_Y \colon \SM_Y \rightarrow \SO_{Y_{\textup{\'et}}}$ inducing an isomorphism $\SM_Y^\ast \xrightarrow{\simeq} \SO_{Y_{\textup{\'et}}}^\ast$. The pair $(Y,\SM_Y)$ is then called a logarithmic Deligne-Mumford stack. If $p_{1.2} \colon Y_1 \rightrightarrows Y_0$ is an atlas of $Y$ by schemes, then a logarithmic structure $\SM_Y$ on $Y$ is equivalent to logarithmic structures $\SM_{Y_i}$ on $Y_i$ (for $i=0,1$) such that $p_1^\ast\SM_{Y_0} = \SM_{Y_1} = p_2^\ast\SM_{Y_0}$. We say a logarithmic Deligne-Mumford stack is \emph{fs}, if for some (and hence any) atlas $p_{1,2} \colon Y_1 \rightrightarrows Y_0$ of $Y$ by schemes, $(Y_0,\SM_{Y_0})$ is fs.

\begin{definition}[Toroidal DM stacks \textup{\cite[3.3.3]{ATW20b}}]\label{B201}
A \emph{toroidal Deligne-Mumford stack over $\kk$} is a fs logarithmic Deligne-Mumford stack $(Y,\SM_Y)$ over $\kk$ admitting an atlas $p_{1,2} \colon Y_1 \rightrightarrows Y_0$ by schemes such that $(Y_0,\SM_{Y_0})$ is a toroidal $\kk$-scheme.
\end{definition}

If $Y$ is a toroidal Deligne-Mumford stack over $\kk$, then $(Y_0,\SM_{Y_0})$ is a toroidal $\kk$-scheme for every atlas $p_{1,2} \colon Y_1 \rightrightarrows Y_0$ of $Y$ by schemes. This follows from \cite[Proposition 12.5.46]{GR19}. Moreover, since \'etale locally every toroidal $\kk$-scheme is a strict toroidal $\kk$-scheme, we may choose the atlas in (\ref{B201}) such that $(Y_0,\SM_{Y_0})$ is a strict toroidal $\kk$-scheme. In this case $Y_1$ is also a strict toroidal $\kk$-scheme.

\section{Proof of \texorpdfstring{(\ref{5101})}{(5.1)}}\label{appC}

The proof of the theorem follows ideas from both \cite[Lemma 5.3.3]{ATW20a} and \cite[Theorem 3.92]{Kol07}. In particular, we need to make a modification to \cite[Proposition 3.94]{Kol07}. Let us first fix some notation: let $\kk$ be a field of characteristic zero, $\kappa/\kk$ be a field extension, $M$ be a sharp monoid (written multiplicatively), and consider the logarithmic $\kk$-algebra $M \rightarrow \kappa\llbracket \NN^n \oplus M \rrbracket = \kappa\llbracket x_1,\dotsc,x_n,M \rrbracket = R$, with maximal ideal $\mathfrak{m} = (x_1,\dotsc,x_n,M \setminus \lbrace 1 \rbrace)$. For a proper ideal $J \subset \mathfrak{m}$ of $R$, we say an automorphism $\psi$ of $R$ is of the form $\mathds{1}+J$, if $\psi$ maps each $x_i$ to $x_i+f_i$ for some $f_i \in J$, and fixes $M$. For an ideal $I \subset R$, we have \begin{align*}
    \SD^{\leq 1}(I) = I+\left(\frac{\partial f}{\partial x_i} \colon f \in I,\ 1 \leq i \leq n\right)
\end{align*}
and inductively, we have $\SD^{\leq \ell}(I) = \SD(\SD^{\leq \ell-1}(I))$ for all $\ell \geq 2$.

\begin{lemma}\label{C01}
Let the notation be as above, and let $I \subset R$ be an ideal. The following statements are equivalent: \begin{enumerate}
    \item $\psi(I)=I$ for every automorphism $\psi$ of the form $\mathds{1}+J$.
    \item $J \cdot \SD^{\leq 1}(I) \subset I$.
    \item $J^\ell \cdot \SD^{\leq \ell}(I) \subset I$ for every $\ell \geq 1$.
\end{enumerate}
\end{lemma}

\begin{proof}[Proof of \emph{(\ref{C01})}]
This proof proceeds in the same way as \cite[Proposition 3.94]{Kol07}, with minor modifications.

Assume (iii). Let $\psi$ be an automorphism of the form $\mathds{1}+J$, and for all $1 \leq i \leq n$, let $b_i \in J$ such that $\psi(x_i) = x_i+b_i$. Then Taylor expansion gives us \begin{align*}
    \psi(f) = f+\sum_{i=1}^n{b_i\frac{\partial f}{\partial x_i}}+\frac{1}{2}\sum_{i,j=1}^n{b_ib_j\frac{\partial^2 f}{\partial x_i \partial x_j}} + \dotsb.
\end{align*}
For any $\ell \geq 1$, we get \begin{align*}
    \psi(f) \in I + J \cdot \SD^{\leq 1}(I) + \dotsb + J^\ell \cdot \SD^{\leq \ell}(I) + \mathfrak{m}^{\ell+1} \subset I+\mathfrak{m}^{\ell+1}.
\end{align*}
By Krull's Intersection Theorem, this implies $\psi(f) \in I$, so we get (i).

Next, assume (i). Let $b \in J$, and let $1 \leq i \leq n$. For general $\lambda \in \kk$, the endomorphism on $R$, which maps $(x_1,\dotsc,x_n) \mapsto (x_1,\dotsc,x_{i-1},x_i+\lambda b,x_{i+1},\dotsc,x_n)$ and fixes $M$, is an automorphism of $R$ of the form $\mathds{1}+J$. Therefore, for every $f \in I$, and every $\ell \geq 1$, \begin{align*}
    \left(f+\lambda b \frac{\partial f}{\partial x_i}+\dotsb+(\lambda b)^\ell\frac{\partial^\ell f}{\partial x_i^\ell}\right) \in \psi(f) + \mathfrak{m}^{\ell+1} \subset I + \mathfrak{m}^{\ell+1}.
\end{align*}
For $\ell+1$ general elements $\lambda  = \lambda_0,\dotsc,\lambda_\ell$ in $\kk$, the column vector obtained from \begin{align*}
    \begin{pmatrix}1&\lambda_0&\lambda_0^2&\dotsb&\lambda_0^\ell\\
    1&\lambda_1&\lambda_1^2&\dotsb&\lambda_1^\ell \\
    \vdots&\vdots&\vdots&\ddots&\vdots \\
    1&\lambda_\ell&\lambda_\ell^2&\dotsb&\lambda_\ell^\ell\end{pmatrix} \cdot \begin{pmatrix}f\\ b \frac{\partial f}{\partial x_i}\\ \vdots \\ b^\ell\frac{\partial^\ell f}{\partial x_i^\ell}\end{pmatrix}
\end{align*} 
has entries in $I+\mathfrak{m}^{\ell+1}$, and the Vandermonde determinant $(\lambda_i^j)$ is invertible. Therefore, $b \cdot \frac{\partial f}{\partial x_i} \in I+\mathfrak{m}^{\ell+1}$. By Krull's intersection theorem again, $b \cdot \frac{\partial f}{\partial x_i} \in I$. Since $J \cdot \SD^{\leq 1}(I)$ is generated by elements of the form $b \cdot f$ or $b \cdot \frac{\partial f}{\partial x_i}$ for $b \in J$, $f \in I$ and $1 \leq i \leq n$, this proves (ii).

Finally, assume (ii). We prove by induction that $J^\ell \cdot \SD^{\leq \ell}(I) \subset I$ for every $\ell \geq 1$. The ideal $J^{\ell+1} \cdot \SD^{\leq \ell+1}(I)$ is generated by elements of the form $b_0 \dotsm b_\ell \cdot \SD^{\leq 1}(g)$ for $g \in \SD^{\leq \ell}(I)$. The product rule says: \begin{align*}
    b_0 \dotsm b_\ell \cdot \SD^{\leq 1}(g) &= b_0 \cdot \SD^{\leq 1}(b_1 \dotsm b_\ell \cdot g) - \sum_{i=1}^\ell{\SD^{\leq 1}(b_i) \cdot (b_0 \dotsm \widehat{b_i} \dotsm b_\ell \cdot g)} \\
    &\in J \cdot \SD^{\leq 1}(J^\ell \cdot \SD^{\leq \ell}(I)) + J^\ell \cdot \SD^{\leq \ell}(I) \subset J \cdot \SD^{\leq 1}(I) + J^\ell \cdot \SD^{\leq \ell}(I) \subset I,
\end{align*}
where the last two inclusions hold by the induction hypothesis. This proves (iii).
\end{proof}

\begin{proof}[Proof of \emph{(\ref{5101})}]
Let $n = \codim_{\fs_y}\overline{\lbrace y \rbrace}$ and $M = \overline{\SM}_{Y,y}$ as in (\ref{B108}). There exist $x_2,\dotsc,x_n \in \SO_{\fs_y,y}$ such that both $x,x_2,\dotsc,x_n$ and $x',x_2,\dotsc,x_n$ form a regular system of parameters of $\SO_{\fs_y,y}$. By (\ref{B107}(ii)), we have \begin{align*}
    \kappa\llbracket x,x_2\dotsc,x_n,M \rrbracket \simeq \widehat{\SO}_{Y,y} \simeq \kappa\llbracket x',x_2,\dotsc,x_n,M \rrbracket, \quad \textup{where } \kappa = \kappa(y).
\end{align*}
Consider the endomorphism $\psi$ of $\widehat{\SO}_{Y,y}$, which maps $(x,x_2,\dotsc,x_n) \mapsto (x'=x+(x'-x),x_2,\dotsc,x_n)$ and fixes $M$. Since $x',x_2,\dotsc,x_n$ are linearly independent modulo $\mathfrak{m}_{Y,y}^2$ (where $\mathfrak{m}_{Y,y}$ is the maximal ideal of $\SO_{Y,y}$), $\psi$ is an automorphism of $\widehat{\SO}_{Y,y}$. Moreover, since $x$ and $x'$ are maximal contact elements at $y$, we have $x'-x \in \widehat{\MC(\mathcal{I})} = \MC(\widehat{\mathcal{I}})$ (note logarithmic derivatives commute with completions), whence $\psi$ is an automorphism of $\widehat{\SO}_{Y,y}$ of the form $\mathds{1}+\MC(\widehat{\mathcal{I}})$. Finally, since $\mathcal{I}$ is MC-invariant, we have $\MC(\widehat{\mathcal{I}}) \cdot \SD^{\leq 1}(\widehat{\mathcal{I}}) \subset \widehat{\mathcal{I}}$, whence (\ref{C01}) implies $\psi(\widehat{\mathcal{I}}) = \widehat{\mathcal{I}}$.

Our goal now is to realize this automorphism $\psi$ on $\widehat{\SO}_{Y,y}$ on some strict, \'etale neighbourhood $\widetilde{U}$ of $y$. We first extend both $(x,x_2,\dotsc,x_n)$ and $(x',x_2,\dotsc,x_n)$ to systems of logarithmic coordinates at $y$ (\ref{B108}): \begin{align*}
    \big((x,x_2\dotsc,x_N),M \xrightarrow{\beta} H^0(U,\SM_Y|_U)\big) \quad \textup{and} \quad \big((x',x_2\dotsc,x_N),M \xrightarrow{\beta} H^0(U,\SM_Y|_U)\big),
\end{align*} 
where $N = \dim(\fs_y)$. We then apply (\ref{B110}): after shrinking $U$ if necessary, $U$ admits strict and \'etale morphisms \begin{align*}
    U \xrightrightarrows[\tau_{x'}]{\tau_x} \Spec(M \rightarrow \kk[X_1,\dotsc,X_N,M])
\end{align*}
induced by: \begin{enumerate}
    \item[(a)] morphisms $U \rightrightarrows \mathbb{A}^N_\kk$, induced by ring morphisms $\kk[X_1,\dotsc,X_n] \rightrightarrows \Gamma(U,\SO_U)$ mapping $(X_1,X_2,\dotsc,X_N) \mapsto (x,x_2,\dotsc,x_N)$ and $(X_1,X_2,\dotsc,X_N) \mapsto (x',x_2,\dotsc,x_N)$ respectively;
    \item[(b)] as well as the chart $M = \overline{\SM}_{Y,y} \xrightarrow{\beta} H^0(U,\SM_Y|_U)$.
\end{enumerate}
Finally, we obtain $\widetilde{U}$ in the statement of (\ref{5101}), by forming the following cartesian square (in the category of fs logarithmic schemes): \begin{equation*}
    \begin{tikzcd}
    \widetilde{U} \arrow[to=1-2, "\phi_{x'}"] \arrow[to=2-1, "\phi_x"] & U \arrow[to=2-2, "\tau_{x'}"] \\
    U \arrow[to=2-2, "\tau_x"] & \Spec(M \rightarrow \kk[X_1,\dotsc,X_N,M])
    \end{tikzcd}
\end{equation*}
Since both $\tau_x$ and $\tau_{x'}$ are strict and \'etale, $\phi_x$ and $\phi_{x'}$ are also strict and \'etale. Moreover, $\phi_x^\ast(x) = \phi_x^\ast(\tau_x^\ast(X_1)) = \phi_{x'}^\ast(\tau_{x'}^\ast(X_1)) = \phi_{x'}^\ast(x')$. Note that $\tau_x$ and $\tau_{x'}$ maps $y$ to the same point in $\Spec(M \rightarrow \kk[X_1,\dotsc,X_N,M])$, so there is a point $\widetilde{y} = (y,y) \in \widetilde{U}$ which is mapped to $y$ via $\phi_x$ and $\phi_{x'}$. Finally, the completion of $\widetilde{U}$ at $\widetilde{y} = (y,y)$ is the graph of the automorphism $\psi$ on $\widehat{\SO}_{Y,y}$, and since $\psi(\widehat{\mathcal{I}}) = \widehat{\mathcal{I}}$, it follows (after shrinking $\widetilde{U}$ if necessary) that $\phi_x^\ast(\mathcal{I}) = \phi_{x'}^\ast(\mathcal{I})$.
\end{proof}

\bibliography{bibfile}

\providecommand{\MR}{\relax\ifhmode\unskip\space\fi MR }
\providecommand{\MRhref}[2]{%
  \href{http://www.ams.org/mathscinet-getitem?mr=#1}{#2}
}
\providecommand{\href}[2]{#2}
\begin{thebibliography}{KKMSD73}

\bibitem[AT17]{AT17}
Dan Abramovich and Michael Temkin, \emph{Torification of diagonalizable group
  actions on toroidal schemes}, J. Algebra \textbf{472} (2017), 279--338.

\bibitem[ATW19]{ATW19}
Dan Abramovich, Michael Temkin, and Jaros{\l}aw W{\l}odarczyk, \emph{Functorial
  embedded resolution via weighted blowings up}, arXiv e-prints, version 3:
  \url{https://arxiv.org/abs/1906.07106v3} (2019).

\bibitem[ATW20a]{ATW20a}
Dan Abramovich, Michael Temkin, and Jaros{\l}aw W{\l}odarczyk,
  \emph{Principalization of ideals on toroidal orbifolds}, JEMS, Electronically
  published, \url{https://doi.org/10.4171/JEMS/997} (2020).

\bibitem[ATW20b]{ATW20b}
Dan Abramovich, Michael Temkin, and Jaros{\l}aw W{\l}odarczyk, \emph{Toroidal
  orbifolds, destackification, and \textit{Kummer} blowings up}, Algebra and
  Number Theory \textbf{14} (2020), no.~8, 2001--2035.

\bibitem[BM08]{BM08}
Edward Bierstone and Pierre~D. Milman, \emph{Functoriality in resolution of
  singularities}, Publ. Res. Inst. Math. Sci. \textbf{44} (2008), no.~2,
  609--639.

\bibitem[BR19]{BR19}
Daniel Bergh and David Rydh, \emph{Functorial destackification and weak
  factorization of orbifolds}, arXiv e-prints:
  \url{https://arxiv.org/abs/1905.00872} (2019).

\bibitem[EV00]{EV00}
Santiago Encinas and Orlando Villamayor, \emph{A course on constructive
  desingularization and equivariance}, Resolution of singularities (Obergurgl,
  1997). Progr. Math. Birkh\"auser, Basel \textbf{181} (2000), 147--227.

\bibitem[EV07]{EV07}
Santiago Encinas and Orlando Villamayor, \emph{Rees algebras and resolution of
  singularities}, Proceedings of the XVI-th Latin American Algebra Colloquium
  (Spanish), Bibl. Rev. Mat. Iberoamericana, Rev. Mat. Iberoamericana, Madrid.
  (2007), 63--85.

\bibitem[Gil11]{Gil11}
William~Daniel Gillam, \emph{Localization of ringed spaces}, Adv. Pure Math.
  1(5) (2011), 250--263.

\bibitem[GR19]{GR19}
O.~Gabber and L.~Ramero, \emph{Foundations for almost ring theory}, arXiv
  e-prints: \url{https://arxiv.org/abs/math/0409584} (April 2019).

\bibitem[Har77]{Har77}
Robin Hartshorne, \emph{Algebraic geometry}, no.~52, Springer-Verlag, New York,
  Graduate Texts in Mathematics, 1977.

\bibitem[Hir64]{Hir64}
Heisuke Hironaka, \emph{Resolution of singularities of an algebraic variety
  over a field of characteristic zero, \textup{I \& II}}, Ann. of Math.
  \textbf{79} (1964), no.~2, 109--326.

\bibitem[Hir77]{Hir77}
Heisuke Hironaka, \emph{Idealistic exponents of singularity}, Algebraic
  geometry (J. J. Sylvester Sympos., Johns Hopkins Univ., Baltimore, Md.,
  1976), Johns Hopkins Univ. Press, Baltimore, Md. MR 0498562 (1977), 52–125.

\bibitem[Kat89]{Kat89}
Kazuya Kato, \emph{Logarithmic structures of \textup{Fontaine-Illusie}},
  Algebraic analysis, geometry, and number theory (Baltimore, MD, 1988), Johns
  Hopkins Univ. Press, Baltimore, MD, 1989. MR 1463703 (99b:14020) (1989),
  191--224.

\bibitem[Kat94]{Kat94}
Kazuya Kato, \emph{Toric singularities}, Amer. J. Math \textbf{116} (1994),
  no.~5, 1073--1099.

\bibitem[KKMSD73]{KKMSD73}
George Kempf, Finn~Faye Knudsen, David Mumford, and Bernard Saint-Donat,
  \emph{Toroidal embeddings 1}, vol. 339, Lecture Notes in Mathematics,
  Springer-Verlag, Berlin, 1973.

\bibitem[Kol07]{Kol07}
J\'anos Koll\'ar, \emph{Lectures on resolution of singularities}, vol. 166,
  Annals of Mathematics Studies, Princeton University Press, Princeton, NJ,
  2007.

\bibitem[Laz04]{Laz04}
Robert Lazarsfeld, \emph{Positivity in algebraic geometry \textup{II}:
  positivity for vector bundles, and multiplier ideals}, Springer-Verlag Berlin
  Heidelberg, 2004.

\bibitem[Mat89]{Mat89}
Hideyuki Matsumura, \emph{Commutative ring theory}, vol.~8, Cambridge Studies
  in Advanced Mathematics, Cambridge University Press, Cambridge, 1989.

\bibitem[Niz06]{Niz06}
Wies{\l}awa Nizio{\l}, \emph{Toric singularities: log-blow-ups and global
  resolutions}, J. Algebraic Geom. \textbf{15} (2006), no.~1, 1--29.

\bibitem[Ogu18]{Ogu18}
Arthur Ogus, \emph{Lectures on logarithmic algebraic geometry}, Cambridge
  studies in advanced mathematics, 178, 2018.

\bibitem[Ols16]{Ols16}
Martin Olsson, \emph{Algebraic spaces and stacks}, vol.~62, American
  Mathematical Society Colloquium Publications, American Mathematical Society,
  Providence, RI, 2016.

\bibitem[QR22]{QR22}
Ming~Hao Quek and David Rydh, \emph{Weighted blow-ups}, In preparation (2022).

\bibitem[Ryd13]{Ryd13}
David Rydh, \emph{Functorial resolution of singularities in characteristic zero
  using \textit{Rees} algebras (draft)},
  \url{https://people.kth.se/~dary/notes/ResSing_Lecture-nr45.pdf} (2013).

\bibitem[W{\l}o05]{Wlo05}
Jaros{\l}aw W{\l}odarczyk, \emph{Simple \textup{Hironaka} resolution in
  characteristic zero}, J. Amer. Math. Soc. \textbf{18} (2005), no.~4,
  779--822.

\bibitem[W{\l}o20]{Wlo19}
Jaros{\l}aw W{\l}odarczyk, \emph{Functorial resolution except of log smooth
  locus. \textit{Toroidal compactification}}, arXiv e-prints:
  \url{https://arxiv.org/abs/2007.13846} (2020).

\end{thebibliography}
\bibliographystyle{dary}

\end{document}